\documentclass[secthm,seceqn,amsthm,ussrhead,reqno, 12pt]{amsart}
\usepackage[utf8]{inputenc}
\usepackage[english]{babel}
\usepackage[symbol]{footmisc}
\usepackage{amssymb,amsmath,amsthm,amsfonts,xcolor,enumerate,hyperref,comment,longtable,cleveref}

\usepackage{times}
\usepackage{cite}
\usepackage{pdflscape}
\usepackage{ulem}
\usepackage[mathcal]{euscript}
\usepackage{tikz}
\usepackage{hyperref}
\usepackage{cancel}
\usepackage{stmaryrd}

\usetikzlibrary{arrows}

\newtheorem{theorem}{Theorem}
\newtheorem{lemma}[theorem]{Lemma}

\newtheorem{definition}[theorem]{Definition}

\newtheorem*{theoremA}{Theorem A}

\usepackage{stmaryrd}
\usepackage{xcolor}

\begin{document}

 \medskip

\noindent{\Large
The algebraic classification of nilpotent
bicommutative algebras}\footnote{
The authors thank
  Thiago Castilho de Mello \& 
  Ivan Gonzales Gargate  for their active stimulation to write this paper.
The first part of this work is supported by 
FCT   UIDB/MAT/00212/2020 and UIDP/MAT/00212/2020.
The second part of this work is supported by the Russian Science Foundation under grant 22-11-00081. 
}

 \medskip

 \medskip

\begin{center}

 {\bf
 Kobiljon Abdurasulov\footnote{ Institute of Mathematics Academy of
Sciences of Uzbekistan, Tashkent, Uzbekistan; \ abdurasulov0505@mail.ru},
Ivan Kaygorodov\footnote{CMA-UBI, Universidade da Beira Interior, Covilh\~{a}, Portugal; 
Moscow Center for Fundamental and Applied Mathematics, Moscow,   Russia; 
Saint Petersburg State University, Russia;\    kaygorodov.ivan@gmail.com} 
\& 
Abror Khudoyberdiyev\footnote{ Institute of Mathematics Academy of
Sciences of Uzbekistan, Tashkent, Uzbekistan; National University of Uzbekistan, Tashkent, Uzbekistan; \ khabror@mail.ru}

}

 \medskip

\end{center}

\

\noindent{\bf Abstract}:
{\it This paper is devoted to the complete algebraic  classification of complex  $5$-dimensional nilpotent bicommutative  algebras.}

\

 \medskip


\

 \medskip


 \medskip
 
\section*{Introduction}
The problem of classification of $n$-dimensional associative algebras was posed by Latyshev in Dniester Notebook, and some achievements were made by Pikhtilkov  and  it has also been discoursed in the paper of Belov \cite{kb}. One of the classical problems in the theory of non-associative algebras is to classify (up to isomorphism) the algebras of dimension $n$ from a certain variety defined by some family of polynomial identities. It is typical to focus on small dimensions, and there are two main directions for the classification: algebraic and geometric. Varieties such as Jordan, Lie, Leibniz or Zinbiel algebras have been studied from these two approaches (\!\cite{cfk18,klp20,  fkk222,fkk223,  fkk22,akk23, fkkv22  } and 
\cite{   akk23, ikp20,   fkkv22}, respectively).
In the present paper, we give the algebraic   classification of
$5$-dimensional nilpotent bicommutative   algebras.

One-sided commutative algebras first appeared in the paper by Cayley \cite{cayley} in 1857. 
The variety of bicommutative algebras is defined by the following identities of right- and left-commutativity:
\[
\begin{array}{rclllrcl}
(xy)z &=& (xz)y, & \ &  x(yz) &=& y(xz).
\end{array} \]
It contains the commutative associative algebras as a subvariety; 
the square of each bicommutative algebra gives a structure of a commutative associative algebra \cite{dt03};
and each bicommutative algebra is Lie admissible (in \cite{dt03,DKS09} there were shown that any  bicommutative algebra  under commutator multiplication  gives a metabelian Lie algebra).
The variety of $2$-dimensional bicommutative algebras is described by Kaygorodov and Volkov;
algebraic and geometric classification of $4$-dimensional nilpotent bicommutative algebras is given by Kaygorodov, 
 Páez-Guillán and  Voronin in \cite{kpv20};
 algebraic classification of one-generated $6$-dimensional nilpotent bicommutative algebras is given by Kaygorodov, 
 Páez-Guillán and  Voronin in \cite{kpv21}.
Bicommutative central extensions of $n$-dimensional restricted polynomial algebras are studied by Kaygorodov,
 Lopes and P\'{a}ez-Guill\'{a}n in \cite{klp20}.
 The structure of the free bicommutative algebra of countable rank and its main numerical invariants were described by
Dzhumadildaev, Ismailov and Tulenbaev \cite{dit11}, see also the announcement \cite{dt03}.
They also proved that  the bicommutative operad is not Koszul \cite{dit11}.
Shestakov and Zhang described     automorphisms of finitely generated relatively free bicommutative algebras \cite{shestak}.
Drensky and  Zhakhayev  proved that every free bicommutative algebra generated by one element is not noetherian, in the sense that it doesn't have finitely generated one-sided ideals and they also
obtained a positive solution of the Specht problem for any variety of bicommutative algebras  over an arbitrary field of any characteristic \cite{drensky1}.
Identities of $2$-dimensional bicommutative algebras
and invariant theory of free bicommutative algebras 
are studied by Drensky in \cite{drensky2,drensky}. 
Dzhumadildaev and Ismailov  prove that every identity satisfied by the commutator multiplication  in all bicommutative algebras is a consequence of anti-commutativity, the Jacobi and the metabelian identities \cite{di}.
They also proved  that in the anti-commutator case
every identity satisfied by the anti-commutator product in all bicommutative algebras is a consequence of
commutativity and two identities obtained in  \cite{di}.
Bai,  Chen and Zhang proved that the Gelfand-Kirillov dimension of a finitely generated bicommutative algebra is a nonnegative integer \cite{bcz21}.    
Bicommutative algebras are also known under the name of LR-algebras in a series of papers of Burde,  Dekimpe and their co-authors \cite{DKS09, DKV10, bdd}.
The studied structures of LR-algebras on a certain Lie algebra.
Burde,  Dekimpe and  Deschamps proved
the existence of an LR-complete structure on a nilpotent Lie algebra of dimension $n$ is equivalent to the existence of an $n$-dimensional abelian subgroup of the affine group ${\rm Aff}(N)$ which acts simply and transitively on $N$, where $N$ is the connected and simply connected Lie group associated with $n$  \cite{DKS09}.
Burde,  Dekimpe and Vercammen  show that if a nilpotent Lie algebra admits an {\rm LR}-structure, then it admits a complete {\rm LR}-structure, i.e., the right multiplication for the {\rm LR}-structure is always nilpotent. Extending this result, it is proven that a meta-solvable Lie algebra with two generators also admits a complete {\rm LR}-structure \cite{DKV10}.

Our method for classifying nilpotent  bicommutative  algebras is based on the calculation of central extensions of nilpotent algebras of smaller dimensions from the same variety. 
The algebraic study of central extensions of   algebras has been an important topic for years \cite{  klp20,hac16,  ss78}.
First, Skjelbred and Sund used central extensions of Lie algebras to obtain a classification of nilpotent Lie algebras  \cite{ss78}.
Note that the Skjelbred-Sund method of central extensions is an important tool in the classification of nilpotent algebras.
Using the same method,  
 small dimensional nilpotent 
(associative, 
 terminal, Jordan,
  Lie, 
 anticommutative) algebras,
and some others have been described. Our main results related to the algebraic classification of the  variety of bicommutative algebras are summarized below.

\begin{theoremA}%\label{teoA}
Up to isomorphism, there are infinitely many isomorphism classes of  
complex  non-split non-one-generated $5$-dimensional   nilpotent (non-2-step nilpotent) 
non-commutative bicommutative  algebras, 
described explicitly  in  section \ref{secteoA} in terms of 
$77$ one-parameter families, 
$20$ two-parameter families, 
$3$ three-parameter families and 
$107$ additional isomorphism classes.

\end{theoremA}

\newpage
\section{The algebraic classification of nilpotent bicommutative algebras}

\subsection{Method of classification of nilpotent algebras}

The objective of this section is to give an analogue of the Skjelbred-Sund method for classifying nilpotent bicommutative algebras. As other analogues of this method were carefully explained in, for example, \cite{hac16,kpv20}, we will give only some important definitions, and refer the interested reader to the previous sources. %We will also employ their notations.

Let $({\bf A}, \cdot)$ be a bicommutative  algebra of dimension $n$ over  $\mathbb C$
and $\mathbb V$ a vector space of dimension $s$ over ${\mathbb C}$. We define the $\mathbb C$-linear space ${\rm Z}^2\left(
\bf A,\mathbb V \right) $  as the set of all  bilinear maps $\theta  \colon {\bf A} \times {\bf A} \longrightarrow {\mathbb V}$
such that
\begin{center}
$\theta(xy,z)=\theta(xz,y)$ and $\theta(x,yz)= \theta(y,xz).$ 
\end{center}
These maps will be called {\it  cocycles}. Consider a
linear map $f$ from $\bf A$ to  $\mathbb V$, and set $\delta f\colon {\bf A} \times
{\bf A} \longrightarrow {\mathbb V}$ with $\delta f  (x,y ) =f(xy )$. Then, $\delta f$ is a cocycle, and we define ${\rm B}^2\left(
{\bf A},{\mathbb V}\right) =\left\{ \theta =\delta f\ : f\in \textup{Hom}\left( {\bf A},{\mathbb V}\right) \right\}$, which is a linear subspace of ${\rm Z}^2\left( {\bf A},{\mathbb V}\right)$. Its elements are called
{\it  coboundaries}. The {\it  second cohomology space} ${\rm H}^2\left( {\bf A},{\mathbb V}\right) $ is defined to be the quotient space ${\rm Z}^2
\left( {\bf A},{\mathbb V}\right) \big/{\rm B}^2\left( {\bf A},{\mathbb V}\right) $.
% The equivalence class of $%
%\theta \in {\rm Z}^2\left( {\bf A},{\mathbb V}\right) $ will be denoted by $\left[
%\theta \right] \in {\rm H}^2\left( {\bf A},{\mathbb V}\right) $.

Let ${\rm Aut}({\bf A}) $ be the automorphism group of the bicommutative algebra ${\bf A} $ and let $\phi \in {\rm Aut}({\bf A})$. Every $\theta \in
{\rm Z}^2\left( {\bf A},{\mathbb V}\right) $ defines $\phi \theta (x,y)
=\theta \left( \phi \left( x\right) ,\phi \left( y\right) \right) $, with $\phi \theta \in {\rm Z}^2\left( {\bf A},{\mathbb V}\right) $. It is easily checked that ${\rm Aut}({\bf A})$
acts on the right on ${\rm Z}^2\left( {\bf A},{\mathbb V}\right) $, and that
 ${\rm B}^2\left( {\bf A},{\mathbb V}\right) $ is invariant under the action of ${\rm Aut}({\bf A}).$  
 So, we have that ${\rm Aut}({\bf A})$ acts on ${\rm H}^2\left( {\bf A},{\mathbb V}\right)$.

Let $\theta$ be a cocycle, and consider the direct sum ${\bf A}_{\theta } = {\bf A}\oplus {\mathbb V}$ with the
bilinear product `` $\left[ -,-\right] _{{\bf A}_{\theta }}$'' defined by $\left[ x+x^{\prime },y+y^{\prime }\right] _{{\bf A}_{\theta }}=
 xy +\theta(x,y) $ for all $x,y\in {\bf A},x^{\prime },y^{\prime }\in {\mathbb V}$.
It is straightforward that ${\bf A_{\theta}}$ is a bicommutative algebra if and only if $\theta \in {\rm Z}^2({\bf A}, {\mathbb V})$; it is  then a $s$-dimensional central extension of ${\bf A}$ by ${\mathbb V}$.

We also call the
set ${\rm Ann}(\theta)=\left\{ x\in {\bf A}:\theta \left( x, {\bf A} \right)+ \theta \left({\bf A} ,x\right) =0\right\} $
the {\it  annihilator} of $\theta $. We recall that the {\it  annihilator} of an  algebra ${\bf A}$ is defined as
the ideal ${\rm Ann}(  {\bf A} ) =\left\{ x\in {\bf A}:  x{\bf A}+ {\bf A}x =0\right\}$. Observe
 that
${\rm Ann}\left( {\bf A}_{\theta }\right) =\big({\rm Ann}(\theta) \cap{\rm Ann}({\bf A})\big)
 \oplus {\mathbb V}$.

\begin{definition}
Let ${\bf A}$ be an algebra and $I$ be a subspace of ${\rm Ann}({\bf A})$. If ${\bf A}={\bf A}_0 \oplus I$ as a direct sum of ideals,
then $I$ is called an {\it annihilator component} of ${\bf A}$.
\end{definition}
\begin{definition}
A central extension of an algebra $\bf A$ without annihilator component is called a {\it non-split central extension}.
\end{definition}

The following result is fundamental for the classification method.

\begin{lemma}\label{l1}
Let ${\bf A}$ be an $n$-dimensional bicommutative algebra such that $\dim({\rm Ann}({\bf A}))=s\neq0$. Then there exists, up to isomorphism, a unique $(n-s)$-dimensional bicommutative  algebra ${\bf A}'$ and a bilinear map $\theta \in {\rm Z}^2({\bf A}, {\mathbb V})$ with ${\rm Ann}({\bf A})\cap{\rm Ann}(\theta)=0$, where $\mathbb V$ is a vector space of dimension $s$, such that ${\bf A} \cong {{\bf A}'}_{\theta}$ and
 ${\bf A}/{\rm Ann}({\bf A})\cong {\bf A}'$.
\end{lemma}

For the proof, we refer the reader to~\cite[Lemma 5]{hac16}.

Then, in order to decide when two bicommutative algebras with nonzero annihilator are isomorphic, it suffices to find conditions in terms of the cocycles.

Let us fix a basis $\{e_{1},\ldots ,e_{s}\}$ of ${\mathbb V}$, and $
\theta \in {\rm Z}^2\left( {\bf A},{\mathbb V}\right) $. Then $\theta $ can be uniquely
written as $\theta \left( x,y\right) =
\displaystyle \sum_{i=1}^{s} \theta _{i}\left( x,y\right) e_{i}$, where $\theta _{i}\in
{\rm Z}^2\left( {\bf A},\mathbb C\right) $. It holds that $\theta \in
{\rm B}^2\left( {\bf A},{\mathbb V}\right)$ if and only if all $\theta _{i}\in {\rm B}^2\left( {\bf A},
\mathbb C\right) $, and it also holds that ${\rm Ann}(\theta)={\rm Ann}(\theta _{1})\cap\ldots \cap{\rm Ann}(\theta _{s})$. 
Furthermore, if ${\rm Ann}(\theta)\cap {\rm Ann}\left( {\bf A}\right) =0$, then ${\bf A}_{\theta }$ has an
annihilator component if and only if $\left[ \theta _{1}\right],\ldots ,\left[ \theta _{s}\right] $ are linearly
dependent in ${\rm H}^2\left( {\bf A},\mathbb C\right)$ (see \cite[Lemma 13]{hac16}).

Recall that, given a finite-dimensional vector space ${\mathbb V}$ over $\mathbb C$, the {\it  Grassmannian} $G_{k}\left( {\mathbb V}\right) $ is the set of all $k$-dimensional
linear subspaces of $ {\mathbb V}$. Let $G_{s}\left( {\rm H}^2\left( {\bf A},\mathbb C\right) \right) $ be the Grassmannian of subspaces of dimension $s$ in
${\rm H}^2\left( {\bf A},\mathbb C\right) $.
 For $W=\left\langle
\left[ \theta _{1}\right],\dots,\left[ \theta _{s}
\right] \right\rangle \in G_{s}\left( {\rm H}^2\left( {\bf A},\mathbb C
\right) \right) $ and $\phi \in {\rm Aut}({\bf A})$, define $\phi W=\left\langle \left[ \phi \theta _{1}\right],\dots,\left[ \phi \theta _{s}\right]
\right\rangle $. It holds that $\phi W\in G_{s}\left( {\rm H}^2\left( {\bf A},\mathbb C \right) \right) $, and this induces an action of ${\rm Aut}({\bf A})$ on $G_{s}\left( {\rm H}^2\left( {\bf A},\mathbb C\right) \right) $. We denote the orbit of $W\in G_{s}\left(
{\rm H}^2\left( {\bf A},\mathbb C\right) \right) $ under this action  by ${\rm Orb}(W)$. Let
\[
W_{1}=\left\langle \left[ \theta _{1}\right],\dots,
\left[ \theta _{s}\right] \right\rangle ,W_{2}=\left\langle \left[ \vartheta
_{1}\right],\dots,\left[ \vartheta _{s}\right]
\right\rangle \in G_{s}\left( {\rm H}^2\left( {\bf A},\mathbb C\right)
\right).
\]
Similarly to~\cite[Lemma 15]{hac16}, in case that $W_{1}=W_{2}$, it holds that \[ \bigcap\limits_{i=1}^{s}{\rm Ann}(\theta _{i})\cap {\rm Ann}\left( {\bf A}\right) = \bigcap\limits_{i=1}^{s}
{\rm Ann}(\vartheta _{i})\cap{\rm Ann}( {\bf A}) ,\] 
and therefore the set
\[
T_{s}({\bf A}) =\left\{ W=\left\langle \left[ \theta _{1}\right],\dots,\left[ \theta _{s}\right] \right\rangle \in
G_{s}\left( {\rm H}^2\left( {\bf A},\mathbb C\right) \right) : \bigcap\limits_{i=1}^{s}{\rm Ann}(\theta _{i})\cap{\rm Ann}({\bf A}) =0\right\}
\]
is well defined, and it is also stable under the action of ${\rm Aut}({\bf A})$ (see~\cite[Lemma 16]{hac16}).

Now, let ${\mathbb V}$ be an $s$-dimensional linear space and let us denote by
$E\left( {\bf A},{\mathbb V}\right) $ the set of all non-split $s$-dimensional central extensions of ${\bf A}$ by
${\mathbb V}$. We can write
\[
E\left( {\bf A},{\mathbb V}\right) =\left\{ {\bf A}_{\theta }:\theta \left( x,y\right) = \sum_{i=1}^{s}\theta _{i}\left( x,y\right) e_{i} \ \ \text{and} \ \ \left\langle \left[ \theta _{1}\right],\dots,
\left[ \theta _{s}\right] \right\rangle \in T_{s}({\bf A}) \right\} .
\]

Having established these results, we can determine whether two $s$-dimensional non-split central extensions ${\bf A}_{\theta },{\bf A}_{\vartheta }$ are isomorphic or not. For the proof, see~\cite[Lemma 17]{hac16}.

\begin{lemma}
 Let ${\bf A}_{\theta },{\bf A}_{\vartheta }\in E\left( {\bf A},{\mathbb V}\right) $. Suppose that $\theta \left( x,y\right) =  \displaystyle \sum_{i=1}^{s}
\theta _{i}\left( x,y\right) e_{i}$ and $\vartheta \left( x,y\right) =
\displaystyle \sum_{i=1}^{s} \vartheta _{i}\left( x,y\right) e_{i}$.
Then the bicommutative algebras ${\bf A}_{\theta }$ and ${\bf A}_{\vartheta } $ are isomorphic
if and only if
$${\rm Orb}\left\langle \left[ \theta _{1}\right],\dots,\left[ \theta _{s}\right] \right\rangle =
{\rm Orb}\left\langle \left[ \vartheta _{1}\right],\dots,\left[ \vartheta _{s}\right] \right\rangle .$$
\end{lemma}

Then, it exists a bijective correspondence between the set of ${\rm Aut}({\bf A})$-orbits on $T_{s}\left( {\bf A}\right) $ and the set of
isomorphism classes of $E\left( {\bf A},{\mathbb V}\right) $. Consequently we have a
procedure that allows us, given a bicommutative algebra ${\bf A'}$ of
dimension $n-s$, to construct all its non-split central extensions.

\; \;

{\centerline {\textsl{Procedure}}}

Let ${\bf A}'$ be a bicommutative algebra of dimension $n-s $.

\begin{enumerate}
\item Determine ${\rm H}^2( {\bf A}',\mathbb {C}) $, ${\rm Ann}({\bf A}')$ and ${\rm Aut}({\bf A}')$.

\item Determine the set of ${\rm Aut}({\bf A}')$-orbits on $T_{s}({\bf A}') $.

\item For each orbit, construct the bicommutative algebra associated with a
representative of it.
\end{enumerate}

It follows that, thanks to this procedure and to Lemma~\ref{l1}, we can classify all the nilpotent bicommutative algebras of dimension $n$, provided that the nilpotent bicommutative algebras of dimension $n-1$ are known.

\subsection{Notations}
Let ${\bf A}$ be a bicommutative algebra and fix
a basis $\{e_{1},\dots,e_{n}\}$. We define the bilinear form
$\Delta _{ij} \colon {\bf A}\times {\bf A}\longrightarrow \mathbb C$
by $\Delta _{ij}\left( e_{l},e_{m}\right) = \delta_{il}\delta_{jm}$.
Then the set $\left\{ \Delta_{ij}:1\leq i, j\leq n\right\} $ is a basis for the linear space of
the bilinear forms on ${\bf A}$, and in particular, every $\theta \in
{\rm Z}^2\left( {\bf A},\mathbb V\right) $ can be uniquely written as $
\theta = \displaystyle \sum_{1\leq i,j\leq n} c_{ij}\Delta _{{i}{j}}$, where $
c_{ij}\in \mathbb C$.
${\rm H}_{com}^2({\mathfrak N})$ is the subspace of commutative cocycles of 
${\rm H}_{bicom}^2({\mathfrak N}),$ 
where ${\rm H}_{bicom}^2({\mathfrak N})$ is the cohomology space  for bicommutative cocycles of algebra ${\mathfrak N}$.
Let us fix the following notations:

$$\begin{array}{lll}
%{\mathcal B}^{i*}_j& \mbox{---}& j\mbox{th }i\mbox{-dimensional nilpotent {\it non-pure} bicommutative algebra (with identity $xyz=0$}); \\
%{\mathcal B}^i_j& \mbox{---}& j\mbox{th }i\mbox{-dimensional nilpotent {\it pure} bicommutative algebra (without identity $xyz=0$)}; \\
%{\mathfrak{N}}_i& \mbox{---}& i\mbox{-dimensional algebra with zero product}; \\
%({\bf A})_{i,j}& \mbox{---}& j\mbox{th }i\mbox{-dimensional central extension of }\bf A. \\
{\mathcal B}^{i*}_j& \mbox{---}& j \mbox{th }i\mbox{-dimensional nilpotent  bicommutative algebra with identity $xyz=0$} \\
{\mathcal B}^i_j& \mbox{---}& j \mbox{th }i\mbox{-dimensional nilpotent "pure" bicommutative algebra (without identity $xyz=0$)} \\
{\mathfrak{N}}_i& \mbox{---}& i\mbox{th }4\mbox{-dimensional $2$-step nilpotent algebra} \\
{\rm{B}}_i& \mbox{---}& i\mbox{th  non-split non-one-generated }  
5\mbox{-dimensional 
 nilpotent} \\ 
&&\mbox{(non-$2$-step nilpotent) non-commutative bicommutative algebra} \\

\end{array}$$

\subsection{$1$-dimensional central extensions of $4$-dimensional $2$-step nilpotent bicommutative algebras}

\subsubsection{The description of second cohomology space.}

In the following table, we give the description
of the second cohomology space of $4$-dimensional $2$-step nilpotent bicommutative algebras (see, \cite{corr}).

\begin{longtable}{ll llllll}
\hline

\multicolumn{8}{c}{{\bf The list of 2-step nilpotent 4-dimensional bicommutative algebras}}  \\
\hline
 
{${\mathfrak N}_{01}$} &$:$ &  $e_1e_1 = e_2$ &&&&\\ 
\multicolumn{8}{l}{
${\rm H}_{com}^2({\mathfrak N}_{01})= 
\Big\langle 
 [\Delta_{ 1 2}+\Delta_{21}], 
[\Delta_{1 3}+\Delta_{ 31}], 
[\Delta_{ 1 4}+\Delta_{ 41}],  
[\Delta_{3 3}], 
[\Delta_{3 4}+\Delta_{ 43}], [ \Delta_{ 4 4}]

\Big\rangle $} \\
\multicolumn{8}{l}{
${\rm H}_{bicom}^2({\mathfrak N}_{01})=  {\rm H}_{com}^2({\mathfrak N}_{01}) \oplus \Big\langle 
 [\Delta_{21}],[\Delta_{31}], 
[\Delta_{41}],[\Delta_{43}]
\Big\rangle  $}\\
 
\hline
{${\mathfrak N}_{02}$} &$:$ & $e_1e_1 = e_3$& $e_2e_2=e_4$  &&&  \\ 

\multicolumn{8}{l}{
${\rm H}_{com}^2({\mathfrak N}_{02})=
\Big\langle [ \Delta_{12}+\Delta_{21}], [\Delta_{13}+\Delta_{31}], [\Delta_{24}+\Delta_{42}]
 \Big\rangle $}\\ 

\multicolumn{8}{l}{
${\rm H}_{bicom}^2({\mathfrak N}_{02})=
 {\rm H}_{com}^2({\mathfrak N}_{02}) \oplus \Big\langle [\Delta_{21}],[\Delta_{31}],[\Delta_{42}]
 \Big\rangle $}\\
\hline
{${\mathfrak N}_{03}$} &$:$ &  $e_1e_2=  e_3$ & $e_2e_1=-e_3$ &&& \\ 

\multicolumn{8}{l}{
${\rm H}^2({\mathfrak N}_{03})=
\Big\langle [\Delta_{11}],[\Delta_{14}], [\Delta_{21}],[\Delta_{22}],[\Delta_{24}],[\Delta_{41}],[\Delta_{42}],[\Delta_{44}]
 \Big\rangle $}\\

\hline
${\mathfrak N}_{04}^{\alpha}$ &$:$ & $e_1e_1=  e_3$ & $e_1e_2=e_3$& $e_2e_2=\alpha e_3$  &&\\  

\multicolumn{8}{l}{
${\rm H}^2({\mathfrak N}_{04}^{\alpha\neq 0})=
\Big\langle  [\Delta_{12}], [\Delta_{14}], [\Delta_{21}], [\Delta_{22}], [\Delta_{24}], [\Delta_{41}],[ \Delta_{42}], [\Delta_{44}]
 \Big\rangle = \Phi_{\alpha}$ }\\
 \multicolumn{8}{l}{
 ${\rm H}^2({\mathfrak N}_{04}^{0})= \Phi_0 \oplus \Big\langle  [\Delta_{13}], [\Delta_{31}+\Delta_{32}] \Big\rangle $}\\

\hline
${\mathfrak N}_{05}$ &$:$ & $e_1e_1=  e_3$& $e_1e_2=e_3$&  $e_2e_1=e_3$ &&\\ 

\multicolumn{8}{l}{
${\rm H}_{com}^2({\mathfrak N}_{05})=
\Big\langle  [\Delta_{11}], [\Delta_{14}], [\Delta_{21}], [\Delta_{22}], [\Delta_{24}], [\Delta_{41}], [\Delta_{42}], [\Delta_{44}]
\Big\rangle $}\\

\hline
${\mathfrak N}_{06}$ &$:$ & $e_1e_2 = e_4$& $e_3e_1 = e_4$   &&&\\ 

\multicolumn{8}{l}{
${\rm H}^2({\mathfrak N}_{06})=
\Big\langle  [\Delta_{11}], [\Delta_{12}], [\Delta_{13}], [\Delta_{21}], [\Delta_{22}], [\Delta_{23}], [\Delta_{32}], [\Delta_{33}]
\Big\rangle $}\\

\hline
{${\mathfrak N}_{07}$} &$:$ & $e_1e_2 = e_3$ & $e_2e_1 = e_4$ &  $e_2e_2 = -e_3$ &&\\ 

 \multicolumn{8}{l}{
${\rm H}^2({\mathfrak N}_{07})=\Big\langle  [\Delta_{11}],[\Delta_{22}],[\Delta_{13}-\Delta_{23}],[\Delta_{24}], [\Delta_{32}],[\Delta_{41}]
\Big\rangle $}\\

\hline
${\mathfrak N}_{08}^{\alpha}$ &$:$ & $e_1e_1 = e_3$ & $e_1e_2 = e_4$ & $e_2e_1 = -\alpha e_3$ & $e_2e_2 = -e_4$& \\

\multicolumn{8}{l}{
${\rm H}^2({\mathfrak N}_{08}^{\alpha \neq 1})=
\Big\langle  [\Delta_{12}],  [\Delta_{21}], [\Delta_{13}-\alpha \Delta_{23}],  [\Delta_{14}- \Delta_{24}], [\Delta_{31}], [\Delta_{42}]
\Big\rangle= \Phi_{\alpha}$ }\\

\multicolumn{8}{l}{
${\rm H}^2({\mathfrak N}_{08}^{1})= \Phi_{1} \oplus 
\Big\langle  [\Delta_{32}+\Delta_{41}]
\Big\rangle$}\\

\hline
${\mathfrak N}_{09}^{\alpha}$ &$:$ & $e_1e_1 = e_4$ & $e_1e_2 = \alpha e_4$ &  $e_2e_1 = -\alpha e_4$ & $e_2e_2 = e_4$ &  $e_3e_3 = e_4$\\

\multicolumn{8}{l}{
${\rm H}^2({\mathfrak N}_{09}^{\alpha})=
\Big\langle  [\Delta_{12}],  [\Delta_{13}],  [\Delta_{21}],  [\Delta_{22}],  [\Delta_{23}],  [\Delta_{31}],  [\Delta_{32}], [\Delta_{33}]
\Big\rangle $}\\

\hline
${\mathfrak N}_{10}$ &$:$ &  $e_1e_2 = e_4$ & $e_1e_3 = e_4$ & $e_2e_1 = -e_4$ & $e_2e_2 = e_4$ & $e_3e_1 = e_4$  \\

\multicolumn{8}{l}{
${\rm H}^2({\mathfrak N}_{10})=
\Big\langle [\Delta_{11}],  [\Delta_{13}],  [\Delta_{21}],  [\Delta_{22}],  [\Delta_{23}],  [\Delta_{31}],  [\Delta_{32}], [\Delta_{33}] 
\Big\rangle $}\\

\hline
${\mathfrak N}_{11}$ &$:$ &  $e_1e_1 = e_4$ & $e_1e_2 = e_4$ & $e_2e_1 = -e_4$ & $e_3e_3 = e_4$&  \\

\multicolumn{8}{l}{
${\rm H}^2({\mathfrak N}_{11})=
\Big\langle  [\Delta_{12}],  [\Delta_{13}],  [\Delta_{21}],  [\Delta_{22}],  [\Delta_{23}],  [\Delta_{31}],  [\Delta_{32}], [\Delta_{33}] 
\Big\rangle $}\\

\hline
{${\mathfrak N}_{12}$} &$:$ & $e_1e_2 = e_3$ & $e_2e_1 = e_4$  &&& \\ 

 \multicolumn{8}{l}{${\rm H}^2({\mathfrak N}_{12})=
\Big\langle  [\Delta_{11}], [\Delta_{13}], [\Delta_{22}], [\Delta_{24}], [\Delta_{32}], [\Delta_{41}]
\Big\rangle $}\\

\hline
${\mathfrak N}_{13}$ &$:$ & $e_1e_1 = e_4$ & $e_1e_2 = e_3$ & $e_2e_1 = -e_3$ & 
\multicolumn{2}{l}{$e_2e_2=2e_3+e_4$} \\
 
\multicolumn{8}{l}{
${\rm H}^2({\mathfrak N}_{13})=
\Big\langle [\Delta_{21}], [\Delta_{22}], [\Delta_{14}+\Delta_{23}], [\Delta_{13}-2\Delta_{14}-\Delta_{24}], [\Delta_{32}-\Delta_{41}], [\Delta_{31}-2\Delta_{32}+\Delta_{42}]
\Big\rangle $}\\

\hline
{${\mathfrak N}_{14}^{\alpha}$} &$:$ &   $e_1e_2 = e_4$ & $e_2e_1 =\alpha e_4$ & $e_2e_2 = e_3$&& \\ 

\multicolumn{8}{l}{
${\rm H}^2({\mathfrak N}_{14}^{\alpha})=
\Big\langle  [\Delta_{11}],  [\Delta_{21}], [\Delta_{23}], [\Delta_{13}+\Delta_{24}], [\Delta_{32}],[\alpha\Delta_{31}+ \Delta_{42}]
\Big\rangle = \Phi_{\alpha}$}\\

\multicolumn{8}{l}{
${\rm H}^2({\mathfrak N}_{14}^{0})= \Phi_{0} \oplus 
\Big\langle [\Delta_{14}]
\Big\rangle $}\\

\hline

${\mathfrak N}_{15}$ &$:$ &  $e_1e_2 = e_4$ & $e_2e_1 = -e_4$ & $e_3e_3 = e_4$ && \\

\multicolumn{8}{l}{
${\rm H}^2({\mathfrak N}_{15})=
\Big\langle  [\Delta_{11}],  [\Delta_{13}], [\Delta_{21}], [\Delta_{22}], [\Delta_{23}], [\Delta_{31}], [\Delta_{32}], [\Delta_{33}]
\Big\rangle $}\\
\hline

\end{longtable} 

\subsubsection{Central extensions of ${\mathfrak N}_{01}$}
	Let us use the following notations:
	\begin{longtable}{lllllll} 
	$\nabla_1 = [\Delta_{12}+\Delta_{21}],$ & $\nabla_2 = [\Delta_{13}+\Delta_{31}],$ & $\nabla_3 = [\Delta_{14}+\Delta_{41}],$ & $\nabla_4 = [\Delta_{33}],$ & $\nabla_5 = [\Delta_{34}+\Delta_{43}],$ \\
$\nabla_6 = [\Delta_{44}],$ & $\nabla_7 = [\Delta_{21}],$ &$\nabla_8 = [\Delta_{31}],$ & $\nabla_9 = [\Delta_{41}],$ & $\nabla_{10}= [\Delta_{43}].$ 
	\end{longtable}	
	
Take $\theta=\sum\limits_{i=1}^{10}\alpha_i\nabla_i\in {\rm H^2}({\mathfrak N}_{01}).$
	The automorphism group of ${\mathfrak N}_{01}$ consists of invertible matrices of the form
	$$\phi=
	\begin{pmatrix}
	x &  0  & 0 & 0\\
	q &  x^2& r & u\\
	w &  0  & t & k\\
    z &  0  & y & l
	\end{pmatrix}.
	$$
	Since
	$$
	\phi^T\begin{pmatrix}
	0 & \alpha_1   & \alpha_2 & \alpha_3\\
	\alpha_1 +\alpha_7  & 0 & 0 & 0\\
	\alpha_2+\alpha_8 &  0  & \alpha_4 & \alpha_5\\
	\alpha_3+\alpha_9 &  0  & \alpha_5+\alpha_{10} & \alpha_6
	\end{pmatrix} \phi=	\begin{pmatrix}
	\alpha^* & \alpha_1^*   & \alpha_2^* & \alpha_3^*\\
	\alpha_1^* +\alpha_7^*  & 0 & 0 & 0\\
	\alpha_2^*+\alpha_8^* &  0  & \alpha_4^* & \alpha_5^*\\
	\alpha_3^*+\alpha_9^* &  0  & \alpha_5^*+\alpha_{10}^* & \alpha_6^*
	\end{pmatrix},
	$$
	 we have that the action of ${\rm Aut} ({\mathfrak N}_{01})$ on the subspace
$\langle \sum\limits_{i=1}^{10}\alpha_i\nabla_i  \rangle$
is given by
$\langle \sum\limits_{i=1}^{10}\alpha_i^{*}\nabla_i\rangle,$
where
\begin{longtable}{lcl}
$\alpha^*_1$&$=$&$x^3 \alpha _1,$ \\
$\alpha^*_2$&$=$&$r x \alpha _1+y \left(x \alpha _3+w \alpha _5+z \alpha _6\right)+t \left(x \alpha _2+w \alpha _4+z \left(\alpha _5+\alpha _{10}\right)\right),$ \\
$\alpha^*_3$&$=$&$u x \alpha _1+l \left(x \alpha _3+w \alpha _5+z \alpha _6\right)+k \left(x \alpha _2+w \alpha _4+z \left(\alpha _5+\alpha _{10}\right)\right),$ \\
$\alpha_4^*$&$=$&$t^2 \alpha _4+y \left(2 t \alpha _5+y \alpha _6+t \alpha _{10}\right),$\\
$\alpha_5^*$&$=$&$k t \alpha _4+(l t+k y) \alpha _5+y \left(l \alpha _6+k \alpha _{10}\right)$\\
$\alpha_6^*$&$=$&$k^2 \alpha _4+l \left(2 k \alpha _5+l \alpha _6+k \alpha _{10}\right),$\\
$\alpha_7^*$&$=$&$x^3 \alpha _7,$\\
$\alpha_8^*$&$=$&$r x \alpha _7+t x \alpha _8+x y \alpha _9+w y \alpha _{10}-t z \alpha _{10},$\\
$\alpha_9^*$&$=$&$u x \alpha _7+k x \alpha _8+l x \alpha _9+l w \alpha _{10}-k z \alpha _{10},$\\
$\alpha_{10}^*$&$=$&$(l t-k y) \alpha _{10}.$\\
\end{longtable}

 We are interested only in the cases with 
 \begin{center}
$(\alpha_1,\alpha_7)\neq (0,0),$ 
$(\alpha_2,\alpha_4,\alpha_5,\alpha_8, \alpha_{10}) \neq (0,0,0,0,0),$\\
$(\alpha_3,\alpha_5,\alpha_6,\alpha_9, \alpha_{10}) \neq (0,0,0,0,0),$
$(\alpha_7,\alpha_8,\alpha_9, \alpha_{10}) \neq (0,0,0,0).$ 
 \end{center} 

\begin{enumerate}
    \item $\alpha_1=0,\ \alpha_{7}\neq0,$ then choosing $r=-\frac{t x \alpha _8+x y \alpha _9+(w y-t z) \alpha _{10}}{x \alpha _7},$ $u=-\frac{k x \alpha _8+l x \alpha _9+(l w-k z) \alpha _{10}}{x \alpha _7},$ we have $\alpha_8^*=\alpha_9^*=0.$ 
  
  The family of orbits $\langle\alpha_4\nabla_4+\alpha_5\nabla_5+\alpha_6\nabla_6+\alpha_{10}\nabla_{10}\rangle$ gives us characterised structure of three dimensional ideal whose a one dimensional extension of two dimensional subalgebra with basis $\{e_3, e_4\}.$ 
Let us remember the classification of algebras of this type.

\begin{longtable}{ll llllllllllll}
${\mathcal B}^{3*}_{01}$ &:&  $e_1 e_1 = e_2$\\

${\mathcal B}^{3*}_{02}$ &:&  $e_1 e_1 = e_3$ &  $e_2 e_2=e_3$ \\

${\mathcal B}^{3*}_{03}$ &:&   $e_1 e_2=e_3$ & $e_2 e_1=-e_3$   \\

${\mathcal B}^{3*}_{04}(\lambda)$ &:&
$e_1 e_1 = \lambda e_3$  & $e_2 e_1=e_3$  & $e_2 e_2=e_3$   \\

\end{longtable}

Using the classification of three dimensional nilpotent algebras, we may consider following cases. 
   \begin{enumerate}
	\item $\alpha_4=\alpha_5=\alpha_6=\alpha_{10}=0,$ i.e., three dimensional ideal is abelian. Then 
	we may suppose $\alpha_2 \neq 0$ and choosing $y=0,$ $l=\alpha_2,$ $k =-\alpha_3,$ we obtain that  $\alpha_3^*=0,$ which implies $(\alpha_3^*,\alpha_5^*,\alpha_6^*, \alpha_9^*, \alpha_{10}^*) = (0,0,0,0,0).$ Thus, in this case we do not have new algebras. 
	
	\item $\alpha_4=1,$ $\alpha_5=\alpha_6=\alpha_{10}=0,$ i.e., three dimensional ideal is isomorphic to ${\mathcal B}_{01}^{3*}$. Then $\alpha_3 \neq 0$ and choosing $x=1,$ $t=1,$ $k=0,$ $y=0,$ $w = -\alpha_2,$ $l = \frac{\alpha_7}{\alpha_3}$ and  $t=\sqrt{\alpha_7},$ we have the representative 
	$\langle \nabla_3+ \nabla_4+\nabla_7\rangle.$

\item $\alpha_4=\alpha_6=1,$ $\alpha_5=\alpha_{10}=0,$ i.e., three dimensional ideal is isomorphic to ${\mathcal B}_{02}^{3*}$. Then choosing 
$x=\frac{1}{\sqrt[3]{\alpha_7}},$ $k=y=0,$ $l=t=1,$ $w = -\frac{\alpha_2}{\sqrt[3]{\alpha_7}},$ $z = -\frac{\alpha_3}{\sqrt[3]{\alpha_7}},$ 
we have the representative 
	$\langle \nabla_4+ \nabla_6+\nabla_7\rangle.$
	
\item $\alpha_4= \alpha_6=0,$ $\alpha_5=1, \alpha_{10}=-2, $ i.e., three dimensional ideal is isomorphic to ${\mathcal B}_{03}^{3*}$.  Then choosing 
$x=\frac{1}{\sqrt[3]{\alpha_7}},$ $k=y=0,$ $l=t=1,$ $w = -\frac{\alpha_3}{\sqrt[3]{\alpha_7}},$ $z = \frac{\alpha_2}{\sqrt[3]{\alpha_7}},$ 
we have the representative 
	$\langle \nabla_5+ \nabla_7-2\nabla_{10}\rangle.$
	
	\item $\alpha_4=\lambda, \alpha_5=0, \alpha_6=1, \alpha_{10}=1,$ i.e., three dimensional ideal is isomorphic to ${\mathcal B}_{04}^{3*}(\lambda)$. 
	\begin{enumerate}

\item If  $\lambda \neq 0,$ then choosing $x=\frac{1}{\sqrt[3]{\alpha_7}},$ $k=0,$ $y=0,$ 
$l=t=1,$
$z = \frac{\alpha_3}{\sqrt[3]{\alpha_7}},$
 and 
$w = \frac{\alpha_2-\alpha_3}{\lambda \sqrt[3]{\alpha_7}},$ we have the family of representatives 
	$\langle \lambda \nabla_4+ \nabla_6+\nabla_7+\nabla_{10}\rangle_{\lambda\neq 0}.$

	\item If $\lambda = 0$  
and  $\alpha_2=\alpha_3,$ then choosing  
$x=\frac{1}{\sqrt[3]{\alpha_7}},$ $k=0,$ $y=0,$ 
$l=t=1$ and
$z = \frac{\alpha_3}{\sqrt[3]{\alpha_7}},$
   we have the representative 
	$\langle \nabla_6+\nabla_7+\nabla_{10}\rangle.$

\item If $\lambda = 0$  
and  $\alpha_2\neq\alpha_3,$ then choosing  
$x=\frac{(\alpha_2-\alpha_3)^2}{\alpha_7},$ $k=0,$ $y=0,$ 
$l=t=\frac{(\alpha_2-\alpha_3)^3}{\alpha_7}$ and
$z = -\frac{\alpha_3(\alpha_2-\alpha_3)^2}{\alpha_7},$
   we have the representative 
	$\langle \nabla_2+ \nabla_6+\nabla_7+\nabla_{10}\rangle.$ 
	
 	\end{enumerate}
	
 \end{enumerate}
 
  \item $\alpha_1\neq 0,$ then choosing 
 \begin{longtable}{lcl}
  $r$&$=$&$ - \frac{t x \alpha _2+x y \alpha _3+t w \alpha _4+w y \alpha _5+t z \alpha _5+y z \alpha _6+t z \alpha _{10}}{x \alpha _1},$\\ 
  $u$&$ =$&$- \frac{k x \alpha _2+l x \alpha _3+k w \alpha _4+l w \alpha _5+k z \alpha _5+l z \alpha _6+k z \alpha _{10}}{x \alpha _1},$ 
\end{longtable}  
  we have $\alpha_2^*=\alpha_3^*=0.$ 
   
   \begin{enumerate}
	\item $\alpha_4=\alpha_5=\alpha_6=\alpha_{10}=0,$ i.e., three dimensional ideal is abelian. Then 
	we may suppose $\alpha_8 \neq 0$ and choosing $y=0,$ $l=\alpha_8,$ $k =-\alpha_9,$ we obtain that  $\alpha_9^*=0,$ which implies $(\alpha_3^*,\alpha_5^*,\alpha_6^*, \alpha_9^*, \alpha_{10}^*) = (0,0,0,0,0).$ Thus, in this case we do not have new algebras.

	\item $\alpha_4=1,$ $\alpha_5=\alpha_6=\alpha_{10}=0,$ i.e., three dimensional ideal is isomorphic to ${\mathcal B}_{01}^{3*}$. 
	Then  $\alpha_9\neq 0,$ and choosing  
$x=1,$ $k=0,$ $t=\sqrt{\alpha_1},$ $y=-\frac{\sqrt{\alpha_1}\alpha_8}{\alpha_9},$ 
$l=\frac{\alpha_1}{\alpha_9}$ and $w=0,$
   we have the family of representatives 
	$\langle \nabla_1+ \nabla_4+\alpha \nabla_7+\nabla_{9}\rangle.$
	
\item $\alpha_4=\alpha_6=1,$ $\alpha_5=\alpha_{10}=0,$ i.e., three dimensional ideal is isomorphic to ${\mathcal B}_{02}^{3*}$. 
	\begin{enumerate}
	\item $\alpha_7=0,$ then  $\alpha _8^2+\alpha _9^2\neq 0,$ then choosing $x=\frac{\alpha _8^2+\alpha _9^2}{\alpha _1},$ $t=\frac{\alpha _8(\alpha _8^2+\alpha _9^2)}{\alpha _1},$ $y=\frac{\alpha _9(\alpha _8^2+\alpha _9^2)}{\alpha _1},$
	$l=\frac{\alpha _8(\alpha _8^2+\alpha _9^2)}{\alpha _1},$
	$k=-\frac{\alpha _9(\alpha _8^2+\alpha _9^2)}{\alpha _1},$ we have the representative 
	$\langle \nabla_1+ \nabla_4+\nabla_6+\nabla_{8}\rangle.$
	
	\item $\alpha_7=0,$ $\alpha _8^2+\alpha _9^2= 0,$ i.e., $\alpha_9 = \pm i \alpha _8 \neq 0,$ 
	then choosing $x=\sqrt{\alpha_8},$ 
	$t=\frac{\alpha_1} 2,$ $y=\pm \frac{\alpha_1} {2i},$
	$l=\pm ix\alpha_8,$ $k=x\alpha_8,$
	have the representative 
	$\langle \nabla_1+ \nabla_5+\nabla_8\rangle.$
	
	\item $\alpha_7\neq0,$ then choosing $x=1,$ $y=k=0,$
	$t=l=\sqrt{\alpha_7},$
	$z=\frac{\alpha _1 \alpha _9}{\alpha _7},$ $w=\frac{\alpha _1 \alpha _8}{\alpha _7},$
	we have the family of representatives 
	$\langle \alpha \nabla_1+ \nabla_4+\nabla_6+\nabla_{7}\rangle_{\alpha\neq0}.$
	\end{enumerate}
\item $\alpha_4= \alpha_6=0,$ $\alpha_5=1, \alpha_{10}=-2, $ i.e., three dimensional ideal is isomorphic to ${\mathcal B}_{03}^{3*}$.
	\begin{enumerate}
	\item $2 \alpha _1+\alpha _7\neq 0,$ then choosing $x=1,$ $y=k=0,$  $t=\sqrt{\alpha_1},$ $l=1,$ $z-\frac{ \alpha _1 \alpha _8}{2 \alpha _1+\alpha _7}$ and $w=\frac{ \alpha _1 \alpha _9}{2 \alpha _1+\alpha _7},$ we have the family of representatives 
	$\langle  \nabla_1+ \nabla_5+\alpha\nabla_{7}-2\nabla_{10}\rangle_{\alpha\neq-2}.$
	\item $2 \alpha _1+\alpha _7= 0,$ then in case of $(\alpha_8,\alpha_9) = (0,0),$ we have the representative 
	$\langle  \nabla_1+ \nabla_5-2\nabla_{7}-2\nabla_{10}\rangle$ and in case of $(\alpha_8,\alpha_9) \neq (0,0),$ without loss of generality we may assume $\alpha_8 \neq 0$ and choosing $x=1,$ $y=0,$
	$l=\alpha_8,$ $k=-\alpha_9,$
	$t=\frac{\alpha_1}{\alpha_8},$  we have the representative 
	$\langle  \nabla_1+ \nabla_5-2\nabla_{7}+\nabla_{8}-2\nabla_{10}\rangle.$

\end{enumerate}
\item $\alpha_4=\lambda, \alpha_5=0, \alpha_6=1, \alpha_{10}=1,$ i.e., three dimensional ideal is isomorphic to ${\mathcal B}_{04}^{3*}(\lambda)$. 
	\begin{enumerate}

\item  $\alpha _1^2+\alpha _1 \alpha _7+\lambda \alpha _7^2 \neq 0,$ then choosing $x=1,$ $y=k=0$ and 
\begin{center}$t=l=\sqrt{\alpha_7},$
$z=\frac{ \alpha _1 \left(\alpha _1 \alpha _8+\alpha _4 \alpha _7 \alpha _9\right)}{\alpha _1^2+\alpha _1 \alpha _7+\alpha _4 \alpha _7^2},$ 
$w=\frac{\alpha _1 \left(\alpha _7 \left(\alpha _8-\alpha _9\right)+\alpha _1 \alpha _9\right)}{\alpha _1^2+\alpha _1 \alpha _7+\alpha _4 \alpha _7^2},$ 
\end{center} we have the representative
	$\langle \alpha \nabla_1+ \lambda \nabla_4+\nabla_{6}+\nabla_{7}+\nabla_{10}\rangle.$

\item  $\alpha _1^2+\alpha _1 \alpha _7+\lambda \alpha _7^2 = 0,$ then choosing $y=k=0,$ $w=\frac{z \alpha _7}{\alpha _1}-x \alpha _9,$ we have 
$\alpha_9^*=0,$ $\alpha_8^*=\frac{tx}{\alpha_1}(\alpha_1\alpha_8 - \lambda \alpha_7 \alpha_9).$ Thus, in this case we have the representatives $\langle \frac{-1 \pm \sqrt{1-4\lambda}}{2} \nabla_1+ \lambda \nabla_4+\nabla_{6}+\nabla_{7}+\nabla_{10}\rangle$  and $\langle \frac{- 1 \pm \sqrt{1-4\lambda}}{2} \nabla_1+ \lambda \nabla_4+\nabla_{6}+\nabla_{7}+\nabla_{8}+\nabla_{10}\rangle$ depending on $\alpha_1\alpha_8 = \lambda \alpha_7 \alpha_9$ or not.

\end{enumerate}
\end{enumerate}

\end{enumerate}

Summarizing all cases, we have the following distinct orbits
\begin{center}
$\langle \nabla_3+ \nabla_4+\nabla_7\rangle,$  
$\langle \nabla_5+ \nabla_7-2\nabla_{10}\rangle,$  
$\langle \nabla_2+ \nabla_6+\nabla_7+\nabla_{10}\rangle$ 
$\langle \nabla_1+ \nabla_4+\alpha \nabla_7+\nabla_{9}\rangle,$ 
$\langle \nabla_1+ \nabla_4+\nabla_6+\nabla_{8}\rangle,$ 
$\langle \nabla_1+ \nabla_5+\nabla_8\rangle,$ 
	$\langle \alpha \nabla_1+ \nabla_4+\nabla_6+\nabla_{7}\rangle,$ 
$\langle  \nabla_1+\nabla_5+\alpha\nabla_{7}-2\nabla_{10}\rangle,$ \\
$\langle  \nabla_1+ \nabla_5-2\nabla_{7}+\nabla_{8}-2\nabla_{10}\rangle,$  
$\langle \alpha \nabla_1+ \lambda \nabla_4+\nabla_{6}+\nabla_{7}+\nabla_{10}\rangle,$ \\
$\langle \frac{-1 \pm \sqrt{1-4\lambda}}{2} \nabla_1+ \lambda\nabla_4+\nabla_{6}+\nabla_{7}+\nabla_{8}+\nabla_{10}\rangle,$
\end{center}
which gives the following new algebras (see section \ref{secteoA}):

\begin{center}

${\rm B}_{01},$ 
${\rm B}_{02},$  
${\rm B}_{03},$
${\rm B}_{04}^{\alpha},$
${\rm B}_{05},$
${\rm B}_{06},$ 
${\rm B}_{07}^{\alpha},$ 
${\rm B}_{08}^{\alpha},$  
${\rm B}_{09},$  
${\rm B}_{10}^{\alpha, \lambda},$  
${\rm B}_{11}^{\lambda},$ 
${\rm B}_{12}^{\lambda\neq \frac{1}{4}}.$ 

\end{center}

\subsubsection{Central extensions of ${\mathfrak N}_{02}$}
	Let us use the following notations:
	\begin{longtable}{lll} 
	$\nabla_1 = [\Delta_{12}+\Delta_{21}],$ & $\nabla_2 = [\Delta_{13}+\Delta_{31}],$ & $\nabla_3 = [\Delta_{24}+\Delta_{42}],$ \\ $\nabla_4 = [\Delta_{21}],$ & $\nabla_5 = [\Delta_{31}],$ & $\nabla_6 = [\Delta_{42}].$
	\end{longtable}	
	
Take $\theta=\sum\limits_{i=1}^{6}\alpha_i\nabla_i\in {\rm H^2}({\mathfrak N}_{02}).$
	The automorphism group of ${\mathfrak N}_{02}$ consists of invertible matrices of the form
		$$\phi_1=
	\begin{pmatrix}
	x &  0  & 0 & 0\\
	0 &  y  & 0 & 0\\
	z &  u  & x^2 & 0\\
    t &  v  & 0 & y^2
	\end{pmatrix}, \quad \phi_2=
	\begin{pmatrix}
	0 &  x  & 0 & 0\\
	y &  0  & 0 & 0\\
	z &  u  & 0 & x^2\\
    t &  v  & y^2 & 0
	\end{pmatrix}.
	$$
	Since
	$$
	\phi_1^T\begin{pmatrix}
	0 & \alpha_1   & \alpha_2 & 0\\
	\alpha_1 +\alpha_4  & 0 & 0 & \alpha_3\\
	\alpha_2+\alpha_5 &  0  & 0 & 0\\
	0 &  \alpha_3+\alpha_6  & 0 & 0
	\end{pmatrix} \phi_1=	\begin{pmatrix}
	\alpha^* & \alpha_1^*   & \alpha_2^* & 0\\
	\alpha_1^* +\alpha_4^*  & \alpha^{**} & 0 & \alpha_3^*\\
	\alpha_2^*+\alpha_5^* &  0  & 0 & 0\\
	0 &  \alpha_3^*+\alpha_6^*  & 0 & 0
	\end{pmatrix},
	$$
	 we have that the action of ${\rm Aut} ({\mathfrak N}_{02})$ on the subspace
$\langle \sum\limits_{i=1}^{10}\alpha_i\nabla_i  \rangle$
is given by
$\langle \sum\limits_{i=1}^{10}\alpha_i^{*}\nabla_i\rangle,$
where
\begin{longtable}{lcllcllcl}
$\alpha^*_1$&$=$&$x y \alpha _1+u x \alpha _2+t y \left(\alpha _3+\alpha _6\right),$ &
$\alpha^*_2$&$=$&$x^3 \alpha _2,$ &
$\alpha^*_3$&$=$&$y^3 \alpha _3,$ \\
$\alpha_4^*$&$=$&$x y \alpha _4+u x \alpha _5-t y \alpha _6,$&
$\alpha_5^*$&$=$&$x^3 \alpha _5,$&
$\alpha_6^*$&$=$&$y^3 \alpha _6.$\\
\end{longtable}

 We are interested only in the cases with 
 \begin{center}
$(\alpha_3,\alpha_6)\neq (0,0),$  $(\alpha_2,\alpha_5)\neq (0,0),$ $(\alpha_4,\alpha_5,\alpha_6) \neq (0,0,0).$ 
 \end{center} 

\begin{enumerate}
    \item $(\alpha_5,\alpha_6)=(0,0),$ then $\alpha_2\alpha_3\alpha_4\neq0 $ and by choosing 
    $ x=\frac{\alpha_4}{\sqrt[3]{\alpha_2^2\alpha_3}}, $ 
    $ y=\frac{\alpha_4}{\sqrt[3]{\alpha_2\alpha_3^2}}$ 
    and $t=-\frac{x \left(y \alpha _1+u \alpha _2\right)}{y \alpha _3},$ 
we have the representative $ \left\langle \nabla_2+\nabla_3+\nabla_4\right\rangle;$

    \item $(\alpha_5,\alpha_6)\neq (0,0),$ then without loss of generality (maybe with an action of a suitable $\phi_2$), we can suppose $\alpha_5\neq0$ and choosing $u=\frac{t y \alpha _6-x y \alpha _4}{x \alpha _5},$ we have $\alpha_4^*=0$.
    \begin{enumerate}
        \item $\alpha_3\alpha_5+(\alpha_2+\alpha_5)\alpha_6=0,$ then $\alpha_6\neq0.$ 
        \begin{enumerate}
        
            \item if $\alpha_1\neq0,$ then choosing $x=\frac{\alpha_1}{\sqrt[3]{\alpha_5^2\alpha_6}},$ 
            $ y=\frac{\alpha_1}{\sqrt[3]{\alpha_5\alpha_6^2}},$ we have the family of representatives
            $\left\langle\nabla_1+\alpha\nabla_2-(1+ \alpha)\nabla_3+ \nabla_5+\nabla_6\right\rangle;$
        
            \item if $\alpha_1=0,$ then choosing $x=\sqrt[3]{\frac{\alpha_6}{\alpha_5}}, \ y=1,$ we have the family of representatives
            $\left\langle \alpha\nabla_2-(1+ \alpha)\nabla_3+ \nabla_5+\nabla_6\right\rangle.$
        \end{enumerate}
        
        \item $\alpha_3\alpha_5+(\alpha_2+\alpha_5)\alpha_6\neq0,$ then choosing $t=-\frac{x \alpha _1 \alpha _5}{\alpha _3 \alpha _5+\left(\alpha _2+\alpha _5\right) \alpha _6},$ we have $\alpha_1^*=0.$

        \begin{enumerate}
            \item if $\alpha_6=0,$ then choosing $x=\sqrt[3]{\frac{\alpha_3}{\alpha_5}}, \ y=1,$ we have the family of representatives
            $\left\langle\alpha\nabla_2+\nabla_3+ \nabla_5\right\rangle;$
            \item if $\alpha_6\neq0,$ then choosing $x=\sqrt[3]{\frac{\alpha_6}{\alpha_5}}, \ y=1,$ we have the family of representatives
            $\left\langle \alpha\nabla_2+\beta\nabla_3+ \nabla_5+\nabla_6\right\rangle_{\beta\neq-(1+\alpha)}.$
        \end{enumerate}
    \end{enumerate}
\end{enumerate}

Summarizing all cases, we have the following distinct orbits 

\begin{center} 
$\langle \nabla_2+\nabla_3+\nabla_4\rangle,$ 
$\langle\nabla_1+\alpha\nabla_2-(1+ \alpha)\nabla_3+ \nabla_5+\nabla_6\rangle,$

$\langle\alpha\nabla_2+\nabla_3+ \nabla_5\rangle,$ 
$\langle \alpha\nabla_2+\beta\nabla_3+ \nabla_5+\nabla_6\rangle^{O(\alpha, \beta) \simeq O(\beta,\alpha)},$
\end{center}
which gives the following new algebras (see section \ref{secteoA}):

\begin{center}

${\rm B}_{13},$ 
${\rm B}_{14}^{\alpha},$ 
${\rm B}_{15}^{\alpha},$ 
${\rm B}_{16}^{\alpha, \beta}.$
 
\end{center}

\subsubsection{Central extensions of ${\mathfrak N}_{04}^0$}
	Let us use the following notations:
	\begin{longtable}{lllllll} 
	$\nabla_1 = [\Delta_{12}],$ & $\nabla_2 = [\Delta_{13}],$ & $\nabla_3 = [\Delta_{14}],$ & $\nabla_4 = [\Delta_{21}],$ & $\nabla_5 = [\Delta_{22}],$\\
    $\nabla_6 = [\Delta_{24}],$	 & $\nabla_7 = [\Delta_{41}],$ & $\nabla_8 = [\Delta_{42}],$ & $\nabla_{9} = [\Delta_{44}],$ & $\nabla_{10} = [\Delta_{31}+\Delta_{32}].$ 
	\end{longtable}	
	
Take $\theta=\sum\limits_{i=1}^{10}\alpha_i\nabla_i\in {\rm H^2}({\mathfrak N}_{04}^0).$
	The automorphism group of ${\mathfrak N}_{04}^0$ consists of invertible matrices of the form
		$$\phi=
	\begin{pmatrix}
	x &  0  & 0 & 0\\
	y &  x+y  & 0 & 0\\
	z &  t  & x(x+y) & w\\
    u &  v  & 0 & r
	\end{pmatrix}. 
	$$
	Since
	$$
	\phi^T\begin{pmatrix}
	0 & \alpha_1   & \alpha_2 & \alpha_3\\
	\alpha_4  & \alpha_5 & 0 & \alpha_6\\
	\alpha_{10} &  \alpha_{10}  & 0 & 0\\
	\alpha_7 &  \alpha_8  & 0 & \alpha_9
	\end{pmatrix} \phi=	\begin{pmatrix}
	\alpha^{*} & \alpha_1^{*}+\alpha^{*}   & \alpha_2^{*} & \alpha_3^{*}\\
	\alpha_4^{*}  & \alpha_5^{*} & 0 & \alpha_6^{*}\\
	\alpha_{10}^{*} &  \alpha_{10}^{*}  & 0 & 0\\
	\alpha_7^{*} &  \alpha_8^{*}  & 0 & \alpha_9^{*}
	\end{pmatrix},
	$$
	 we have that the action of ${\rm Aut} ({\mathfrak N}_{01})$ on the subspace
$\langle \sum\limits_{i=1}^{10}\alpha_i\nabla_i  \rangle$
is given by
$\langle \sum\limits_{i=1}^{10}\alpha_i^{*}\nabla_i\rangle,$
where
\begin{longtable}{lcl}
$\alpha^*_1$&$=$&$x^2 \alpha _1+x (t-z) \alpha _2-x(u-v)\alpha _3-x y \alpha_4+x y\alpha_5$\\
&&\multicolumn{1}{r}{$-y( u-v)\alpha _6-u x \alpha _7+u x \alpha _8-u(u-v)\alpha _9,$ } \\
$\alpha^*_2$&$=$&$x^2 (x+y) \alpha _2,$ \\
$\alpha^*_3$&$=$&$wx\alpha_2+r(x\alpha_3+y\alpha_6+u\alpha_9),$ \\
$\alpha_4^*$&$=$&$u \left((x+y) \alpha _6+v \alpha _9\right)+x \left((x+y) \alpha _4+v \alpha _7+t \alpha _{10}\right)+y((x+y) \alpha _5+v \alpha _8+t \alpha _{10}),$\\
$\alpha_5^*$&$=$&$v \left((x+y) \alpha _6+v \alpha _9\right)+(x+y) \left((x+y) \alpha _5+v \alpha _8+t \alpha _{10}\right),$\\
$\alpha_6^*$&$=$&$r \left((x+y) \alpha _6+v \alpha _9\right),$\\
$\alpha_7^*$&$=$&$r(x\alpha _7+y\alpha _8+u\alpha _9)+w(x+y) \alpha _{10},$\\
$\alpha_8^*$&$=$&$r v \alpha _9+(x+y)(r \alpha _8+w \alpha _{10}),$\\
$\alpha_9^*$&$=$&$r^2 \alpha _9,$\\
$\alpha_{10}^*$&$=$&$x (x+y)^2 \alpha _{10}.$\\
\end{longtable}

Since we are interested only in the cases with 
$$(\alpha_2,\alpha_{10})\neq (0,0),\quad  (\alpha_3,\alpha_6,\alpha_7,\alpha_8,\alpha_9) \neq (0,0,0,0,0),$$
consider the following subcases:

\begin{enumerate}
    \item $\alpha_{10}=0,$ then $\alpha_2\neq0$ and
  choosing  $w=-\frac{r \left(x \alpha _3+y \alpha _6+u \alpha _9\right)}{x \alpha _2}$ and 
  \begin{center}$t=\frac{-x^2 \alpha _1+x z \alpha _2+u x \alpha _3-v x \alpha _3+x y \alpha _4-x y \alpha _5+u y \alpha _6-v y \alpha _6+u x \alpha _7-u x \alpha _8+u^2 \alpha _9-u v \alpha _9}{x \alpha _2},$ 
 \end{center} we have $\alpha_1^{*}=\alpha_3^{*}=0.$  
    \begin{enumerate}
        \item $\alpha_9\neq0,$ then choosing $u=-\frac{x \alpha _7+y \alpha _8}{\alpha _9},\ v=-\frac{(x+y) \alpha _8}{\alpha _9},$ we have $\alpha_7^*=\alpha_8^*=0.$ 
        \begin{enumerate}
            \item $\alpha_5=\alpha_4=\alpha_6=0,$ then choosing $x=1,\ y=0,\ r=\sqrt{\frac{\alpha_2}{\alpha_9}},$ we have the representative $ \left\langle \nabla_2+ \nabla_9\right\rangle;$ 
            \item $\alpha_5=\alpha_4=0,\ \alpha_6\neq0,$ then choosing $x=1,\ y=\frac{\alpha _2 \alpha _9 - \alpha _6^2}{\alpha _6^2}, \ r=\frac{\alpha _2}{\alpha _6},$ we have the representative $ \left\langle \nabla_2+\nabla_6+ \nabla_9 \right\rangle;$ 
            \item $\alpha_5=0,\ \alpha_4\neq0,\ \alpha_6=0$ then choosing $x=\frac{ \alpha _4}{\alpha _2},\ y=0,\ r=\frac{\sqrt{\alpha_4^{3}}}{\alpha _2 \sqrt{\alpha _9}},$ we have the representative $ \left\langle \nabla_2+\nabla_4+ \nabla_9 \right\rangle;$ 
            \item $\alpha_5=0,\ \alpha_4\neq0,\ \alpha_6\neq0,$ then choosing $x=\frac{\alpha _4}{\alpha _2},\ y=\frac{\alpha _4 \left(\alpha _4 \alpha _9-\alpha _6^2\right)}{\alpha _2 \alpha _6^2},\ r=\frac{\alpha _4^2}{\alpha _2 \alpha _6},$ we have the representative $ \left\langle \nabla_2+\nabla_4+\nabla_6+\nabla_9 \right\rangle;$ 
            \item $\alpha_5\neq0,\ \alpha_4=\alpha_5,$ then choosing $x=1,\ y= \frac{\alpha _2-\alpha _5}{\alpha _5},\ r=\frac{\alpha _2}{\sqrt{\alpha _5\alpha _9}},$ we have the representative $ \left\langle \nabla_2+\nabla_4+\nabla_5+\alpha\nabla_6+ \nabla_9 \right \rangle;$ 
            \item $\alpha_5\neq0,\ \alpha_4\neq\alpha_5,$ then choosing $x=\frac{\alpha _5-\alpha _4}{\alpha _2},\ y=\frac{\alpha _4 \left(\alpha _4-\alpha _5\right)}{\alpha _2 \alpha _5},\ r=\frac{(\alpha _4-\alpha _5)^2}{\alpha _2 \sqrt{\alpha_5\alpha _9}},$  we have the representative $ \left\langle \nabla_2+\nabla_5+\alpha\nabla_6+\nabla_9\right \rangle;$
        \end{enumerate}

\item $\alpha_9=0,\ \alpha_8\neq0,\ \alpha_7=\alpha_8$, then choosing $v=-\frac{x \alpha _4+y \alpha _5+u \alpha _6}{\alpha _8},$ we have $\alpha_4^*=0.$ 
    \begin{enumerate}
         \item $\alpha_6=-\alpha_8,$ $\alpha_5=0,$ then choosing $x=1,\ y=0,\ r=\frac{\alpha_2}{\alpha_5},$ we have the representative $\left\langle \nabla_2-\nabla_6 +\nabla_7+\nabla_8\right\rangle;$
        \item $\alpha_6=-\alpha_8, \alpha_5\neq0,$ then choosing $x=1,\ y=\frac{\alpha _2-\alpha _5}{\alpha _5},\ r=\frac{\alpha_2}{\alpha_8},$ we have the representative $\left\langle \nabla_2+\nabla_5 -\nabla_6 +\nabla_7+\nabla_8\right\rangle;$
        \item $\alpha_6\neq-\alpha_8$, $\alpha_6=0, \ \alpha_5=0,$ then choosing $x=1,\ y=0,\ r=\frac{\alpha_2}{\alpha _8},$  we have the representative $ \left\langle \nabla_2+\nabla_7+ \nabla_8\right\rangle;$
        \item $\alpha_6\neq-\alpha_8$, $\alpha_6=0, \ \alpha_5\neq0,$ then choosing $x=\frac{\alpha _5}{\alpha _2},\ r=\frac{\alpha _5^2}{\alpha_2\alpha_8},$ we have the representative $\left\langle \nabla_2+ \nabla_5+\nabla_7 +\nabla_8 \right\rangle;$ 
        \item $\alpha_6\neq-\alpha_8$, $\alpha_6\neq0,$ then choosing $x=1,\ y=0,\ r=\frac{\alpha_2}{\alpha_8},$ we have the family of representatives $\left\langle \nabla_2+\alpha\nabla_6+\nabla_7+\nabla_8\right \rangle_{\alpha\neq0, -1}.$
          \end{enumerate}

    \item $\alpha_9=0,\ \alpha_8\neq0,\ \alpha_7\neq\alpha_8$, then choosing $y=-\frac{x \alpha _7}{\alpha _8},$ we have $\alpha_7^*=0.$ Hence,
    \begin{enumerate}
        \item $\alpha_6=-\alpha_8,\ \alpha_5=0,$ then choosing $x=1,\ u=\frac{\alpha _4}{\alpha _8},\ r=\frac{\alpha_2}{\alpha_8},$ we have the representative $\left\langle\nabla_2- \nabla_6+\nabla_8 \right\rangle;$
        \item $\alpha_6=-\alpha_8,\ \alpha_5\neq0,$ then choosing $x=\frac{\alpha_5}{\alpha_2},\ u=\frac{\alpha_4\alpha_5}{\alpha_2\alpha_8},\ r=\frac{\alpha_5^2}{\alpha_2\alpha_8},$ we have the representative $\left\langle\nabla_2+\nabla_5- \nabla_6+\nabla_8\right\rangle;$
        \item $\alpha_6\neq-\alpha_8,$ $\alpha_6=\alpha_4=0,$ then choosing $x=1,\ v=-\frac{\alpha_5}{\alpha_8},\ r=\frac{\alpha_2}{\alpha_8},$ we have the representative $\left\langle\nabla_2+\nabla_8 \right\rangle;$
        \item $\alpha_6\neq-\alpha_8,$ $\alpha_6=0,\ \alpha_4\neq0,$ then choosing $x=\frac{\alpha _4}{\alpha _2},\ v=-\frac{\alpha _4 \alpha _5}{\alpha _2 \alpha _8},\ r=\frac{\alpha _4^2}{\alpha _2 \alpha _8},$ we have the representative $\left\langle\nabla_2+\nabla_4 +\nabla_8 \right\rangle;$
        \item $\alpha_6\neq-\alpha_8,$ $\alpha_6\neq0,$ then choosing $x=1,\ u=-\frac{\alpha_4}{\alpha_6},\ v=-\frac{\alpha _5}{\alpha _6+\alpha _8},\ r=\frac{\alpha_2}{\alpha_8},$ we have the representative $\left\langle\nabla_2+\alpha \nabla_6 +\nabla_8\right\rangle_{\alpha\neq 0, -1}.$
        \end{enumerate}
    \item $\alpha_9=\alpha_8=0,\ \alpha_7\neq0,$  then choosing $v=-\frac{(x+y) \left(x \alpha _4+y \alpha _5+u \alpha _6\right)}{x \alpha _7}$, we have $\alpha_4^*=0.$ Hence, 
    \begin{enumerate}
        \item $\alpha_6=\alpha_5=0,$ then choosing $x=1,\ y=0,\ r=\frac{\alpha_2}{\alpha_7},$ we have the representative $\left\langle\nabla_2+\nabla_7\right \rangle;$
        \item $\alpha_6=0,\ \alpha_5\neq0,$ then choosing $x=\frac{\alpha_5}{\alpha_2},\ y=0,\ r=\frac{\alpha_5^2}{\alpha_2\alpha_7},$ we have the representative $\left\langle\nabla_2+\nabla_5 +\nabla_7\right\rangle;$
        \item $\alpha_6\neq0,$ then choosing $x=1,\ y=\frac{\alpha_7-\alpha_6}{\alpha_6},\ u=\frac{\alpha_5}{\alpha_6},\ r=\frac{\alpha_2}{\alpha_6},$ we have the representative $\left\langle\nabla_2+\nabla_6 +\nabla_7\right\rangle.$
    \end{enumerate}
    \item $\alpha_9=\alpha_8=\alpha_7=0,\ \alpha_6\neq0,$ then choosing $x=1,\ y=0,\ u=-\frac{\alpha _4}{\alpha _6},\ v=-\frac{\alpha_5}{\alpha _6},\ r=\frac{\alpha_2}{\alpha_6},$ we have the representative $\left\langle\nabla_2+\nabla_6\right\rangle.$
    \end{enumerate}
\item $\alpha_{10}\neq0,$ then choosing $ w=-\frac{r \left((x+y) \alpha _8+v \alpha _9\right)}{(x+y) \alpha _{10}}$ and 
\begin{center}$t=-\frac{x (x+y) \alpha _4+y (x+y) \alpha _5+u x \alpha _6+u y \alpha _6+v x \alpha _7+v y \alpha _8+u v \alpha _9}{(x+y) \alpha _{10}},$ \end{center} we have  $\alpha_{4}^*=\alpha_{8}^*=0$. Now we consider following subcases:
    \begin{enumerate}
    \item $\alpha_9\neq0,$ then  choosing $u=-\frac{(x+y) \alpha _6+x \alpha _7}{\alpha _9},\ v=-\frac{(x+y) \alpha _6}{\alpha _9},$ we have  $\alpha_{6}^*=\alpha_{7}^*=0$. Hence, we can suppose $\alpha_4=\alpha_6=\alpha_7=\alpha_8=0$ and consider following cases:
        \begin{enumerate}
        \item $\alpha_2=\alpha_5=\alpha_1=\alpha_3=0,$ then choosing $x=1,\ y=0,\  r=\sqrt{\frac{\alpha_{10}}{\alpha_9}},$ we have the representative $\left\langle\nabla_9+\nabla_{10}\right \rangle;$
        \item $\alpha_2=\alpha_5=\alpha_1=0,\ \alpha_3\neq0,$ then choosing $x=1,\ y=\frac{\alpha_3}{\sqrt{\alpha_9\alpha_{10}}}-1, \  r=\frac{\alpha_3}{\alpha _9},$ we have the representative $\left\langle\nabla_3+\nabla_9+\nabla_{10}\right \rangle;$
        \item $\alpha_2=\alpha_5=0,\ \alpha_1\neq0,$ then choosing $x=1,\ y=\sqrt{\frac{\alpha _1}{\alpha _{10}}}-1,\ r=\sqrt{\frac{\alpha _1}{\alpha _{9}}},$ we have the family of  representatives $\left\langle\nabla_1 +\alpha\nabla_3 +\nabla_9+\nabla_{10}\right \rangle^{O(\alpha)\simeq O(-\alpha)};$ 
        \item $\alpha_2=0,\ \alpha_5\neq0,\ \alpha_1=\alpha_5,\ \alpha_3=0,$ then choosing $x=\frac{\alpha _5}{\alpha _{10}},\ y=0,\  r=\frac{\alpha_5 \sqrt{\alpha_5}}{\alpha_{10}\sqrt{\alpha _9}},$ we have the representative $\left\langle\nabla_1+ \nabla_5+ \nabla_9 + \nabla_{10}\right \rangle;$
        \item $\alpha_2=0,\ \alpha_5\neq0,\ \alpha_1=\alpha_5,\ \alpha_3\neq0,$ then choosing $x=\frac{\alpha _5^2 \alpha _9}{\alpha _3^2 \alpha _{10}},$ $y=\frac{\alpha _5 \left(\alpha _3^2-\alpha _5 \alpha _9\right)}{\alpha _3^2 \alpha _{10}},$ $r=\frac{\alpha _5^2}{\alpha _3 \alpha _{10}},$ we have the family of representatives $\left\langle\nabla_1+\nabla_3+ \nabla_5+ \nabla_9 +\nabla_{10}\right\rangle;$
        \item $\alpha_2=0,\ \alpha_5\neq0,\ \alpha_1\neq\alpha_5,$ then choosing \begin{center}$x=-\frac{\alpha _5^2}{\left(\alpha _1-\alpha _5\right) \alpha _{10}},$ $y=\frac{\alpha _1 \alpha _5}{\left(\alpha _1-\alpha _5\right) \alpha _{10}},$  $r=\frac{\alpha _5^2}{\sqrt{(\alpha _5-\alpha _1)\alpha _9} \alpha _{10}},$\end{center} we have the family of  representatives $\left\langle\alpha \nabla_3+\nabla_5+ \nabla_9 +\nabla_{10}\right\rangle^{O(\alpha)\simeq O(-\alpha)};$ 
        \item $\alpha_2\neq0,\ \alpha_5=\alpha_3=0,$ then choosing $x=1,\ y=\frac{\alpha_2-\alpha_{10}}{\alpha _{10}},\ z=\frac{\alpha_1}{\alpha_2},\ r=\frac{\alpha _2} {\sqrt{\alpha_9\alpha_{10}}},$ we have the representative 
 $\left\langle\nabla_2+\nabla_9+\nabla_{10}\right \rangle;$ 
        \item $\alpha_2\neq0,\ \alpha_5=0,\ \alpha_3\neq0,$ then choosing \begin{center}$x=\frac{\alpha _3^2 \alpha _{10}}{\alpha _2^2 \alpha _9},\ y=\frac{\alpha _3^2 \left(\alpha _2-\alpha _{10}\right)}{\alpha _2^2 \alpha _9},\ z=\frac{\alpha _1 \alpha _3^2 \alpha _{10}}{\alpha _2^3 \alpha _9},\   r=\frac{\alpha _3^3 \alpha _{10}}{\alpha _2^2 \alpha _9^2},$\end{center} we have the representative $\left\langle\nabla_2+ \nabla_3+\nabla_9+ \nabla_{10} \right\rangle;$ 
        \item $\alpha_2\neq0,\ \alpha_5\neq0,$ then choosing   
        \begin{center} $x=\frac{\alpha _5}{\alpha _2},$ $ y=\frac{(\alpha_2-\alpha_{10})\alpha_5}{\alpha_2\alpha _{10}},$ 
        $ z=-\frac{\alpha _5 \left(\alpha _2^2 \alpha _5-2 \alpha _2 \alpha _5 \alpha _{10}-\left(\alpha _1-\alpha _5\right) \alpha _{10}^2\right)}{\alpha _2^2 \alpha _{10}^2},$  
        $r=\frac{\alpha _5\sqrt{\alpha_5}} {\sqrt{\alpha_2\alpha_9\alpha_{10}}},$\end{center} we have the family of representatives $\left\langle\nabla_2+\alpha\nabla_3+\nabla_5+\nabla_9+ \nabla_{10} \right\rangle^{O(\alpha)\simeq O(-\alpha)}.$
        \end{enumerate}
    \item $\alpha_9=0,\ \alpha_6\neq0,$ then  choosing $u=\frac{x (x+y) \alpha _5+v \left((x+y) \alpha _6-x \alpha _7\right)}{(x+y) \alpha _6},$ we have  $\alpha_{5}^*=0$. Hence, we have $\alpha_4=\alpha_5=\alpha_8=\alpha_9=0$ and consider following cases:
    \begin{enumerate}
        \item $\alpha_2=0,\ \alpha_3=\alpha_{6},\ \alpha_7=0,\ \alpha_1=0,$ then choosing $x=1,\  y=0,\ r=\frac{\alpha_{10}}{\alpha_{6}},$ we have the representative $\left\langle\nabla_3+\nabla_6 +\nabla_{10}\right\rangle;$
        \item $\alpha_2=0,\ \alpha_3=\alpha_{6},\ \alpha_7=0,\ \alpha_1\neq0,$ then choosing $x=1,$  $y=\sqrt{\frac{\alpha _1}{\alpha _{10}}}-1, $ $r=\frac{\sqrt{\alpha_1\alpha_{10}}}{\alpha _6},$ we have the representative $\left\langle\nabla_1+\nabla_3 +\nabla_6 +\nabla_{10}\right\rangle;$
        \item $\alpha_2=0,\ \alpha_3=\alpha_{6},\ \alpha_7\neq0,$ then choosing $x=1,\  y=\frac{\alpha_7-\alpha_6}{\alpha_6},\ u=\frac{\alpha_1}{\alpha_7},\  r=\frac{\alpha _7 \alpha _{10}}{\alpha _6^2},$ we have the representative $\left\langle\nabla_3 +\nabla_6 +\nabla_7+ \nabla_{10}\right\rangle;$
        \item $\alpha_2=0,\ \alpha_3\neq\alpha_{6},\ \alpha _7=\alpha_1=0,$ then choosing $x=1,\  y=-\frac{\alpha_3}{\alpha_6},\ r=-\frac{\left(\alpha _3-\alpha _6\right) \alpha _{10}}{\alpha _6^2},$ we have the representative $\left\langle\nabla_6+ \nabla_{10}\right\rangle;$
        \item $\alpha_2=0,\ \alpha_3\neq\alpha_{6},\ \alpha _7=0,\ \alpha_1\neq0,$ then choosing \begin{center} $x=\frac{\alpha _1 \alpha _6^2}{\left(\alpha _6-\alpha _3\right){}^2 \alpha _{10}},\  y=-\frac{\alpha _1 \alpha_3\alpha _6}{\left(\alpha _6-\alpha _3\right){}^2 \alpha _{10}},\ r=\frac{\alpha _1^2 \alpha _6^2}{\left(\alpha _6-\alpha _3\right){}^3 \alpha _{10}},$ \end{center} we have the representative $\left\langle\nabla_1+\nabla_6 +\nabla_{10}\right\rangle;$
        \item $\alpha_2=0,\ \alpha_3\neq\alpha_{6},\ \alpha_7\neq0, \alpha _3+\alpha _7=\alpha _6,\ \alpha_1=0,$ then choosing $x=1,$   
        $y=-\frac{\alpha_3}{\alpha _6},$ 
        $r=\frac{\left(\alpha _6-\alpha _3\right) \alpha _{10}}{\alpha _6^2},$ we have the representative  $\left\langle\nabla_6 +\nabla_7 +\nabla_{10}\right\rangle;$
        \item $\alpha_2=0,\ \alpha_3\neq\alpha_{6},\ \alpha_7\neq0, \ \alpha _3+\alpha _7=\alpha _6,$  $\alpha_1\neq0,$ then choosing $x=\frac{\alpha _1 \alpha _6^2}{\left(\alpha _6-\alpha _3\right)^2 \alpha _{10}},$ $y=-\frac{\alpha _1 \alpha_3\alpha _6}{\left(\alpha _6-\alpha _3\right)^2\alpha _{10}},$ $r=\frac{\alpha _1^2 \alpha _6^2}{\left(\alpha _6-\alpha _3\right)^3 \alpha _{10}},$ we have the representative   $\left\langle\nabla_1+\nabla_6 +\nabla_7 +\nabla_{10}\right\rangle;$
        \item $\alpha_2=0,\ \alpha_3\neq\alpha_{6},\ \alpha_7\neq0, \alpha _3+\alpha _7\neq\alpha _6,$ then choosing $x=1,$  $y=-\frac{\alpha_3}{\alpha_{6}},$ $r=\frac{(\alpha_6-\alpha _3)\alpha_{10}}{\alpha _6^2},$
         we have the family of representatives  $\left\langle\nabla_6 +\alpha\nabla_7 +\nabla_{10}\right\rangle_{\alpha\neq0,1};$
        \item $\alpha_2\neq0,$ then choosing $z=\frac{x \alpha _1}{\alpha _2},\ v=0,$ we have $\alpha_1^*=0.$
        \begin{enumerate}
            \item $\alpha_3=\alpha_6,$ then choosing $x=1,\ y=\frac{\alpha_2-\alpha_{10}}{\alpha_{10}}, \  r=\frac{\alpha_{2}}{\alpha _6},$ we have the family of representatives $\left\langle\nabla_2+\nabla_3+ \nabla_6+\alpha\nabla_7+\nabla_{10}\right\rangle;$
            \item $\alpha_3\neq\alpha_6,$ then choosing $x=1,\ y=-\frac{\alpha_3}{\alpha_{6}}, r=\frac{(\alpha_{6}-\alpha_3)\alpha_{10}}{\alpha _6^2},$ we have the family of representatives $\left\langle\beta\nabla_2+\nabla_6+ \alpha\nabla_7+ \nabla_{10}\right\rangle_{\beta\neq0}.$
        \end{enumerate}
    \end{enumerate}
    \item $\alpha_9=\alpha_6=0,\ \alpha_7\neq0,$ then choosing $v=\frac{(x+y) \alpha _5}{\alpha _7},$ we have $\alpha_5^*=0.$ Hence, we have $\alpha_4=\alpha_5=\alpha_6=\alpha_8=\alpha_9=0$ and consider following cases:
    \begin{enumerate}
        \item $\alpha_2=0,\ \alpha_3+\alpha_7=0,\ \alpha_1=0,$ then choosing $x=1,\ y=0, \ r=\frac{\alpha_{10}}{\alpha_7},$ we have the representative $\left\langle-\nabla_3+\nabla_7+\nabla_{10}\right\rangle;$
        \item $\alpha_{2}=0,\ \alpha_3+\alpha_{7}=0,\ \alpha_1\neq0,$ then choosing $x=\frac{\alpha_{1}}{\alpha_{10}},\ y=0, \ r=\frac{\alpha _1^2}{\alpha_7\alpha_{10}},$ we have the representative $\left\langle \nabla_1-\nabla_3+\nabla_7 +\nabla_{10}\right\rangle;$
        \item $\alpha_2=0,\ \alpha_3+\alpha_{7}\neq0,$ then choosing $x=1,\ y=0, \ u=\frac{\alpha_1}{\alpha_3+\alpha_7},\  r=\frac{\alpha_{10}}{\alpha_7},$ we have the family of representatives $\left\langle \alpha\nabla_3+\nabla_7+ \nabla_{10}\right\rangle_{\alpha\neq -1};$
        \item $\alpha_2\neq0,$ then choosing $x=1,\ y=\frac{\alpha_2- \alpha_{10}}{\alpha_{10}},\ z=\frac{\alpha _1}{\alpha_2},\ r=\frac{\alpha _2^2}{\alpha _7 \alpha _{10}},\ u=0,$ we have the family of representatives $\left\langle\nabla_2+ \alpha\nabla_3+ \nabla_7+ \nabla_{10}\right\rangle.$
    \end{enumerate}
    \item $\alpha_9=\alpha_7=\alpha_6=0,$ then $\alpha_3\neq0,$ and choosing $u=\frac{\left(x \alpha _1+v \alpha _3+y \alpha _5\right) \alpha _{10}-\alpha _2 \left(y \alpha _5+z \alpha _{10}\right)}{\alpha _3 \alpha _{10}}$, we obtain $\alpha_1^{*}=0$. Hence, we have $\alpha_1= \alpha_4=\alpha_6=\alpha_7=\alpha_8=\alpha_9=0$ and consider following cases: \begin{enumerate}
        \item $\alpha_2=\alpha_5=0,$ then choosing $x=1,\ y=0,\ r=\frac{\alpha_{10}}{\alpha_3},$ we have the representative $\left\langle\nabla_3+\nabla_{10} \right\rangle;$
        \item $\alpha_2=0,\ \alpha_5\neq0,$ then choosing $x=\frac{\alpha_{5}}{\alpha_{10}},\ y=0,\ r=\frac{\alpha _5^2}{\alpha _3 \alpha _{10}},$ we have the representative $\left\langle\nabla_3+\nabla_5+\nabla_{10} \right\rangle;$
        \item $\alpha_2\neq0,\ \alpha_5=0,$ then choosing $x=1,\ y=\frac{\alpha_2-\alpha_{10}}{\alpha _{10}},\ r=\frac{\alpha _2^2}{\alpha _3\alpha _{10}},$ we have the representative $\left\langle\nabla_2+\nabla_3+\nabla_{10} \right\rangle;$
        \item $\alpha_2\neq0,\ \alpha_5\neq0,$ then choosing $x=\frac{\alpha _5}{\alpha _2},\ y=\frac{ (\alpha_2-\alpha_{10})\alpha_5}{\alpha_2\alpha_{10}},\ r=\frac{\alpha _5^2}{\alpha _3 \alpha _{10}},$ we have the representative $\left\langle\nabla_2+\nabla_3+ \nabla_5+\nabla_{10} \right\rangle.$
    \end{enumerate}
  \end{enumerate}
\end{enumerate}

Summarizing all cases, we have the following distinct orbits 

\begin{center} 
$\langle \nabla_2+\nabla_9\rangle,$ 
$\langle\nabla_2+ \nabla_6+\nabla_9\rangle,$  
$\langle\nabla_2+ \nabla_4+\nabla_9\rangle,$
$\langle \nabla_2+\nabla_4+\nabla_6+\nabla_9\rangle,$ 
$\langle \nabla_2+\nabla_4+\nabla_5+\alpha\nabla_6+\nabla_9\rangle,$  $\langle\nabla_2+\nabla_5+\alpha\nabla_6+\nabla_9\rangle,$
$\langle \nabla_2+\nabla_5+\nabla_7+\nabla_8\rangle,$ 
$\langle \nabla_2+\nabla_5-\nabla_6+\nabla_7+\nabla_8\rangle,$  $\langle\nabla_2+\alpha\nabla_6+\nabla_7+\nabla_8\rangle,$
$\langle \nabla_2+\nabla_4+\nabla_8\rangle,$ 
$\langle \nabla_2+\nabla_5-\nabla_6+\nabla_8\rangle,$  $\langle\nabla_2+\alpha\nabla_6+\nabla_8\rangle,$
$\langle \nabla_2+\nabla_7\rangle,$ 
$\langle\nabla_2+\nabla_5+\nabla_7\rangle,$ 
$\langle\nabla_2+\nabla_6+\nabla_7\rangle,$
$\langle \nabla_2+\nabla_6\rangle,$ 
 $\left\langle\nabla_9+\nabla_{10}\right \rangle,$ 
$\left\langle\nabla_3+\nabla_9+\nabla_{10}\right \rangle,$ 
$\left\langle\nabla_1 +\alpha\nabla_3 +\nabla_9+\nabla_{10}\right \rangle^{O(\alpha)\simeq O(-\alpha)},$  
$\left\langle\nabla_1+ \nabla_5+ \nabla_9 + \nabla_{10}\right \rangle,$  
$\left\langle\nabla_1+\nabla_3+ \nabla_5+ \nabla_9 +\nabla_{10}\right\rangle,$ 
$\left\langle\alpha \nabla_3+\nabla_5+ \nabla_9 +\nabla_{10}\right\rangle^{O(\alpha)\simeq O(-\alpha)},$ 
$\left\langle\nabla_2+\nabla_9+\nabla_{10}\right \rangle,$ 
$\left\langle\nabla_2+ \nabla_3+\nabla_9+ \nabla_{10} \right\rangle,$ 
$\left\langle\nabla_2+\alpha\nabla_3+\nabla_5+\nabla_9+ \nabla_{10} \right\rangle^{O(\alpha)\simeq O(-\alpha)},$ 
$\left\langle\nabla_3+\nabla_6 +\nabla_{10}\right\rangle,$ 
$\left\langle\nabla_1+\nabla_3 +\nabla_6 +\nabla_{10}\right\rangle,$ 
$\left\langle\nabla_3 +\nabla_6 +\nabla_7+ \nabla_{10}\right\rangle,$ 
$\left\langle\nabla_1+\nabla_6 +\nabla_{10}\right\rangle,$ 
$\left\langle\nabla_1+\nabla_6 +\nabla_7 +\nabla_{10}\right\rangle,$ 
$\left\langle\nabla_2+\nabla_3+ \nabla_6+\alpha\nabla_7+\nabla_{10}\right\rangle,$ 
$\left\langle\beta\nabla_2+\nabla_6+ \alpha\nabla_7+ \nabla_{10}\right\rangle,$ 
$\left\langle \nabla_1-\nabla_3+\nabla_7 +\nabla_{10}\right\rangle,$ 
$\left\langle \alpha\nabla_3+\nabla_7+ \nabla_{10}\right\rangle,$ 
$\left\langle\nabla_2+ \alpha\nabla_3+ \nabla_7+ \nabla_{10}\right\rangle,$ 
$\left\langle\nabla_3+\nabla_{10} \right\rangle,$ 
$\left\langle\nabla_3+\nabla_5+\nabla_{10} \right\rangle,$ 
$\left\langle\nabla_2+\nabla_3+\nabla_{10} \right\rangle,$ 
$\left\langle\nabla_2+\nabla_3+ \nabla_5+\nabla_{10} \right\rangle,$
\end{center}
which gives the following new algebras (see section \ref{secteoA}):

\begin{center}

${\rm B}_{17},$ 
${\rm B}_{18},$ 
${\rm B}_{19},$ 
${\rm B}_{20},$ 
${\rm B}_{21}^{\alpha},$ 
${\rm B}_{22}^{\alpha},$ 
${\rm B}_{23},$ 
${\rm B}_{24},$ 
${\rm B}_{25}^{\alpha},$ 
${\rm B}_{26},$ 
${\rm B}_{27},$ 
${\rm B}_{28}^{\alpha},$ 
${\rm B}_{29},$ 
${\rm B}_{30},$ 
${\rm B}_{31},$ 
${\rm B}_{32},$ 
${\rm B}_{33},$ 
${\rm B}_{34},$ 
${\rm B}_{35}^{\alpha},$ 
${\rm B}_{36},$ 
${\rm B}_{37},$ 
${\rm B}_{38}^{\alpha},$ 
${\rm B}_{39},$ 
${\rm B}_{40},$ 
${\rm B}_{41}^{\alpha},$ 
${\rm B}_{42},$ 
${\rm B}_{43},$ 
${\rm B}_{44},$
${\rm B}_{45},$ 
${\rm B}_{46},$ 
${\rm B}_{47}^{\alpha},$ 
${\rm B}_{48}^{\alpha, \beta},$ 
${\rm B}_{49},$ 
${\rm B}_{50}^{\alpha},$
${\rm B}_{51}^{\alpha},$ 
${\rm B}_{52},$ 
${\rm B}_{53},$ 
${\rm B}_{54},$ 
${\rm B}_{55}.$

\end{center}

\subsubsection{Central extensions of ${\mathfrak N}_{07}$}
	Let us use the following notations:
	\begin{longtable}{lllllll} 
	$\nabla_1 = [\Delta_{11}],$ & $\nabla_2 = [\Delta_{22}],$ & $\nabla_3 = [\Delta_{13}-\Delta_{23}],$ \\
	$\nabla_4 = [\Delta_{24}],$ & $\nabla_5 = [\Delta_{32}],$ & $\nabla_6 = [\Delta_{41}].$
	\end{longtable}	
	
Take $\theta=\sum\limits_{i=1}^{6}\alpha_i\nabla_i\in {\rm H^2}({\mathfrak N}_{07}).$
	The automorphism group of ${\mathfrak N}_{07}$ consists of invertible matrices of the form
		$$\phi=
	\begin{pmatrix}
	x &  0  & 0 & 0\\
	0 &  x  & 0 & 0\\
	z &  u  & x^2 & 0\\
    t &  v  & 0 & x^2
	\end{pmatrix}.
	$$
	Since
	$$
	\phi^T\begin{pmatrix}
	\alpha_1& 0   & \alpha_3 &0\\
	0 & \alpha_2 & -\alpha_3 & \alpha_4\\
	0 &  \alpha_5  & 0 & 0\\
	\alpha_6 &  0  & 0 & 0
	\end{pmatrix} \phi=	\begin{pmatrix}
	\alpha_1^* & \alpha^*   & \alpha_3^* & 0\\
	\alpha^{**}  & -\alpha^*+\alpha_2^* & -\alpha_3^* & \alpha_4^*\\
	0 &  \alpha_5^*  & 0 & 0\\
	\alpha_6^* &  0  & 0 & 0
	\end{pmatrix},
	$$
	 we have that the action of ${\rm Aut} ({\mathfrak N}_{07})$ on the subspace
$\langle \sum\limits_{i=1}^{6}\alpha_i\nabla_i  \rangle$
is given by
$\langle \sum\limits_{i=1}^{6}\alpha_i^{*}\nabla_i\rangle,$
where
\begin{longtable}{lcllcllcl}
$\alpha^*_1$&$=$&$x (x \alpha _1+z \alpha _3+t \alpha _6),$ &
$\alpha^*_3$&$=$&$x^3 \alpha _3,$ &
$\alpha_5^*$&$=$&$x^3 \alpha _5,$\\

$\alpha^*_2$&$=$&$x (x \alpha _2+v \alpha _4+(u+z) \alpha _5),$ &
$\alpha_4^*$&$=$&$x^3 \alpha _4,$&
$\alpha_6^*$&$=$&$x^3 \alpha _6.$
\end{longtable}

 We are interested only in the cases with 
 $(\alpha_3,\alpha_5)\neq (0,0),$  $(\alpha_4,\alpha_6)\neq (0,0).$

\begin{enumerate}
    \item $\alpha_5\neq0,$ then choosing $u=-\frac{x \alpha _2+v \alpha _4+z \alpha _5}{\alpha _5},$ we have $\alpha_2^*=0.$ Now we consider following subcases:
    \begin{enumerate}
        \item $\alpha_6\neq0,$ then choosing $x=1,\ z=0,\ t=-\frac{\alpha _1}{\alpha _6},$ we have the family of representatives $\left\langle\beta\nabla_3+\alpha\nabla_4+\nabla_5 +\gamma \nabla_6\right\rangle_{\gamma\neq0};$
        \item $\alpha_6=0, \alpha_3\neq0,$ then choosing $x=1,\ z=-\frac{\alpha _1}{\alpha _3},$ we have the family of representatives $\left\langle\beta\nabla_3+\alpha\nabla_4+\nabla_5 \right\rangle_{\alpha\beta\neq0};$ 
        \item $\alpha_6=0, \alpha_3=0,\ \alpha_1\neq0,$ then choosing $x=\frac{\alpha_1}{\alpha _3},$ we have the family of representatives $\left\langle\nabla_1+\alpha\nabla_4+\nabla_5 \right\rangle_{\alpha\neq0};$
        \item $\alpha_6=0, \alpha_3=0,\ \alpha_1=0,$ then we have the family of representatives $\left\langle\alpha\nabla_4+\nabla_5\right\rangle_{\alpha\neq0}.$ 
    \end{enumerate}
    \item $\alpha_5=0,\ \alpha_3\neq0,$ then choosing $z=-\frac{x \alpha _1+t \alpha _6}{\alpha _3},$ we have $\alpha_1^*=0.$ 
    \begin{enumerate}
        \item $\alpha_4\neq0,$ then choosing $x=1,\ v=-\frac{\alpha_2}{\alpha_4},$ we have the family of representatives 
        $\left\langle\nabla_3+\beta\nabla_4+ \alpha\nabla_6\right\rangle_{\beta\neq0};$
        \item $\alpha_4=0, \alpha_2\neq0,$ then choosing $x=\frac{\alpha_2}{\alpha_3},$ we have the family of representatives $\left\langle\nabla_2+\nabla_3+\alpha\nabla_6 \right\rangle_{\alpha\neq0};$
        \item $\alpha_4=0, \alpha_2=0,$ then we have the family of representatives $\left\langle\nabla_3+\alpha\nabla_6\right\rangle_{\alpha\neq0}.$ 
    \end{enumerate}
\end{enumerate}

Summarizing all cases, we have the following distinct orbits 

\begin{center} 
$\langle\nabla_1+\alpha\nabla_4+\nabla_5\rangle_{\alpha\neq0},$  $\langle\gamma\nabla_3+\alpha\nabla_4+\nabla_5+\beta\nabla_6\rangle_{(\alpha,\beta)\neq(0,0)},$
$\langle\nabla_3+\alpha\nabla_4+\beta\nabla_6\rangle_{(\alpha,\beta)\neq(0,0)},$  $\langle\nabla_2+\nabla_3+\alpha\nabla_6\rangle_{\alpha\neq0},$
\end{center}
which gives the following new algebras (see section \ref{secteoA}):

\begin{center}
${\rm B}_{56}^{\alpha \neq 0},$
${\rm  B}_{57}^{(\alpha, \beta, \gamma) \neq (0,0,\gamma)},$
${\rm  B}_{58}^{(\alpha, \beta)\neq (0,0)},$
${\rm  B}_{59}^{\alpha \neq 0}.$
\end{center}

\subsubsection{Central extensions of ${\mathfrak N}_{08}^{\alpha\neq1}$}
	Let us use the following notations:
	\begin{longtable}{lllllll} 
	$\nabla_1 = [\Delta_{12}],$ & $\nabla_2 = [\Delta_{21}],$ & $\nabla_3 = [\Delta_{13}-\alpha\Delta_{23}],$ \\
	$\nabla_4 = [\Delta_{14}-\Delta_{24}],$ & $\nabla_5 = [\Delta_{31}],$ & $\nabla_6 = [\Delta_{42}].$
	\end{longtable}	
	
Take $\theta=\sum\limits_{i=1}^{6}\alpha_i\nabla_i\in {\rm H^2}({\mathfrak N}_{08}^{\alpha\neq1}).$
	The automorphism group of ${\mathfrak N}_{08}^{\alpha\neq1}$ consists of invertible matrices of the form
	$$\phi_1=
	\begin{pmatrix}
	x &  0  & 0 & 0\\
	0 &  x  & 0 & 0\\
	t &  v  & x^2 & 0\\
    u &  w  & 0 & x^2
	\end{pmatrix}, \quad 
	\phi_2(\alpha \neq 0)=
	\begin{pmatrix}
	0 &  \alpha x  & 0 & 0\\
	x &  0  & 0 & 0\\
	t &  v  & 0 & -\alpha^2x^2\\
    u &  w  & -x^2 & 0
	\end{pmatrix}.
	$$
	Since
	$$
	\phi_1^T\begin{pmatrix}
	0 & \alpha_1   & \alpha_3 & \alpha_4\\
	\alpha_2 & 0 & -\alpha\alpha_3 & -\alpha_4\\
	\alpha_5 &  0  & 0 & 0\\
	0 &  \alpha_6  & 0 & 0
	\end{pmatrix} \phi_1=	\begin{pmatrix}
\alpha^* & \alpha_1^*+\alpha^{**}   & \alpha_3^* & \alpha_4^*\\
\alpha_2^*-\alpha\alpha^* & -\alpha^{**} & -\alpha\alpha_3^* & -\alpha_4^*\\
\alpha_5^* &  0  & 0 & 0\\
0 &  \alpha_6^*  & 0 & 0
	\end{pmatrix},
	$$
	 we have that the action of ${\rm Aut} ({\mathfrak N}_{08}^{\alpha\neq1})$ on the subspace
$\langle \sum\limits_{i=1}^{6}\alpha_i\nabla_i  \rangle$
is given by
$\langle \sum\limits_{i=1}^{6}\alpha_i^{*}\nabla_i\rangle,$
where
\begin{longtable}{lcllcllcl}
$\alpha^*_1$&$=$&$x \left(x \alpha _1+(v-v \alpha ) \alpha _3+(u+w) \alpha _6\right),$& 
$\alpha^*_3$&$=$&$x^3 \alpha _3,$ &
$\alpha_5^*$&$=$&$x^3 \alpha _5,$\\

$\alpha^*_2$&$=$&$x(x \alpha _2+u (\alpha-1 ) \alpha _4+(v+t \alpha ) \alpha _5),$&
$\alpha_4^*$&$=$&$x^3 \alpha _4,$&
$\alpha_6^*$&$=$&$x^3 \alpha _6.$
\end{longtable}

 We are interested only in the cases with 
 $(\alpha_3,\alpha_5)\neq (0,0),$  $(\alpha_4,\alpha_6)\neq (0,0).$ 

\begin{enumerate}
    \item $\alpha_5 = \alpha_6= 0,$ then $\alpha_3\alpha_4\neq 0,$ and choosing $x=1,\ u=\frac{\alpha _2}{(1-\alpha)\alpha_4},$ $v=-\frac{\alpha_1}{(1-\alpha)\alpha_3},$ we have $\alpha_1^*=\alpha_2^*=0$ and  obtain the family of representatives  $\left\langle\nabla_3+ \beta \nabla_4 \right\rangle_{\beta \neq0}.$
\item $(\alpha_5,\alpha_6)\neq(0,0),$ $\alpha\neq 0,$ then with an action of a suitable $\phi_2$, we can suppose $\alpha_5\neq0$ and choosing $v=-\frac{x \alpha _2+u (\alpha -1) \alpha _4+t \alpha  \alpha _5}{\alpha _5},$ we can suppose $\alpha_2^*=0.$ Now we consider following subcases: 
\begin{enumerate}
    \item $\alpha_3=\alpha_6=\alpha_1=0,$ then we have the family of representatives $\langle\beta\nabla_4+\nabla_5\rangle_{\beta\neq0};$
    \item $\alpha_3=\alpha_6=0,\ \alpha_1\neq 0,$ then choosing $x=\frac{\alpha_1}{\alpha_5},$ we have the family of representatives $\langle\nabla_1+\beta\nabla_4+\nabla_5\rangle_{\beta\neq0};$
    \item $\alpha_3=0,\ \alpha_6\neq0,$ then choosing $x=1,\ u=-\frac{\alpha_1}{\alpha_6},\ w=0,$ we have the family of representatives $\langle\gamma\nabla_4+\nabla_5+\beta\nabla_6 \rangle_{\beta\neq0};$
    \item $\alpha_3\neq0,$ then choosing $x=1,\ t=-\frac{(\alpha -1 ) \alpha _2 \alpha _3+\alpha _1 \alpha _5}{(\alpha-1 ) \alpha  \alpha _3 \alpha _5},\  u=0,\ w=0, $ we have the family of representatives $\langle\gamma\nabla_3+\delta\nabla_4+\nabla_5+\beta\nabla_6 \rangle_{\gamma\neq0,\ (\beta,\delta)\neq(0,0)}.$
\end{enumerate}
\item $(\alpha_5, \alpha_6) \neq (0,0),$ $\alpha= 0.$ If $\alpha_5\neq0,$ then we obtain the previous cases. Thus, we consider the case of $\alpha_5=0.$ Then $\alpha_3\alpha_6 \neq 0$ and choosing $v=-\frac{x \alpha _1+(u+w) \alpha _6}{\alpha _3},$ we can suppose $\alpha_1^*=0.$ Now we consider following
subcases: 
    \begin{enumerate}
    \item $\alpha_4=0,\ \alpha_2=0,$ then  we have the family of representatives $\langle\beta\nabla_3+\nabla_6 \rangle_{\beta\neq0},$
    \item $\alpha_4=0,$ $\alpha_2\neq0,$ then choosing $x=\frac{\alpha_2}{\alpha_6},$ we have the family of representatives $\langle\nabla_2+\beta\nabla_3+\nabla_6 \rangle_{\beta\neq0},$
    \item $\alpha_4\neq0,$ then choosing $x=1,\ u=\frac{\alpha_2}{\alpha_4},$ we obtain $\alpha_2^*=0$ and obtain the family of representatives $\langle\beta\nabla_3+\gamma\nabla_4+\nabla_6\rangle_{\alpha=0, \beta\neq0,  \gamma \neq 0}.$
    \end{enumerate}
\end{enumerate}

\subsubsection{Central extensions of ${\mathfrak N}_{08}^{1}$}
	Let us use the following notations:
	\begin{longtable}{lllllll} 
	$\nabla_1 = [\Delta_{12}],$ & $\nabla_2 = [\Delta_{21}],$ & $\nabla_3 = [\Delta_{13}-\Delta_{23}],$&
	$\nabla_4 = [\Delta_{14}-\Delta_{24}],$ \\ 
 $\nabla_5 = [\Delta_{31}],$ & $\nabla_6 = [\Delta_{42}],$ & $\nabla_7 = [\Delta_{32}+\Delta_{41}].$
	\end{longtable}	
	
Take $\theta=\sum\limits_{i=1}^{7}\alpha_i\nabla_i\in {\rm H^2}({\mathfrak N}_{08}^{1}).$
	The automorphism group of ${\mathfrak N}_{08}^{1}$ consists of invertible matrices of the form
	$$\phi=
	\begin{pmatrix}
	x &  y  & 0 & 0\\
	x+y-z &  z  & 0 & 0\\
	t &  v  & x(z-y) & y(z-y)\\
    u &  w  & (x+y-z)(z-y) & z(z-y)
	\end{pmatrix}.
	$$
	Since
	$$
	\phi^T\begin{pmatrix}
	0 &  \alpha_1   & \alpha_3 & \alpha_4\\
	\alpha_2 & 0 & -\alpha_3 & -\alpha_4\\
	\alpha_5 &  \alpha_7  & 0 & 0\\
	\alpha_7 &  \alpha_6  & 0 & 0
	\end{pmatrix} \phi=	\begin{pmatrix}
\alpha^* & \alpha_1^*+\alpha^{**} & \alpha_3^* & \alpha_4^*\\
\alpha_2^*-\alpha^* & -\alpha^{**} & -\alpha_3^*& -\alpha_4^*\\
\alpha_5^* &  \alpha_7^*  & 0 & 0\\
\alpha_7^* &  \alpha_6^*  & 0 & 0
	\end{pmatrix},
	$$
we have that the action of ${\rm Aut} ({\mathfrak N}_{08}^{1})$ on the subspace
$\langle \sum\limits_{i=1}^{7}\alpha_i\nabla_i  \rangle$
is given by
$\langle \sum\limits_{i=1}^{7}\alpha_i^{*}\nabla_i\rangle,$
where
\begin{longtable}{lcl}
$\alpha^*_1$&$=$&$(x+y) z \alpha _1+y (x+y) \alpha _2+y(t+v)\alpha _5+z(u+w)\alpha _6+  (y(u+w)+z(t+v))\alpha _7,$\\
$\alpha^*_2$&$=$&$(x+y) (x+y-z) \alpha _1+x (x+y) \alpha _2+x(t+v)\alpha _5+(x+y-z)(u+w)\alpha_6$\\ &&\multicolumn{1}{r}{$+((x+y-z)(t+v)+x(u+w))\alpha_7,$}\\
$\alpha^*_3$&$=$&$(y-z)^2(x \alpha _3+(x+y-z) \alpha _4),$ \\
$\alpha_4^*$&$=$&$(y-z)^2(y \alpha _3+z \alpha_4),$\\
$\alpha_5^*$&$=$&$(z-y)(x^2 \alpha _5+(x+y-z)((x+y-z) \alpha _6+2 x \alpha _7)),$\\
$\alpha_6^*$&$=$&$(z-y)(y^2 \alpha _5+z(z \alpha _6+2 y \alpha _7)),$\\
$\alpha_7^*$&$=$&$(z-y)(x y \alpha _5+(x+y-z) z \alpha _6+(y (y-z)+x (y+z)) \alpha _7).$\\
\end{longtable}

 We are interested only in the cases with 
 \begin{center}
$(\alpha_3,\alpha_5,\alpha_7)\neq (0,0,0),$  $(\alpha_4,\alpha_6,\alpha_7)\neq (0,0,0).$ 
 \end{center} 

\begin{enumerate}
    \item $(\alpha_5,\alpha_6,\alpha_7)=(0,0,0),$ then $\alpha_3\neq0,\ \alpha_4\neq0.$ If $\alpha_3\neq-\alpha_4,$ then choosing $z=-\frac{y \alpha _3}{\alpha _4},$ we obtain that $\alpha_4^*=0,$ which implies $(\alpha_4^*,\alpha_6^*,\alpha_7^*) = (0,0,0).$ Thus, we have that $\alpha_3=-\alpha_4.$
    \begin{enumerate}
        \item $(\alpha_1,\alpha_2)=(0,0),$ then we have the representative $\left\langle\nabla_3-\nabla_4\right\rangle;$
        \item $(\alpha_1,\alpha_2)\neq(0,0),$ without loss of generality, we can suppose $\alpha_1\neq0.$
        \begin{enumerate}
            \item $\alpha_1=-\alpha_2,$ then choosing $x=\frac{\alpha_3}{\alpha_1},$ $y=0,$ $z=1,$ we have the representative $\left\langle\nabla_1-\nabla_2+\nabla_3-\nabla_4\right\rangle;$
            \item $\alpha_1\neq-\alpha_2,$ then choosing $x=\frac{\alpha _1^3}{\left(\alpha _1+\alpha _2\right){}^2 \alpha _3},\ y=0,\ z= \frac{\alpha _1^2}{\left(\alpha _1+\alpha _2\right) \alpha _3},$ we have the representative $\left\langle\nabla_1+\nabla_3-\nabla_4\right\rangle.$
        \end{enumerate}
    \end{enumerate}
    \item $(\alpha_5,\alpha_6,\alpha_7)\neq(0,0,0),$ then without loss of generality we can assume $\alpha_5\neq0$ and  consider following subcases:
    \begin{enumerate}
        \item $\alpha_6\alpha_5 =\alpha_7^2, \alpha_7=-\alpha_5,\ \alpha_4=-\alpha_3,\ \alpha_2 =-\alpha_1,$ then taking
        $x=1, \ y=u=v=w=0,$ $t=\frac{\alpha _1}{\alpha _5},$ we have the family of representatives $\left\langle\beta\nabla_3 -\beta\nabla_4+\nabla_5+\nabla_6-\nabla_7 \right\rangle;$
        \item $\alpha_6\alpha_5 =\alpha_7^2,  \alpha_7=-\alpha_5,\  \alpha_4=-\alpha_3,\ \alpha_1\neq- \alpha_2,$ then  taking \begin{center}$x=z=\frac{(\alpha _1+\alpha _2)}{\alpha _5},$ 
        $y=u=v=w=0,$  $t=\frac{\alpha_1(\alpha _1+\alpha _2)}{\alpha _5^2},$ \end{center} we have the family of representatives $\left\langle\nabla_2 + \beta\nabla_3 -\beta\nabla_4+\nabla_5+\nabla_6-\nabla_7 \right\rangle;$
        \item $\alpha_6\alpha_5 =\alpha_7^2, \  \alpha_7=-\alpha_5,\ \alpha_3\neq-\alpha_4,$ then we can suppose $\alpha_3 \neq 0$ and choosing 
     $u=v=w=0,$ $y=-\frac{z \alpha _4}{\alpha _3},$ $t=\frac{\left(x \alpha _3-z \alpha _4\right) \left(\alpha _1 \alpha _3-\alpha _2 \alpha _4\right)}{\alpha _3 \left(\alpha _3+\alpha _4\right) \alpha _5},$ we can suppose  $\alpha_1=\alpha_4=0.$
            \begin{enumerate}
            \item if $\alpha_2=0,$ then choosing $x=1,$ $z=\frac{\alpha_3}{\alpha_5},$ we have the representative $\left\langle \nabla_3+\nabla_5+\nabla_6-\nabla_7 \right\rangle;$
            \item if $\alpha_2\neq0,$ then choosing $x=\frac{\alpha_2\alpha _5^2}{\alpha_3^3},$ $z=\frac{\alpha_2 \alpha_5}{\alpha_3^2},$ we have the representative 
            $\left\langle \nabla_2+ \nabla_3+\nabla_5+\nabla_6-\nabla_7 \right\rangle.$ 
            \end{enumerate}
        \item $\alpha_6\alpha_5 =\alpha_7^2, \alpha_5\neq-\alpha_7,$ then choosing 
        \begin{center}$y=-\frac {z\alpha_7}{\alpha_5},$ $t=u=w=0,$ $v=\frac{\alpha _1 \left((z-x) \alpha _5+z \alpha _7\right)-x \alpha _2 \alpha _5}{\alpha _5(\alpha _5+\alpha _7)},$
        \end{center} we can suppose $\alpha_2=\alpha_6=\alpha_7=0.$ Since $(\alpha_4,\alpha_6,\alpha_7)\neq (0,0,0),$ we have that $\alpha_4\neq 0.$
            \begin{enumerate}
            \item $\alpha_1=0,$ then choosing $x=\sqrt{\frac{\alpha_4}{\alpha_5}},\ z=1,$  we have the family of representatives $\left\langle \beta \nabla_3+ \nabla_4+\nabla_5\right\rangle;$
            \item $\alpha_1\neq0,$ then choosing $x=\frac{\alpha_1}{\alpha_5},$  $z=\frac{\alpha_1}{\sqrt{\alpha_4\alpha_5}},$ we have the family of representatives $\left\langle \nabla_1+\beta\nabla_3+\nabla_4+\nabla_5\right\rangle.$
            \end{enumerate}
        \item $\alpha_6\alpha_5 \neq \alpha_7^2,$ then choosing suitable value of $z$ and $y$ such that $y\neq z,$ we can suppose $\alpha_6=0$ and $\alpha_7^*\neq 0.$ Then choosing $t=-\frac{x \alpha _1}{\alpha _7},\ y=u=v=0,\ w=\frac{x \left(\alpha _1 \alpha _5-\alpha _2 \alpha _7\right)}{\alpha _7^2},$ we have $\alpha_1^*=\alpha_2^*=0.$ 

\begin{enumerate}
            \item
if $\alpha_5\neq-2\alpha_7,$ then choosing $x=1,$ $z=\frac{\alpha _5+2 \alpha _7}{3 \alpha _7},$ we have the family of representatives $\left\langle\beta\nabla_3+\gamma\nabla_4+\nabla_5+\nabla_7\right\rangle;$
           
    \item    
        if $\alpha_5=-2\alpha_7,$ then then we have the family of representatives $\left\langle\beta\nabla_3-2\nabla_5+\nabla_7\right\rangle$
           and $\left\langle\beta\nabla_3+\nabla_4-2\nabla_5 +\nabla_7\right\rangle$ depending on whether $\alpha_4=0$ or not.
           
     \end{enumerate}
    \end{enumerate}
\end{enumerate}

Summarizing all cases of the central extension of the algebra ${\mathfrak N}_{08}^{\alpha},$ we have the following distinct orbits,

in case of $\alpha\neq1$: 
\begin{center} 
$\langle\nabla_3+\beta\nabla_4\rangle_{\beta\neq0},$ 
$\langle\nabla_1+\beta\nabla_4+\nabla_5\rangle_{\beta\neq0},$  $\langle\delta\nabla_3+\gamma\nabla_4+\nabla_5+\beta\nabla_6\rangle_{(\beta,\gamma)\neq(0,0)},$
\end{center}

in case of $\alpha=0$: 
\begin{center} 
$\langle\nabla_2+\beta\nabla_3+\nabla_6\rangle_{\beta\neq0}$, 
$\langle\beta\nabla_3+\gamma\nabla_4+\nabla_6\rangle_{\beta\neq0},$
\end{center}

in case of $\alpha=1$:
\begin{center} 
$\langle \nabla_3-\nabla_4\rangle,$ 
$\langle\nabla_1+\nabla_3-\nabla_4\rangle,$  
$\langle\nabla_1-\nabla_2+\nabla_3-\nabla_4\rangle,$ 
$\langle\beta\nabla_3-\beta\nabla_4+\nabla_5+\nabla_6-\nabla_7\rangle,$ 
$\langle\nabla_2+\beta\nabla_3-\beta\nabla_4+\nabla_5+\nabla_6-\nabla_7\rangle,$  $\langle\nabla_3+\nabla_5+\nabla_6-\nabla_7\rangle,$ 
$\langle\nabla_2+\nabla_3+\nabla_5+\nabla_6-\nabla_7\rangle,$ 
$\langle\beta\nabla_3+\nabla_4+\nabla_5\rangle,$ $\langle\nabla_1+\beta\nabla_3+\nabla_4+\nabla_5\rangle,$ 
$\langle\beta\nabla_3+\gamma\nabla_4+\nabla_5+\nabla_7\rangle,$  $\langle\beta\nabla_3-2\nabla_5+\nabla_7\rangle,$ 
$\langle\beta\nabla_3+\nabla_4-2\nabla_5+\nabla_7\rangle,$ 
\end{center}
which gives the following new algebras (see section \ref{secteoA}):

\begin{center}
${\rm B}_{60}^{\alpha \neq 1, \beta\neq0},$
${\rm B}_{61}^{\alpha \neq 1, \beta\neq0},$ 
${\rm B}_{62}^{\alpha \neq 1, (\beta,\gamma)\neq(0,0), \delta},$ 
${\rm B}_{63}^{\beta\neq0},$ 
${\rm B}_{64}^{\beta \neq 0,\gamma},$ 
${\rm B}_{65},$ 

${\rm B}_{66},$ 
${\rm B}_{67},$ 
${\rm B}_{68}^{\beta},$ 
${\rm B}_{69}^{\beta},$ 
${\rm B}_{70},$ 
${\rm B}_{71},$ 
${\rm B}_{72}^{\beta},$ 
${\rm B}_{73}^{\beta},$ 
${\rm B}_{74}^{\beta,\gamma},$ 
${\rm B}_{75}^{\beta},$ 
${\rm B}_{76}^{\beta}.$
\end{center}

\subsubsection{Central extensions of ${\mathfrak N}_{12}$}
	Let us use the following notations:
	\begin{longtable}{lllllll} 
	$\nabla_1 = [\Delta_{11}],$ & $\nabla_2 = [\Delta_{13}],$ & 
	$\nabla_3 = [\Delta_{22}],$ \\
	$\nabla_4 = [\Delta_{24}],$ & 
	$\nabla_5 = [\Delta_{32}],$ & 
	$\nabla_6 = [\Delta_{41}].$ 
	\end{longtable}	
	
Take $\theta=\sum\limits_{i=1}^{6}\alpha_i\nabla_i\in {\rm H^2}({\mathfrak N}_{12}).$
	The automorphism group of ${\mathfrak N}_{12}$ consists of invertible matrices of the form
		$$\phi_1=
	\begin{pmatrix}
	x &  0  & 0 & 0\\
	0 &  y  & 0 & 0\\
	z &  v  & xy & 0\\
    u &  t  & 0 & xy
	\end{pmatrix}, \quad 
	\phi_2=
	\begin{pmatrix}
	0 &  x  & 0 & 0\\
	y &  0  & 0 & 0\\
	z &  v  & 0 & x y\\
    u &  t  & x y & 0
	\end{pmatrix}.
	$$
	Since
	$$
	\phi^T_1\begin{pmatrix}
	\alpha_1 & 0 & \alpha_2 & 0\\
	0 & \alpha_3 & 0 & \alpha_4\\
	0 &  \alpha_5  & 0 & 0\\
	\alpha_6 &  0  & 0 & 0
	\end{pmatrix} \phi_1=\begin{pmatrix}
    \alpha_1^* & \alpha^* & \alpha_2^* & 0\\
	\alpha^{**} & \alpha_3^* & 0 & \alpha_4^*\\
	0 &  \alpha_5^*  & 0 & 0\\
	\alpha_6^* &  0  & 0 & 0
	\end{pmatrix},
	$$
	 we have that the action of ${\rm Aut} ({\mathfrak N}_{12})$ on the subspace
$\langle \sum\limits_{i=1}^{6}\alpha_i\nabla_i  \rangle$
is given by
$\langle \sum\limits_{i=1}^{6}\alpha_i^{*}\nabla_i\rangle,$
where for $\phi_1:$
\begin{longtable}{lcllcllcl}
$\alpha^*_1$&$=$&$x \left(x \alpha _1+z \alpha _2+t \alpha _6\right)$ & $\alpha^*_3$&$=$&$y \left(y \alpha _3+v \alpha _4+u \alpha _5\right),$ &
$\alpha_5^*$&$=$&$x y^2\alpha_5,$\\

$\alpha^*_2$&$=$&$x^2 y \alpha _2$&
$\alpha_4^*$&$=$&$x y^2 \alpha _4,$&
$\alpha_6^*$&$=$&$x^2 y\alpha_6,$
\end{longtable}
for $\phi_2:$
\begin{longtable}{lcllcllcl}
$\alpha^*_1$&$=$&$y \left(y \alpha _3+u \alpha _4+z \alpha _5\right)$ &
$\alpha^*_3$&$=$&$x \left(x \alpha _1+v \alpha _2+t \alpha _6\right),$&
$\alpha_5^*$&$=$&$x^2 y \alpha _6,$\\

$\alpha^*_2$&$=$&$x y^2 \alpha _4$ &
$\alpha_4^*$&$=$&$x^2 y \alpha _2,$&
$\alpha_6^*$&$=$&$x y^2 \alpha _5,$
\end{longtable}

 We are interested only in the cases with 
 \begin{center}
$(\alpha_2,\alpha_5)\neq (0,0),$  $(\alpha_4,\alpha_6)\neq (0,0).$ 
 \end{center} 
\begin{enumerate}
    \item $(\alpha_2,\alpha_4)=(0,0),$ then $\alpha_5\neq0,\ \alpha_6\neq0$ and choosing $x=1, \ y=\frac{\alpha _6}{\alpha _5},$ $t=-\frac{\alpha _1}{\alpha _6},$ $u=-\frac{\alpha_3 \alpha _6}{\alpha _5^2},$ we have the representative $\left\langle\nabla_5+\nabla_6\right\rangle;$
    \item $(\alpha_2,\alpha_4)\neq(0,0),$ then without loss of generality, we can suppose $\alpha_4\neq0$ and choosing $v=-\frac{y\alpha_3+u \alpha_5}{\alpha_4},$ we have $\alpha_3^*=0.$
    \begin{enumerate}
        \item $\alpha_6=\alpha_2=\alpha_1=0,$ then we have the family of representatives $\left\langle\nabla_4+\alpha\nabla_5\right\rangle_{\alpha\neq0};$
        \item $\alpha_6=\alpha_2=0,\ \alpha_1\neq0,$ then choosing $x=\frac{\alpha_4}{\alpha_1},\ y=1,$ we have the family of representatives $\left\langle\nabla_1+\nabla_4+\alpha\nabla_5\right\rangle_{\alpha\neq0};$
        \item $\alpha_6=0,\ \alpha_2\neq0,$ then choosing $x=1,\ y=\frac{\alpha_2}{\alpha_4}, \  z=-\frac{\alpha_1}{\alpha_2},$ we have the family of representatives $\left\langle\nabla_2+\nabla_4+\alpha\nabla_5\right\rangle;$
        \item $\alpha_6\neq0,$ then  choosing $x=1,\ y=\frac{\alpha_6}{\alpha_4},\ z=0, \ t=-\frac{\alpha_1}{\alpha_6},$ we have the family of representatives $\left\langle\beta\nabla_2+\nabla_4+\alpha\nabla_5 +\nabla_6\right\rangle_{(\alpha,\beta)\neq(0,0)}.$
    \end{enumerate}
\end{enumerate}

Summarizing all cases, we have the following distinct orbits 

\begin{center} 
$\langle\nabla_5+\nabla_6\rangle,$ 
$\langle\nabla_4+\alpha\nabla_5\rangle_{\alpha\neq0},$ $\langle\nabla_1+\nabla_4+\alpha\nabla_5\rangle_{\alpha\neq0},$
$\langle\nabla_2+\nabla_4+\alpha\nabla_5\rangle,$ 
$\langle\beta\nabla_2+\nabla_4+\alpha\nabla_5+\nabla_6\rangle_{(\alpha,\beta)\neq(0,0)}^{O(\alpha,\beta)\simeq O(\beta^{-1},  {\alpha^{-1}})},$
\end{center}
which gives the following new algebras (see section \ref{secteoA}):

\begin{center}
${\rm B}_{77},$
${\rm B}_{78}^{\alpha\neq0},$
${\rm B}_{79}^{\alpha\neq0},$
${\rm B}_{80},$
${\rm B}_{81}^{(\alpha,\beta)\neq(0,0)}.$
\end{center}

\subsubsection{Central extensions of ${\mathfrak N}_{13}$}
	Let us use the following notations:
	\begin{longtable}{lllllll} 
	$\nabla_1 = [\Delta_{21}],$ & $\nabla_2 = [\Delta_{22}],$ & 
	$\nabla_3 = [\Delta_{14}+\Delta_{23}],$ \\
	$\nabla_4 = [\Delta_{13}-2\Delta_{14}-\Delta_{24}],$ & 
	$\nabla_5 = [\Delta_{32}-\Delta_{41}],$ & 
	$\nabla_6 = [\Delta_{31}-2\Delta_{32}+\Delta_{42}].$ 
	\end{longtable}	
	
Take $\theta=\sum\limits_{i=1}^{6}\alpha_i\nabla_i\in {\rm H^2}({\mathfrak N}_{13}).$
	The automorphism group of ${\mathfrak N}_{13}$ consists of invertible matrices of the form
	$$\phi_1=
	\begin{pmatrix}
	x &  0  & 0 & 0\\
	0 & x  & 0 & 0\\
	z &  u  & x^2 & 0\\
    t &  v  & 0 & x^2
	\end{pmatrix}, \quad 
	\phi_2=
	\begin{pmatrix}
	0 &  x  & 0 & 0\\
	x &  0  & 0 & 0\\
	z &  u  & -x^2 & 2x^2\\
    t &  v  & 0 & x^2
	\end{pmatrix}.
	$$
	Since
	$$
	\phi^T_1\begin{pmatrix}
	0 & 0 & \alpha_4 & \alpha_3-2\alpha_4\\
	\alpha_1 & \alpha_2 & \alpha_3 & -\alpha_4\\
	\alpha_6 &  \alpha_5-\alpha_6  & 0 & 0\\
	-\alpha_5 &  \alpha_6  & 0 & 0
	\end{pmatrix} \phi_1=\begin{pmatrix}
    \alpha^* & \alpha^{**} & \alpha_4^*& \alpha_3^*-2\alpha_4^*\\
	\alpha_1^*-\alpha^{**} & \alpha_2^*+\alpha^*+2\alpha^{**} & \alpha_3^* & -\alpha_4^*\\
	\alpha_6^* &  \alpha_5^*-\alpha_6^*  & 0 & 0\\
	-\alpha_5^* &  \alpha_6^*  & 0 & 0
	\end{pmatrix},
	$$
	 we have that the action of ${\rm Aut} ({\mathfrak N}_{13})$ on the subspace
$\langle \sum\limits_{i=1}^{6}\alpha_i\nabla_i  \rangle$
is given by
$\langle \sum\limits_{i=1}^{6}\alpha_i^{*}\nabla_i\rangle,$
where
for $\phi_1$: \begin{longtable}{lcl}
$\alpha^*_1$&$=$&$x(x\alpha_1+(v+z)\alpha_3-(t-u+2v)\alpha _4-(v-z)\alpha_5+(t+u-z)\alpha_6),$ \\
$\alpha^*_2$&$=$&$x(x\alpha _2+(-t+u-2 v) \alpha _3+(2t-2u+3v-z)\alpha _4+$\\
&&\multicolumn{1}{r}{$(t+u-2z)\alpha_5-(2t +u -v-z)\alpha _6),$}\\
$\alpha^*_3$&$=$&$x^3 \alpha _3,$ \\
$\alpha_4^*$&$=$&$x^3 \alpha _4,$\\
$\alpha_5^*$&$=$&$x^3 \alpha _5,$\\
$\alpha_6^*$&$=$&$x^3 \alpha _6.$\\
\end{longtable}

for $\phi_2$: \begin{longtable}{lcl}
$\alpha^*_1$&$=$&$x(x\alpha_1+(t+u)\alpha_3-(2t+v-z)\alpha_4-(t-u)\alpha _5-(u-v-z)\alpha_6),$ \\
$\alpha^*_2$&$=$&$-x(2x\alpha_1+x\alpha_2+(2u-v+z)\alpha_3-(t+u)\alpha _4-(2t-v-z)\alpha _5+ (t-u+z)\alpha_6),$\\
$\alpha^*_3$&$=$&$-x^3 \alpha _4,$ \\
$\alpha_4^*$&$=$&$-x^3 \alpha _3,$\\
$\alpha_5^*$&$=$&$-x^3(2\alpha_5-\alpha_6),$\\
$\alpha_6^*$&$=$&$-x^3(\alpha_5-2\alpha_6).$\\
\end{longtable}

 We are interested only in the cases with 
 \begin{center}
$(\alpha_3,\alpha_4,\alpha_5,\alpha_6)\neq (0,0,0,0).$ 
 \end{center} 
\begin{enumerate}
    \item $(\alpha_5,\alpha_6)=(0,0),$ then without loss of generality, we can suppose $\alpha_3 \neq 0$. Let us consider the following subcases:
    \begin{enumerate}
        \item $\alpha_3=\alpha_4,\ \alpha_1=-\alpha_2,$ then choosing $x=1,\ u=v=z=0,$ $t=\frac{\alpha _1}{\alpha _3},$ we have the representative  $\left\langle\nabla_3+\nabla_4\right\rangle;$
        \item $\alpha_3=\alpha_4,\ \alpha_1\neq-\alpha_2,$ then choosing $x=\frac{\alpha_1+\alpha_2}{\alpha_3},$ $u=v=z=0,$ $t=\frac{(\alpha_1+\alpha_2)\alpha_1}{\alpha_3^2},$  we have the representative  $\left\langle\nabla_1+\nabla_3+\nabla_4\right\rangle;$
        \item $\alpha_3\neq\alpha_4,$ then choosing $x=1, $ $u=v=0,$ $ z=\frac{\alpha _1 \left(2 \alpha _4-\alpha _3\right)+\alpha _2 \alpha _4}{\left(\alpha _3-\alpha _4\right){}^2},$ $t=\frac{\alpha _2 \alpha _3+\alpha _1 \alpha _4}{\left(\alpha _3-\alpha _4\right){}^2},$ we have the family of representatives  $\left\langle\nabla_3+\alpha\nabla_4\right\rangle_{\alpha \neq1};$
    \end{enumerate}
    \item $(\alpha_5,\alpha_6)\neq (0,0),$ then without loss of generality, we can suppose $\alpha_6\neq0.$ Let us consider the following subcases:  
    \begin{enumerate}
        \item $\alpha_4=\alpha_6,\ \alpha_3=\alpha_5,$ then choosing 
        \begin{center} $x=1,$  $z=0,$ 
        $u=\frac{ \alpha _1 \left(\alpha _3-2 \alpha _6\right)-\alpha _2 \alpha _6}{\alpha _6^2},$
        $v=\frac{\alpha _1 \left(2 \alpha _3-3 \alpha _6\right)-2 \alpha _2 \alpha _6}{2 \alpha _6^2},$ \end{center}  we have the family of representatives  $\left\langle\alpha \nabla_3+ \nabla_4+ \alpha\nabla_5+\nabla_6\right\rangle;$
        \item $\alpha_4=\alpha_6,\ \alpha_3\neq\alpha_5,$ then choosing $x=1,\ z=v=0,\ t=-\frac{x \left(\alpha _1 \left(\alpha _3+\alpha _5-3 \alpha _6\right)-2 \alpha _2 \alpha _6\right)}{2 \left(\alpha _3-\alpha _5\right) \alpha _6}$ $u=-\frac{\alpha _1}{2 \alpha _6},$  we have the family of representatives  $\left\langle\alpha \nabla_3+ \nabla_4+ \beta\nabla_5+\nabla_6\right\rangle_{\beta\neq\alpha};$
        \item $\alpha_4-\alpha_6\neq0,$ then choosing $t=\frac{x \alpha _1+(v+z) \alpha _3+u \alpha _4-2 v \alpha _4-v \alpha _5+z \alpha _5+u \alpha _6-z \alpha _6}{\alpha _4-\alpha _6},$ we can suppose $\alpha_1^*=0,$ and consider following subcases:
        \begin{enumerate}
            \item $\alpha_6(2\alpha _3+\alpha_6)=\alpha_4(2 \alpha _5+\alpha _6).$ 
            \begin{enumerate}
                \item $\alpha_6=2\alpha_5,$ then choosing $x=1,\ z=0,\ v=-\frac{\alpha _2}{2 \alpha _5},$ we have the family of representatives $\left\langle(\alpha-\frac12)\nabla_3+\alpha\nabla_4+\frac12\nabla_5+\nabla_6\right\rangle_{\alpha\neq1};$
                \item $\alpha_6\neq2\alpha_5,\ \alpha _4=0, \ 4 \alpha _5^2-4 \alpha_5\alpha_6+5\alpha_6^2=0, \alpha_2=0,$ then we have the representative $\left\langle-\frac12\nabla_3+(\frac12\pm i)\nabla_5+\nabla_6\right\rangle;$
                \item $\alpha_6\neq2\alpha_5,\ \alpha _4=0, \ 4 \alpha _5^2-4 \alpha_5\alpha_6+5\alpha_6^2=0, \  \alpha_2\neq0,$ then choosing $x=\frac{\alpha_2}{\alpha_6},$ we have the representative $\left\langle\nabla_2-\frac12\nabla_3+(\frac12\pm i)\nabla_5+\nabla_6\right\rangle;$
  \item $\alpha_6\neq2\alpha_5,\ \alpha _4(\alpha_6-2\alpha_5)^2=\alpha_6(4 \alpha _5 \alpha _6-4\alpha _5^2-5 \alpha _6^2), \ 4 \alpha _5^2\neq4 \alpha_5\alpha_6+5\alpha_6^2,$ then choosing $x=1,\ v=-\frac{2\alpha _2 \alpha _6}{4 \alpha _5^2-4 \alpha _5 \alpha _6+5 \alpha _6^2},$ we have the family of representatives 
                \begin{center}$\left\langle-\frac{(3+\alpha+4\alpha^3)}{(1-2\alpha)^2}\nabla_3- \frac{(5-4\alpha+4\alpha^2)}{(1-2\alpha)^2}\nabla_4+ \alpha\nabla_5+\nabla_6\right\rangle_{\alpha\neq\frac12,\frac{1}{2}\pm i};$\end{center}
                 \item $\alpha_6\neq2\alpha_5,\ \alpha _4(\alpha_6-2\alpha_5)^2+\alpha_6(4\alpha _5^2-4 \alpha _5 \alpha _6+5 \alpha _6^2)\neq0,$ then choosing $x=1,$ $z=\frac{4\alpha _2 \alpha _6^2}{\alpha _4 \left(2 \alpha _5-\alpha _6\right){}^2+\alpha _6 \left(4 \alpha _5^2-4 \alpha _5 \alpha _6+5 \alpha _6^2\right)},$ $v=0,$ we have the family of representatives \begin{center}$\left\langle\frac12(\beta(1+2\alpha)+\alpha)\nabla_3+\beta\nabla_4+  \alpha\nabla_5+\nabla_6\right\rangle_{\alpha\neq\frac{1}{2},\ \beta\neq1, -\frac{(5-4\alpha+4\alpha^2)}{(1-2\alpha)^2}}.$\end{center}
            \end{enumerate}
            \item $\alpha_6(2\alpha _3+\alpha_6)\neq\alpha_4(2 \alpha _5+\alpha _6),$ then choosing  
            $u=\frac{\alpha _2 \left(\alpha _4-\alpha _6\right)}{\alpha _6 \left(2 \alpha _3+\alpha _6\right)-\alpha _4 \left(2 \alpha _5+\alpha _6\right)},$ $x=1,$  $z=v=0,$  we have the family of representatives \begin{center}$\left\langle\gamma\nabla_3+\beta\nabla_4+\alpha\nabla_5+\nabla_6\right\rangle_{\gamma\neq\frac{\beta(2\alpha+1)-1}{2},\ \beta\neq1}.$\end{center}
        \end{enumerate}
    \end{enumerate}
\end{enumerate}

Summarizing all cases, we have the following distinct orbits 
\begin{center} 
$\langle\nabla_4\rangle,$  $\langle\nabla_3+\alpha\nabla_4\rangle^{O(\alpha)\simeq O({\alpha^{-1}})},$  $\langle\nabla_1+\nabla_3+\nabla_4\rangle,$ 
$\langle\nabla_2-\frac12\nabla_3+(\frac12 \pm i)\nabla_5+\nabla_6\rangle,$
$\langle\gamma\nabla_3+\alpha\nabla_4+\beta\nabla_5+\nabla_6\rangle^{O(\alpha,\beta,\gamma)\simeq O(\frac{\gamma}{\beta-2}, \frac{1-2\beta}{2-\beta},\frac{\alpha}{\beta-2})},$
\end{center}
which gives the following new algebras (see section \ref{secteoA}):

\begin{center}
${\rm B}_{82},$
${\rm B}_{83}^{\alpha},$
${\rm B}_{84},$ 
${\rm B}_{85},$
${\rm B}_{86},$
${\rm B}_{87}^{\alpha,\beta,\gamma}.$
\end{center}

\subsubsection{Central extensions of ${\mathfrak N}_{14}^0$}
	Let us use the following notations:
	\begin{longtable}{lllllll} 
	$\nabla_1 = [\Delta_{11}],$ & $\nabla_2 = [\Delta_{21}],$ &
	$\nabla_3 = [\Delta_{23}],$ &
	$\nabla_4 = [\Delta_{13}+\Delta_{24}],$ \\ 
	$\nabla_5 = [\Delta_{32}],$ &
	$\nabla_6 = [\Delta_{42}],$ & 
	$\nabla_7 = [\Delta_{14}].$
	\end{longtable}	
	
Take $\theta=\sum\limits_{i=1}^{7}\alpha_i\nabla_i\in {\rm H^2}({\mathfrak N}_{14}^0).$
	The automorphism group of ${\mathfrak N}_{14}^0$ consists of invertible matrices of the form
	$$\phi=
	\begin{pmatrix}
	x &  z  & 0 & 0\\
	0 &  y  & 0 & 0\\
	w &  u  & y^2 & 0\\
    t &  v  & yz & xy
	\end{pmatrix}.
	$$
	Since
	$$
	\phi^T\begin{pmatrix}
	\alpha_1 & 0 & \alpha_4 & \alpha_7\\
	\alpha_2 & 0 & \alpha_3 & \alpha_4\\
	0 &  \alpha_5  & 0 & 0\\
	0 &  \alpha_6  & 0 & 0
	\end{pmatrix} \phi=\begin{pmatrix}
    \alpha_1^* & \alpha^* & \alpha_4^* & \alpha_7^*\\
	\alpha_2^* & \alpha^{**} & \alpha_3^* & \alpha_4^*\\
	0 &  \alpha_5^*  & 0 & 0\\
	0 &  \alpha_6^*  & 0 & 0
	\end{pmatrix},
	$$
	 we have that the action of ${\rm Aut} ({\mathfrak N}_{14}^0)$ on the subspace
$\langle \sum\limits_{i=1}^{7}\alpha_i\nabla_i  \rangle$
is given by
$\langle \sum\limits_{i=1}^{7}\alpha_i^{*}\nabla_i\rangle,$
where
\begin{longtable}{lcl}
$\alpha^*_1$&$=$&$x \left(x \alpha _1+w \alpha _4+t \alpha _7\right),$ \\
$\alpha^*_2$&$=$&$x z \alpha _1+x y \alpha _2+w y \alpha _3+(ty+wz)\alpha_4+t z \alpha _7,$\\
$\alpha^*_3$&$=$&$y (y^2\alpha _3+2yz\alpha_4+z^2\alpha _7)),$ \\
$\alpha_4^*$&$=$&$x y \left(y \alpha _4+z \alpha _7\right),$\\
$\alpha_5^*$&$=$&$y^2 (y \alpha _5+z \alpha _6),$\\
$\alpha_6^*$&$=$&$x y^2 \alpha _6,$\\
$\alpha_7^*$&$=$&$x^2 y \alpha _7.$\\
\end{longtable}

We are interested only in the cases with 
 \begin{center}
$(\alpha_3,\alpha_4,\alpha_5)\neq (0,0,0),$ $(\alpha_4,\alpha_6,\alpha_7)\neq (0,0,0).$ 
 \end{center} 

\begin{enumerate}
    \item $\alpha_7\neq0,$ then choosing $z=-\frac{y \alpha _4}{\alpha _7}, \ t= -\frac{x \alpha _1+w \alpha _4}{\alpha _7},$ we have $\alpha_1^*=\alpha_4^*=0.$ Thus, we can suppose $\alpha_1=\alpha_4=0$ and consider following subcases:
    \begin{enumerate}
        \item $\alpha_3\neq0,$ then choosing $x=1, \ y = \sqrt{\frac{\alpha _7}{\alpha _3}}, \ w=-\frac{\alpha_2}{\alpha_3},$ 
        we have the family of representatives  $\left\langle \nabla_3+\beta\nabla_4+\gamma\nabla_6+\nabla_7\right\rangle;$
        \item $\alpha_3=0,$ then $\alpha_5\neq0.$ 
        \begin{enumerate}
        \item $\alpha_2=0,$ then choosing $x=1, \ y = \sqrt{\frac{\alpha _7}{\alpha _5}},$ we have the family of representatives  $\left\langle\nabla_5+\beta\nabla_6+\nabla_7\right\rangle;$
        \item $\alpha_2\neq0,$ then choosing $x=\frac{\alpha_2}{\alpha_7}, \ y=\frac{\alpha _2}{\sqrt{\alpha _5\alpha _7}},$ we have the family of representatives $\left\langle\nabla_2+\nabla_5+\beta\nabla_6+\nabla_7\right\rangle.$
        \end{enumerate}
    \end{enumerate}
    \item $\alpha_7=0,$ $\alpha_4\neq0,$ then choosing $z=-\frac{y \alpha _3}{2\alpha _4},$ $t=\frac{x \left(\alpha _1 \alpha _3-\alpha _2 \alpha _4\right)}{\alpha _4^2},$ $w=-\frac{x \alpha _1}{\alpha _4},$ we have  
    $\alpha_1^*=\alpha_2^*=\alpha_3^*=0$ and consider following subcases: 
\begin{enumerate}
        \item $\alpha_5=0,$ then we have the family of representatives  $\left\langle\nabla_4+\beta\nabla_6\right\rangle;$
        \item  $\alpha_5\neq 0,$ then choosing $x=1, \ y = \frac{\alpha_4}{\alpha_5},$  have the family of representatives  $\left\langle\nabla_4+\nabla_5+\beta\nabla_6\right\rangle.$
        \end{enumerate}
    
    \item $\alpha_7=\alpha_4=0,$ then $\alpha_6\neq0$ and choosing $z=-\frac{y\alpha_5}{\alpha_6},$ we have $\alpha_5^*=0.$ Thus we obtain that $\alpha_5 =0,$ which implies $\alpha_3 \neq 0.$
    Then choosing $w=-\frac{x \alpha _2}{\alpha _3},$ we have $\alpha_2 =0$
    and obtain the representatives  $\left\langle\nabla_3+\nabla_6\right\rangle$ and $\left\langle\nabla_1+ \nabla_3+\nabla_6\right\rangle$ depending on whether $\alpha_1 =0$ or not.
    \end{enumerate}

Summarizing all cases, we have the following distinct orbits 

\begin{center} 
$\left\langle \nabla_3+\beta\nabla_4+\gamma\nabla_6+\nabla_7\right\rangle^{O(\beta, \gamma) \simeq O(\beta, - \gamma)},$ 
$\left\langle\nabla_5+\beta\nabla_6+\nabla_7\right\rangle^{O(\beta) \simeq O(-\beta)},$ 
$\left\langle\nabla_2+\nabla_5+\beta\nabla_6+\nabla_7\right\rangle^{O(\beta) \simeq O(-\beta)},$ 
$\left\langle\nabla_4+\beta\nabla_6\right\rangle,$ 
$\left\langle\nabla_4+\nabla_5+\beta\nabla_6\right\rangle,$
$\left\langle\nabla_3+\nabla_6\right\rangle,$  $\left\langle\nabla_1+ \nabla_3+\nabla_6\right\rangle.$ 
\end{center}

\subsubsection{Central extensions of ${\mathfrak N}_{14}^{\alpha\neq0}$}
	Let us use the following notations:
	\begin{longtable}{lllllll} 
	$\nabla_1 = [\Delta_{11}],$ & $\nabla_2 = [\Delta_{21}],$ & 
	$\nabla_3 = [\Delta_{23}],$ \\
	$\nabla_4 = [\Delta_{13}+\Delta_{24}],$ & 
	$\nabla_5 = [\Delta_{32}],$ &
$\nabla_6=[\alpha\Delta_{31}+\Delta_{42}].$
	\end{longtable}	
	
Take $\theta=\sum\limits_{i=1}^{6}\alpha_i\nabla_i\in {\rm H^2}({\mathfrak N}_{14}^{\alpha\neq0}).$
	The automorphism group of ${\mathfrak N}_{14}^{\alpha\neq0}$ consists of invertible matrices of the form
	$$\phi=
	\begin{pmatrix}
	x &  z  & 0 & 0\\
	0 &  y  & 0 & 0\\
	w &  u  & y^2 & 0\\
    t &  v  & (1+\alpha)yz & xy
	\end{pmatrix}.
	$$
	Since
	$$
	\phi^T\begin{pmatrix}
	\alpha_1 & 0 & \alpha_4 & 0\\
	\alpha_2 & 0 & \alpha_3 & \alpha_4\\
	\alpha\alpha_6 &  \alpha_5  & 0 & 0\\
	0 &  \alpha_6  & 0 & 0
	\end{pmatrix} \phi=\begin{pmatrix}
    \alpha_1^* & \alpha^* & \alpha_4^* & 0\\
	\alpha_2^*+\alpha\alpha^* & \alpha^{**} & \alpha_3^* & \alpha_4^*\\
	\alpha\alpha_6^* &  \alpha_5^*  & 0 & 0\\
	0 &  \alpha_6^*  & 0 & 0
	\end{pmatrix},
	$$
	 we have that the action of ${\rm Aut} ({\mathfrak N}_{14}^{\alpha\neq0})$ on the subspace
$\langle \sum\limits_{i=1}^{6}\alpha_i\nabla_i  \rangle$
is given by
$\langle \sum\limits_{i=1}^{6}\alpha_i^{*}\nabla_i\rangle,$
where
\begin{longtable}{lcl}
$\alpha^*_1$&$=$&$x(x\alpha _1+w(\alpha _4+\alpha \alpha _6)),$ \\
$\alpha^*_2$&$=$&$x z (1-\alpha)\alpha _1+x y \alpha _2+w y \alpha _3+$\\
&&$(ty+wz-ux\alpha)\alpha_4-wy\alpha\alpha_5-\alpha(ty-ux+wz\alpha)\alpha _6,$\\
$\alpha^*_3$&$=$&$y^2(y\alpha_3+z(2+\alpha)\alpha_4),$ \\
$\alpha_4^*$&$=$&$x y^2 \alpha _4,$\\
$\alpha_5^*$&$=$&$y^2(y\alpha _5+z(1+2\alpha)\alpha_6),$\\
$\alpha_6^*$&$=$&$x y^2 \alpha _6.$\\
\end{longtable}

We are interested only in the cases with 
 \begin{center}
$(\alpha_4,\alpha_6)\neq (0,0).$
 \end{center} 

\begin{enumerate}
    \item $\alpha_4=0,$ then $\alpha_6\neq0$ and choosing $w= -\frac{x \alpha _1}{\alpha \alpha _6},$ $t=\frac{x(\alpha\alpha_6(y\alpha_2+u\alpha\alpha _6)-\alpha _1(y\alpha_3 -y\alpha \alpha _5-z \alpha  \alpha _6))}{y \alpha ^2 \alpha _6^2},$ we have $\alpha_1^*=\alpha_2^*=0.$ Thus we can suppose $\alpha_1=\alpha_2=0$ and consider following subcases:  
\begin{enumerate}
    \item $\alpha = - \frac 1 2,$ $\alpha_5=0,$ $\alpha_3=0,$ then we have the representative $\left\langle\nabla_6\right\rangle_{\alpha = - \frac 1 2 };$  
    \item $\alpha = - \frac 1 2,$ $\alpha_5=0,$ $\alpha_3\neq0,$ then choosing $x=\frac{\alpha_3}{\alpha_6}, \ y=1,$  we have the representative $\left\langle\nabla_3+\nabla_6\right\rangle_{\alpha = - \frac 1 2 };$
    \item $\alpha = - \frac 1 2,$ $\alpha_5\neq0,$ then choosing $x=\frac{\alpha_5}{\alpha_6}, \ y=1,$  we have the family of representatives $\left\langle\beta\nabla_3+\nabla_5+\nabla_6\right\rangle_{\alpha = - \frac 1 2 };$
 \item $\alpha \neq - \frac 1 2,$ $\alpha_3=0,$ then choosing $y=1, \ z=-\frac{\alpha_5}{\alpha_6(1+2\alpha)}, $  we have the representative $\left\langle\nabla_6\right\rangle_{\alpha \neq - \frac 1 2 };$
 \item $\alpha \neq - \frac 1 2,$ $\alpha_3\neq0,$ then choosing $x=\frac{\alpha_3}{\alpha_6}, \ y=1, \ z=-\frac{\alpha_5}{\alpha_6(1+2\alpha)}, $  we have the representative $\left\langle\nabla_3+\nabla_6\right\rangle_{\alpha \neq - \frac 1 2 };$
    \end{enumerate}
\item $\alpha_4\neq 0,$ then consider following subcases:
\begin{enumerate}
    \item $\alpha = - 2,$ $\alpha_4=2\alpha_6,$ then choosing $z=\frac{2 y \alpha _1}{3 \alpha _4},$ $t=-\frac{x(\alpha _1\alpha _5+\alpha _2\alpha _6)}{4 \alpha _6^2},$ $u=w=0,$ we can suppose $\alpha_2^*=\alpha_5^*=0$ and consider following subcases:
\begin{enumerate}
    \item $\alpha_1=\alpha_3=0,$ we have the representative $\left\langle\nabla_4+ \frac 1 2\nabla_6\right\rangle_{\alpha =- 2};$
    \item $\alpha_1=0, \ \alpha_3\neq0,$ then choosing $x = \frac{\alpha_3}{\alpha_4}, \ y=1,$ we have the representative $\left\langle\nabla_3+\nabla_4+\frac 1 2\nabla_6\right\rangle_{\alpha =- 2};$
    \item $\alpha_1\neq0, \ \alpha_3=0,$ then choosing $x = -\frac{\alpha_6}{\alpha_1}, \ y=1,$ we have the representative $\left\langle\nabla_1-2\nabla_4-\nabla_6\right\rangle_{\alpha =- 2};$
    \item $\alpha_1\neq0, \ \alpha_3\neq0,$ then choosing $x = \frac{\alpha _1 \alpha _3^2}{\alpha _6^3}, \ y=\frac{\alpha _1 \alpha _3}{\alpha _6^2},$ we have the representative $\left\langle\nabla_1+\nabla_3+2\nabla_4+\nabla_6\right\rangle_{\alpha =- 2}.$
    \end{enumerate}
    \item $\alpha = - 2,$ $\alpha_4\neq2\alpha_6,$ then choosing $w=-\frac{x \alpha _1}{\alpha _4-2 \alpha _6},$ we can suppose $\alpha_1^*=0$ and consider following subcases: 
\begin{enumerate}
    \item $\alpha_4 = \alpha_6,$ $\alpha_3 = 0,$ then choosing $x=1,\ z=\frac{y \alpha _5}{3 \alpha _4}, \ t=-\frac{\alpha _2}{3\alpha _4},$ we have the representative $\left\langle\nabla_4+\nabla_6\right\rangle_{\alpha =- 2};$
    \item $\alpha_4 = \alpha_6,$ $\alpha_3 \neq 0,$ then choosing $x=1, \ y=\frac{\alpha _3}{\alpha _4}\ z=\frac{y \alpha _5}{3 \alpha _4}, \ t=-\frac{\alpha _2}{3\alpha _4},$ we have the representative $\left\langle\nabla_3+\nabla_4+\nabla_6\right\rangle_{\alpha =- 2};$
    \item $\alpha_4 \neq \alpha_6,$ then choosing $t=0,\ u=-\frac{y \alpha _2}{2 \left(\alpha _4-\alpha _6\right)},$ we have $\alpha_2^*=0.$ 
    \begin{enumerate}
        \item $\alpha_6\neq0,\ \alpha_3=0,$ then choosing $y=1,\ z=\frac{\alpha_5}{3\alpha_6},$ we have the family of representatives $\left\langle\nabla_4+\beta\nabla_6\right\rangle_{\alpha =- 2, \beta\neq 0,1};$
        \item $\alpha_6\neq0,\ \alpha_3\neq0,$ then choosing $x= \frac{\alpha _3}{\alpha _4}, \ y=1,\ z=\frac{\alpha _5}{3 \alpha _6},$ we have the family of representatives $\left\langle\nabla_3+\nabla_4+\beta\nabla_6\right\rangle_{\alpha =- 2, \beta\neq 0,1};$
        \item $\alpha_6=0,\ \alpha_5=\alpha_3=0,$ then we have the representative $\left\langle\nabla_4\right\rangle_{\alpha =- 2};$
        \item $\alpha_6=0,\ \alpha_5=0,\ \alpha_3\neq0, $ then choosing $x=\frac{\alpha _3}{\alpha _4}, \ y=1,$ we have the representative $\left\langle\nabla_3+\nabla_4\right\rangle_{\alpha =- 2};$
        \item $\alpha_6=0,\ \alpha_5\neq0,$ then choosing $x=\frac{\alpha _5}{\alpha _4}, \ y=1,$ we have the family of representatives $\left\langle\beta\nabla_3+\nabla_4+ \nabla_5\right\rangle_{\alpha =- 2}.$
    \end{enumerate}
    \end{enumerate}
    \item $\alpha \neq - 2,$ then choosing $z=-\frac{y \alpha _3}{(2+\alpha ) \alpha _4},$ we can suppose $\alpha_3^*=0$ and consider following subcases: 
    \begin{enumerate}
      
    \item $\alpha_4 + \alpha \alpha_6=0,$ then choosing $t=-\frac{x \alpha _2}{2\alpha _4}, \ u=0, \  w=0,$ we can suppose $\alpha_2^*=0$ and consider following subcases: 
\begin{enumerate}
    \item $\alpha_1 = \alpha_5=0,$ then we have the family of representatives $\left\langle \nabla_4-\frac{1}{\alpha}\nabla_6\right\rangle;$ 
    \item $\alpha_1 =0, \ \alpha_5\neq0,$ then choosing $x=1, \ y=\frac{\alpha _4}{\alpha _5},$ we have the family of representatives $\left\langle \nabla_4+\nabla_5-\frac{1}{\alpha}\nabla_6\right\rangle;$
    \item $\alpha_1 \neq 0, \ \alpha_5=0,$ then choosing $ x=-\frac{\alpha _6}{\alpha _1}, \ y=1,$ we have the family of representatives $\left\langle\nabla_1+\alpha\nabla_4-\nabla_6\right\rangle;$
    \item $\alpha_1 \neq 0, \ \alpha_5\neq0,$ then choosing $x=-\frac{\alpha _1 \alpha _5^2}{\alpha _6^3}, \ y=\frac{\alpha _1 \alpha _5}{\alpha _6^2},$ we have the family of representatives $\left\langle\nabla_1+\alpha\nabla_4+\nabla_5-\nabla_6\right\rangle.$
     \end{enumerate}
     
    \item  $\alpha_4 + \alpha \alpha_6\neq 0,$ then choosing $w=-\frac{x \alpha _1}{\alpha _4+\alpha  \alpha _6},$ we can suppose $\alpha_1^*=0$ and consider following subcases:
    
    \begin{enumerate}
    \item $\alpha_6 = \alpha_4,$ $\alpha = 1,$ $\alpha_2 = \alpha_5=0,$
    then we have the representative $\left\langle \nabla_4+\nabla_6\right\rangle_{\alpha = 1};$
    \item $\alpha_6 = \alpha_4,$ $\alpha = 1,$ $\alpha_2 =0,$ $\alpha_5\neq 0,$
    then choosing $x=1, \ y=\frac{\alpha _4}{\alpha _5},$ we have the representative $\left\langle \nabla_4+\nabla_5+\nabla_6\right\rangle_{\alpha = 1};$
\item $\alpha_6 = \alpha_4,$ $\alpha = 1,$ $\alpha_2 \neq 0,$ $\alpha_5=0,$
    then choosing $x=1, \ y=\frac{\alpha _2}{\alpha _4},$ we have the representative $\left\langle \nabla_2+\nabla_4+\nabla_6\right\rangle_{\alpha = 1};$
\item $\alpha_6 = \alpha_4,$ $\alpha = 1,$ $\alpha_2 \neq 0,$ $\alpha_5\neq 0,$
    then choosing $x=\frac{\alpha _2 \alpha _5}{\alpha _4^2}, \ y=\frac{\alpha _2}{\alpha _4},$ we have the representative $\left\langle \nabla_2+\nabla_4+\nabla_5+\nabla_6\right\rangle_{\alpha = 1};$
    \item $\alpha_6 = \alpha_4,$ $\alpha \neq 1,$ $\alpha_5 = 0,$
    then choosing $x=1, \ t=\frac{\alpha _2}{(\alpha-1 ) \alpha _4},$ we have the representative $\left\langle \nabla_4+\nabla_6\right\rangle_{\alpha \neq -1,1};$
    \item $\alpha_6 = \alpha_4,$ $\alpha \neq 1,$ $\alpha_5 \neq  0,$
    then choosing  $x=1, \ y=\frac{\alpha _4}{\alpha _5}, \ t=\frac{\alpha _2}{(\alpha-1) \alpha _4},$
    we have the representative $\left\langle \nabla_4+\nabla_5+\nabla_6\right\rangle_{\alpha \neq -1,1};$
    \item $\alpha_6 \neq \alpha_4,$ $\alpha_5 = 0,$
    then choosing $y=1, \ t=0, \ u = \frac{\alpha _2}{\alpha(\alpha _4-\alpha _6)},$ we have the family of representatives $\left\langle \nabla_4+\beta \nabla_6\right\rangle_{\beta \neq 1};$
    \item $\alpha_6 \neq \alpha_4,$ $\alpha_5 \neq 0,$
    then choosing $x=1,\ y=\frac{\alpha _4}{\alpha _5}, \ t=0, \ u = \frac{\alpha _2\alpha _4}{\alpha \alpha _5(\alpha _4-\alpha _6)},$ we have the family of representatives $\left\langle \nabla_4+\nabla_5+\beta \nabla_6\right\rangle_{\beta \neq 1};$
    \end{enumerate}
    \end{enumerate}
    \end{enumerate}
\end{enumerate}

 Summarizing all cases of the central extension of the algebra ${\mathfrak N}_{14}^{\alpha},$ we have the following distinct orbits:

in case of $\alpha=-\frac12$: \begin{center} 
$\langle\beta\nabla_3+\nabla_5+\nabla_6\rangle,$
\end{center}

in case of $\alpha=1$: \begin{center} 
$\langle\nabla_2+\nabla_4+\nabla_6\rangle,$  $\langle\nabla_2+\nabla_4+\nabla_5+\nabla_6\rangle,$
\end{center}

in case of $\alpha=0$: \begin{center} 
$\langle\nabla_3+\beta\nabla_4+\gamma\nabla_6+\nabla_7\rangle^{O(\beta, \gamma) \simeq O(\beta, - \gamma)}$,  
$\langle\nabla_2+\nabla_5+\beta\nabla_6+\nabla_7\rangle^{O(\beta) \simeq O(-\beta)},$ 
$\langle\nabla_5+\beta\nabla_6+\nabla_7\rangle^{O(\beta) \simeq O(-\beta)},$ 
$\langle\nabla_1+\nabla_3+\nabla_6\rangle,$
\end{center}

in case of $\alpha=-2$: \begin{center} 
$\langle\nabla_1+\nabla_3+2\nabla_4+\nabla_6\rangle$, 
$\langle\nabla_3+\nabla_4+\beta\nabla_6\rangle,$ 
$\langle\beta\nabla_3+\nabla_4+\nabla_5\rangle,$ 
\end{center}

for any $\alpha\in \mathbb{C}$:
\begin{center}  
$\langle \nabla_3+\nabla_6\rangle,$  
$\langle\nabla_4+\beta\nabla_6\rangle,$ $\langle\nabla_4+\nabla_5+\beta\nabla_6\rangle_{\alpha\neq-2},$ 
$\langle\nabla_6\rangle_{\alpha\neq0},$  $\langle\nabla_4-\frac{1}{\alpha}\nabla_6\rangle_{\alpha\neq0},$  $\langle\nabla_4+\nabla_5-\frac{1}{\alpha}\nabla_6\rangle_{\alpha\neq0},$
$\langle\nabla_1+\alpha\nabla_4-\nabla_6\rangle_{\alpha\neq0},$ 
$\langle\nabla_1+\alpha\nabla_4+\nabla_5-\nabla_6\rangle_{\alpha\neq0,-2},$  
\end{center}
which gives the following new algebras (see section \ref{secteoA}):

\begin{center}
${\rm B}_{88}^{\beta},$ 
${\rm B}_{89},$
${\rm B}_{90},$
${\rm B}_{91}^{\beta,\gamma},$
${\rm B}_{92}^{\beta},$ 
${\rm B}_{93}^{\beta},$ 
${\rm B}_{94},$ 
${\rm B}_{95},$ 
${\rm B}_{96}^{\beta},$
${\rm B}_{97}^{\beta},$ 
${\rm B}_{98}^{\alpha},$ 
${\rm B}_{99}^{\alpha,\beta},$ 

${\rm B}_{100}^{\alpha\neq-2, \beta},$ 
${\rm B}_{101}^{\alpha\neq0},$ 
${\rm B}_{102}^{\alpha\neq0},$ 
${\rm B}_{103}^{\alpha\neq0},$ 
${\rm B}_{104}^{\alpha\neq0},$ 
${\rm B}_{105}^{\alpha\neq0,-2}.$
\end{center}

\subsection{$1$-dimensional central extensions of $4$-dimensional  $3$-step nilpotent bicommutative algebras}

\subsubsection{The description of second cohomology space.}

In the following table, we give the description
of the second cohomology space of $4$-dimensional $3$-step nilpotent bicommutative algebras.

\begin{longtable}{ll llllll}
\hline 

${\mathcal B}^4_{01}$ &$:$     & $e_1 e_1 = e_2$  & $e_2 e_1=e_3$ &&&& \\
\multicolumn{8}{l}{
${\rm H}^2({\mathcal B}_{01}^4)=
\Big\langle  [\Delta_{12}], [\Delta_{14}], [\Delta_{31}], [\Delta_{41}], [\Delta_{44}]
\Big\rangle $}\\
\hline

${\mathcal B}^4_{02}(\alpha)$ &$:$  & $e_1 e_1 = e_2$ & $e_1 e_2=e_3$ & $e_2 e_1=\alpha e_3$ & &&\\

\multicolumn{8}{l}{
${\rm H}^2({\mathcal B}_{02}^4(\alpha))=
\Big\langle  [\Delta_{14}], [\Delta_{21}], [\Delta_{13}+\alpha\Delta_{22}+\alpha\Delta_{31}], [\Delta_{41}], [\Delta_{44}]
\Big\rangle $}\\

\hline

${\mathcal B}^4_{04}(\alpha)$  &$:$    &
$e_1 e_1 = e_2$ &   $e_1e_2=e_4$ & $e_2e_1=\alpha e_4$ & $e_3e_3=e_4$&&\\

\multicolumn{8}{l}{
${\rm H}^2({\mathcal B}_{04}^4(\alpha))=
\Big\langle  [\Delta_{13}], [\Delta_{21}], [\Delta_{31}], [\Delta_{33}]
\Big\rangle $}\\

\hline

${\mathcal B}^4_{05}$ &$:$  &
$e_1 e_1 = e_2$ &   $e_1e_2=e_4$ & $e_1e_3= e_4$ & $e_2e_1= e_4$ & $e_3e_3=e_4$ & \\

\multicolumn{8}{l}{
${\rm H}^2({\mathcal B}_{05}^4)=
\Big\langle  [\Delta_{13}], [\Delta_{21}], [\Delta_{31}], [\Delta_{33}]
\Big\rangle $}\\

\hline

${\mathcal B}^4_{06}(\alpha\neq0)$   &$:$    &
$e_1 e_1 = e_2$ &   $e_1e_2=e_4$ & $e_1e_3=e_4$ & $e_2e_1=\alpha e_4$&& \\

\multicolumn{8}{l}{
${\rm H}^2({\mathcal B}_{06}^4(\alpha))=
\Big\langle  [\Delta_{13}], [\Delta_{21}], [\Delta_{31}],  [\Delta_{14}+\alpha\Delta_{22}+\alpha\Delta_{23}+ \alpha\Delta_{41}], [\Delta_{33}]
\Big\rangle $}\\

\hline

${\mathcal B}^4_{07}$  &$:$   &
$e_1 e_1 = e_2$ & $e_2e_1= e_4$ & $e_3e_3=e_4$ &&& \\

\multicolumn{8}{l}{
${\rm H}^2({\mathcal B}_{07}^4)=
\Big\langle  [\Delta_{12}], [\Delta_{13}], [\Delta_{31}], [\Delta_{33}]
\Big\rangle $}\\

\hline

${\mathcal B}^4_{08}$  &$:$  &
$e_1 e_1 = e_2$ & $e_1e_3=e_4$ & $e_2e_1= e_4$ &&& \\

\multicolumn{8}{l}{
${\rm H}^2({\mathcal B}_{08}^4)=
\Big\langle  [\Delta_{12}], [\Delta_{21}], [\Delta_{31}], [\Delta_{33}], [\Delta_{23}+\Delta_{41}]
\Big\rangle $}\\

\hline

${\mathcal B}^4_{09}$   &$:$  &
$e_1 e_1 = e_2$ & $e_1e_2=e_4$ &  $e_3e_1= e_4$ &&& \\

\multicolumn{8}{l}{
${\rm H}^2({\mathcal B}_{09}^4)=
\Big\langle  [\Delta_{13}], [\Delta_{21}], [\Delta_{31}],  [\Delta_{33}], [\Delta_{14}+\Delta_{32}]
\Big\rangle $}\\

\hline

${\mathcal B}^{4}_{10}$   &$:$   & $e_1 e_2=e_3$  & $e_1 e_3=e_4$ & $e_2e_1=e_4$ & $e_3e_2=e_4$&& \\

\multicolumn{8}{l}{
${\rm H}^2({\mathcal B}_{10}^4)=
\Big\langle  [\Delta_{11}], [\Delta_{13}], [\Delta_{22}],  [\Delta_{32}]
\Big\rangle $}\\

\hline

${\mathcal B}^{4}_{11}$   &$:$  &  $e_1e_2=e_3$ & $e_1 e_3=e_4$  & $e_3 e_2=e_4$ &&& \\

\multicolumn{8}{l}{
${\rm H}^2({\mathcal B}_{11}^4)=
\Big\langle  [\Delta_{11}], [\Delta_{21}], [\Delta_{22}],  [\Delta_{32}], [\Delta_{14}+\Delta_{33}+\Delta_{42}]
\Big\rangle $}\\

\hline

${\mathcal B}^{4}_{12}$  &$:$    &   $e_1 e_1=e_4$ &  $e_1e_2=e_3$ & $e_2 e_1=e_4$ & $e_3 e_2=e_4$&&\\

\multicolumn{8}{l}{
${\rm H}^2({\mathcal B}_{12}^4)=
\Big\langle  [\Delta_{11}], [\Delta_{13}], [\Delta_{22}],  [\Delta_{32}]
\Big\rangle $}\\

\hline

${\mathcal B}^{4}_{13}$   &$:$   &  $e_1e_2=e_3$ & $e_2 e_1=e_4$  & $e_3 e_2=e_4$ &&& \\

\multicolumn{8}{l}{
${\rm H}^2({\mathcal B}_{13}^4)=
\Big\langle  [\Delta_{11}], [\Delta_{13}], [\Delta_{22}],  [\Delta_{32}]
\Big\rangle $}\\

\hline

${\mathcal B}^{4}_{14}$  &$:$  &  $e_1e_2=e_3$ & $e_1 e_3=e_4$ & $e_2 e_1=e_4$  & $e_2 e_2=e_4$&& \\

\multicolumn{8}{l}{
${\rm H}^2({\mathcal B}_{14}^4)=
\Big\langle  [\Delta_{11}], [\Delta_{13}], [\Delta_{22}],  [\Delta_{32}]
\Big\rangle $}\\

\hline

${\mathcal B}^{4}_{15}$  &$:$   &  $e_1e_2=e_3$ & $e_1 e_3=e_4$ & $e_2 e_1=e_4$ &&& \\

\multicolumn{8}{l}{
${\rm H}^2({\mathcal B}_{15}^4)=
\Big\langle  [\Delta_{11}], [\Delta_{13}], [\Delta_{22}],  [\Delta_{32}]
\Big\rangle $}\\

\hline

${\mathcal B}^{4}_{16}$  &$:$   &  $e_1e_2=e_3$ &  $e_1 e_3=e_4$ & $e_2 e_2=e_4$ &&& \\

\multicolumn{8}{l}{
${\rm H}^2({\mathcal B}_{16}^4)=
\Big\langle  [\Delta_{11}], [\Delta_{21}], [\Delta_{22}],  [\Delta_{32}], [\Delta_{14}+\Delta_{23}]
\Big\rangle $}\\
\hline

${\mathcal B}^{4}_{17}$   &$:$   &  $e_1e_2=e_3$ &  $e_1 e_3=e_4$ & &&& \\

\multicolumn{8}{l}{
${\rm H}^2({\mathcal B}_{17}^4)=
\Big\langle  [\Delta_{11}], [\Delta_{14}], [\Delta_{21}], [\Delta_{22}],  [\Delta_{32}]
\Big\rangle $}\\

\hline

${\mathcal B}^{4}_{18}$   &$:$  &  $e_1 e_1=e_4$  & $e_1e_2=e_3$ &   $e_3 e_2=e_4$ &&& \\

\multicolumn{8}{l}{
${\rm H}^2({\mathcal B}_{18}^4)=
\Big\langle  [\Delta_{13}], [\Delta_{21}], [\Delta_{22}],  [\Delta_{32}], [\Delta_{31}+\Delta_{42}]
\Big\rangle $}\\

\hline

${\mathcal B}^{4}_{19}$   &$:$  &  $e_1e_2=e_3$  & $e_3 e_2=e_4$ & &&& \\

\multicolumn{8}{l}{
${\rm H}^2({\mathcal B}_{19}^4)=
\Big\langle  [\Delta_{11}], [\Delta_{13}], [\Delta_{21}],  [\Delta_{22}], [\Delta_{42}]
\Big\rangle $}\\
\hline

\end{longtable}

\subsubsection{Central extensions of ${\mathcal B}_{01}^4$}
	Let us use the following notations:
	\begin{longtable}{lllllll} 
	$\nabla_1 = [\Delta_{12}],$ & $\nabla_2 = [\Delta_{14}],$ &
	$\nabla_3 = [\Delta_{31}],$ &
	$\nabla_4 = [\Delta_{41}],$&
	$\nabla_5 = [\Delta_{44}].$ 
	\end{longtable}	
	
Take $\theta=\sum\limits_{i=1}^{5}\alpha_i\nabla_i\in {\rm H^2}({\mathcal B}_{01}^4).$
	The automorphism group of ${\mathcal B}_{01}^4$ consists of invertible matrices of the form
	$$\phi=
	\begin{pmatrix}
	x &  0  & 0 & 0\\
	y &  x^2  & 0 & 0\\
	z &  xy  & x^3 & t\\
    u &  0  & 0 & r
	\end{pmatrix}.
	$$
	Since
	$$
	\phi^T\begin{pmatrix}
	0 & \alpha_1 & 0 & \alpha_2\\
	0 &	0 & 0 & 0\\
	\alpha_3 &  0  & 0 & 0\\
	\alpha_4 & 0  & 0 & \alpha_5
	\end{pmatrix} \phi=\begin{pmatrix}
	\alpha^* & \alpha_1^* & 0 & \alpha_2^*\\
	\alpha^{**} &	0 & 0 & 0\\
	\alpha_3^* &  0  & 0 & 0\\
	\alpha_4^* & 0  & 0 & \alpha_5^*
	\end{pmatrix},
	$$
	 we have that the action of ${\rm Aut} ({\mathcal B}_{01}^4)$ on the subspace
$\langle \sum\limits_{i=1}^{5}\alpha_i\nabla_i  \rangle$
is given by
$\langle \sum\limits_{i=1}^{5}\alpha_i^{*}\nabla_i\rangle,$
where
\begin{longtable}{lcl}
$\alpha^*_1$&$=$&$x^3 \alpha _1,$ \\
$\alpha_2^*$&$=$&$r \left(x \alpha _2+u \alpha _5\right),$\\
$\alpha^*_3$&$=$&$x^4 \alpha _3,$ \\
$\alpha^*_4$&$=$&$t x \alpha _3+r x \alpha _4+r u \alpha _5,$\\
$\alpha_5^*$&$=$&$r^2 \alpha _5.$\\
\end{longtable}

We are interested only in the cases with 
 \begin{center}
$(\alpha_2,\alpha_4,\alpha_5)\neq (0,0,0),$  $\alpha_3\neq 0.$ 
 \end{center} 
Since $\alpha_3\neq0,$ then choosing $t=-\frac{r \left(x \alpha _4+u \alpha _5\right)}{x \alpha _3},$ we have $\alpha_4^*=0.$
\begin{enumerate}
    \item If $\alpha_5\neq0,$ then choosing $u=-\frac{x \alpha _2}{\alpha _5},$ we have $\alpha_2^*=0.$
    \begin{enumerate}
        \item $\alpha_1=0,$ then choosing $x=1,\ r=\sqrt{\frac{\alpha_3}{\alpha_5}},$ we have the representative $\left\langle\nabla_3+\nabla_5\right\rangle;$
        \item $\alpha_1\neq0,$ then choosing $x=\frac{\alpha_1}{\alpha_3},\ r=\frac{\alpha_1^2}{\alpha_3\sqrt{\alpha_3\alpha_5}},$ we have the representative $\left\langle\nabla_1+\nabla_3+ \nabla_5\right\rangle.$
    \end{enumerate}
    \item If $\alpha_5=0,$ then $\alpha_2\neq0.$
     \begin{enumerate}
        \item $\alpha_1=0,$ then choosing $x=1,\ r=\frac{\alpha_3}{\alpha_2},$ we have the representative $\left\langle\nabla_2+\nabla_3\right\rangle;$
        \item $\alpha_1\neq0,$ then choosing $x=\frac{\alpha_1}{\alpha_3},\ r=\frac{\alpha_1^3}{\alpha_2\alpha_3^2},$ we have the representative $\left\langle\nabla_1+\nabla_2+ \nabla_3\right\rangle.$
    \end{enumerate}
\end{enumerate}

Therefore, we have the following distinct orbits 
\begin{longtable} {lllll}
$\langle\nabla_3+\nabla_5\rangle,$ & $\langle\nabla_1+\nabla_3+\nabla_5\rangle,$ & $\langle\nabla_2+\nabla_3\rangle,$ &  $\langle\nabla_1+\nabla_2+\nabla_3\rangle,$ \\
\end{longtable}
which gives the following new algebras (see section \ref{secteoA}):

\begin{center}
${\rm B}_{106},$
${\rm B}_{107},$
${\rm B}_{108},$
${\rm B}_{109}.$
\end{center}

\subsubsection{Central extensions of ${\mathcal B}_{02}^{\alpha}$}
	Let us use the following notations:
	\begin{longtable}{lllllll} 
	$\nabla_1 = [\Delta_{14}],$ & $\nabla_2 = [\Delta_{21}],$ &
	$\nabla_3 = [\Delta_{13}+\alpha\Delta_{22}+\alpha\Delta_{31}],$ &
	$\nabla_4 = [\Delta_{41}],$&
	$\nabla_5 = [\Delta_{44}].$ 
	\end{longtable}	
	
Take $\theta=\sum\limits_{i=1}^{5}\alpha_i\nabla_i\in {\rm H^2}({\mathcal B}_{02}^{\alpha}).$
	The automorphism group of ${\mathcal B}_{02}^{\alpha}$ consists of invertible matrices of the form
	$$\phi=
	\begin{pmatrix}
	x &  0  & 0 & 0\\
	y &  x^2  & 0 & 0\\
	z &  (1+\alpha)xy  & x^3 & t\\
    u &  0  & 0 & r
	\end{pmatrix}.
	$$
	Since
	$$
	\phi^T\begin{pmatrix}
	0 & 0 & \alpha_3 & \alpha_1\\
	\alpha_2 &  \alpha\alpha_3  & 0 & 0\\
	\alpha\alpha_3 &	0 & 0 & 0\\
	\alpha_4 & 0  & 0 & \alpha_5
	\end{pmatrix} \phi=\begin{pmatrix}
	\alpha^* & \alpha^{**} & \alpha_3^* & \alpha_1^*\\
	\alpha_2^{*}+\alpha \alpha^{**} &	\alpha\alpha_3^* & 0 & 0\\
	\alpha\alpha_3^* &  0  & 0 & 0\\
	\alpha_4^* & 0  & 0 & \alpha_5^*
	\end{pmatrix},
	$$
	 we have that the action of ${\rm Aut} ({\mathcal B}_{02}^{\alpha})$ on the subspace
$\langle \sum\limits_{i=1}^{5}\alpha_i\nabla_i  \rangle$
is given by
$\langle \sum\limits_{i=1}^{5}\alpha_i^{*}\nabla_i\rangle,$
where
\begin{longtable}{lcl}
$\alpha^*_1$&$=$&$r x \alpha _1+t x \alpha _3+r u \alpha _5,$ \\
$\alpha_2^*$&$=$&$x^2(x\alpha_2-y\alpha(\alpha-1) \alpha _3),$\\
$\alpha^*_3$&$=$&$x^4 \alpha _3,$ \\
$\alpha^*_4$&$=$&$tx\alpha\alpha_3+rx\alpha_4+ru\alpha _5,$\\
$\alpha_5^*$&$=$&$r^2 \alpha _5.$\\
\end{longtable}

We are interested only in the cases with 
 \begin{center}
$\alpha_3\neq 0, \ $   $(\alpha_1,\alpha_4,\alpha_5)\neq (0,0,0).$
 \end{center}
$\alpha_3\neq0,$ then choosing $t=-\frac{r \left(x \alpha _1+u \alpha _5\right)}{x \alpha _3},$ we have $\alpha_1^*=0.$ Now we consider following cases:
\begin{enumerate}
\item $\alpha_5=0,$ then $\alpha_4\neq0.$
\begin{enumerate}
    \item $\alpha\in\{0,1\},\ \alpha_2=0,$ then choosing $x=1,\ r=\frac{\alpha_3}{\alpha_4},$ we have the representative  $\left\langle\nabla_3+\nabla_4\right\rangle;$
    \item $\alpha\in\{0,1\},\ \alpha_2\neq0,$ then choosing $x=\frac{\alpha_2}{\alpha_3},\ r=\frac{\alpha_2^3}{\alpha_3^2\alpha_4},$ we have the representative  $\left\langle\nabla_2+\nabla_3+ \nabla_4\right\rangle;$
    \item $\alpha\not \in \{0,1\},$ then choosing $x=1,\ y=\frac{\alpha _2}{(\alpha-1)\alpha\alpha_3},\  r=\frac{\alpha_3}{\alpha_4},$ we have the representative  $\left\langle\nabla_3+\nabla_4\right\rangle.$
\end{enumerate}
\item $\alpha_5\neq0.$
\begin{enumerate}
    \item $\alpha=1,\ \alpha_4=\alpha_2=0,$ then choosing $x=1,\  r=\sqrt{\frac{\alpha_3}{\alpha_5}},$ we have the representative  $\left\langle\nabla_3+\nabla_5\right\rangle;$
    \item $\alpha=1,\ \alpha_4=0,\ \alpha_2\neq0,$ then choosing $x=\frac{\alpha_2}{\alpha_3},\  r=\frac{\alpha_2^2}{\alpha_3\sqrt{\alpha_3\alpha_5}},$ we have the representative $\left\langle\nabla_2+\nabla_3+ \nabla_5\right\rangle;$
    \item $\alpha=1,\ \alpha_4\neq0,$ then choosing $x=\frac{\alpha_4}{\sqrt{\alpha_3\alpha_5}},\  r=\frac{\alpha_4^2}{\alpha_5\sqrt{\alpha_3\alpha_5}},$ we have the family of representatives  $\left\langle\beta\nabla_2+\nabla_3+\nabla_4+ \nabla_5\right\rangle^{ O(\beta)\simeq O(-\beta)};$
    \item $\alpha=0,\ \alpha_2=0,$ then choosing $x=1,\  r=\sqrt{\frac{\alpha_3}{\alpha_5}},\ u=-\frac{\alpha_4}{\alpha_5},$ we have the representative  $\left\langle\nabla_3+\nabla_5\right\rangle;$
    \item $\alpha=0,\ \alpha_2\neq0,$ then choosing $x=\frac{\alpha_2}{\alpha_3},\  r=\frac{\alpha_2^2}{\alpha_3\sqrt{\alpha_3\alpha_5}}\ u=-\frac{\alpha_2\alpha_4}{\alpha_3\alpha_5},$ we have the representative  $\left\langle\nabla_2+\nabla_3+ \nabla_5\right\rangle;$
    \item $\alpha\not \in \{0,1\},  $ then choosing $x=1,\ y=\frac{\alpha_2}{(-1+\alpha)\alpha\alpha_3}, \  r=\sqrt{\frac{\alpha_3}{\alpha_5}},\ u=\frac{\alpha_4}{(-1+\alpha)\alpha_5},$ we have the representative  $\left\langle\nabla_3+ \nabla_5\right\rangle.$
\end{enumerate}
\end{enumerate}

Summarizing all cases, we have the following distinct orbits

in case of $\alpha=0:$
\begin{longtable} {llll}
$\langle\nabla_2+\nabla_3+\nabla_4\rangle,$ & $\langle\nabla_2+\nabla_3+\nabla_5\rangle,$ \\
\end{longtable}
in case of $\alpha=1$
\begin{longtable} {llll}
$\langle\beta\nabla_2+\nabla_3+\nabla_4+\nabla_5\rangle^{O(\beta)\simeq O(-\beta)}$, &
$\langle\nabla_2+\nabla_3+\nabla_4\rangle,$ & $\langle\nabla_2+\nabla_3+\nabla_5\rangle,$ \\
\end{longtable}

in case of $\alpha \in \mathbb{C}$
\begin{longtable} {llll}
$\langle\nabla_3+\nabla_4\rangle,$ & $\langle\nabla_3+\nabla_5\rangle,$\\
\end{longtable}
which gives the following new algebras (see section \ref{secteoA}, 
as we are interested in non-commutative algebras, ew do not consider ${\rm B}_{115}^{1}$):

\begin{center}

${\rm B}_{110},$
${\rm B}_{111},$
${\rm B}_{112}^{\beta},$ 
${\rm B}_{113},$
${\rm B}_{114},$ 
${\rm B}_{115}^{\alpha\neq1},$
${\rm B}_{116}^{\alpha}.$ 

\end{center}

\subsubsection{Central extensions of ${\mathcal B}_{06}^4(\alpha\neq0)$}
Let us use the following notations:
	\begin{longtable}{lllllll} 
$\nabla_1 = [\Delta_{13}],$&
$\nabla_2 = [\Delta_{21}],$ & $\nabla_3=[\Delta_{31}],$ &
$\nabla_4=[\Delta_{14}+\alpha\Delta_{22}+\alpha\Delta_{23}+ \alpha\Delta_{41}],$ & $\nabla_5 = [\Delta_{33}].$ 
	\end{longtable}	
	
Take $\theta=\sum\limits_{i=1}^{5}\alpha_i\nabla_i\in {\rm H^2}({\mathcal B}_{06}^4(\alpha\neq0)).$
The automorphism group of ${\mathcal B}_{06}^4(\alpha\neq0)$ consists of invertible matrices of the form
$$\phi=
	\begin{pmatrix}
	x &  0  & 0 & 0\\
	y &  x^2  & 0 & 0\\
	z &  0  & x^2 & 0\\
    u &  x((1+\alpha)y+z)  & v & x^3
	\end{pmatrix}.
	$$
	Since
	$$
	\phi^T\begin{pmatrix}
	0 & 0 & \alpha_1 & \alpha_4\\
	\alpha_2 &  \alpha\alpha_4  & \alpha\alpha_4 & 0\\
	\alpha_3 &	0 & \alpha_5 & 0\\
	\alpha\alpha_4 & 0  & 0 & 0
	\end{pmatrix} \phi=\begin{pmatrix}
	\alpha^* & \alpha^{**} & \alpha_1^*+\alpha^{**} & \alpha_4^*\\
	\alpha_2^*+\alpha \alpha^{**} &  \alpha\alpha_4^* & \alpha\alpha_4^* & 0\\
	\alpha_3^* &	0& \alpha_5^* & 0\\
	\alpha\alpha_4^* & 0  & 0 & 0
	\end{pmatrix},
$$
we have that the action of ${\rm Aut} ({\mathcal B}_{06}^4(\alpha\neq0))$ on the subspace
$\langle \sum\limits_{i=1}^{5}\alpha_i\nabla_i  \rangle$
is given by
$\langle \sum\limits_{i=1}^{5}\alpha_i^{*}\nabla_i\rangle,$
where
\begin{longtable}{lcl}
$\alpha^*_1$&$=$&$x(x^2\alpha_1+(v-x((1+\alpha)y+z))\alpha _4+x z \alpha _5),$ \\
$\alpha_2^*$&$=$&$x^2(x\alpha _2+\alpha ((1-\alpha)y+ z)\alpha_4),$\\
$\alpha^*_3$&$=$&$x \left(x^2 \alpha _3+v \alpha  \alpha _4+x z \alpha _5\right),$ \\
$\alpha^*_4$&$=$&$x^4 \alpha _4,$\\
$\alpha_5^*$&$=$&$x^4 \alpha_5.$\\
\end{longtable}
We are interested only in the cases with $\alpha_4\neq 0.$

Choosing $v=\frac{x(y \alpha\alpha_4(2\alpha\alpha_4+\alpha_5-\alpha  \alpha_5)-x\alpha\alpha_1\alpha_4-x \alpha _2(\alpha _4-\alpha_5))}{\alpha\alpha_4^2},$ 
$z=y (\alpha-1)-\frac{x \alpha _2}{\alpha  \alpha _4},$ we have $\alpha_1^*=\alpha_2^*=0.$
\begin{enumerate}
    \item $2\alpha^2\alpha_4=(\alpha-1)^2\alpha_5,\ \alpha_3=0,$ then $\alpha\neq1$ and we have the family of representatives  $\left\langle\nabla_4+\frac{2\alpha^2}{(\alpha-1)^2}\nabla_5\right\rangle;$
    \item $2\alpha^2\alpha_4=(\alpha-1)^2\alpha_5,\ \alpha_3\neq0,$ then $\alpha\neq1$ and choosing $x=\frac{\alpha_3}{\alpha_4},$ we have the family of representatives $\left\langle\nabla_3+\nabla_4+\frac{2\alpha^2}{(\alpha-1)^2}\nabla_5\right\rangle;$
    \item $2\alpha^2\alpha_4\neq(\alpha-1)^2\alpha_5,$ then choosing $x=1,\ y=-\frac{x\alpha_3}{2\alpha^2\alpha _4-(\alpha-1)^2\alpha _5},$ we have the family of representatives  $\left\langle\nabla_4+\beta\nabla_5\right\rangle_{\beta \neq \frac{2\alpha^2}{(\alpha-1)^2}}.$
\end{enumerate}

Summarizing all cases, we have the following distinct orbits
\begin{center}  
$\langle\nabla_4+\beta\nabla_5\rangle,$  $\langle\nabla_3+\nabla_4+\frac{2\alpha^2}{(\alpha-1)^2}\nabla_5\rangle_{\alpha\neq1},$ 
\end{center}
which gives the following new algebras (see section \ref{secteoA}):

\begin{center}
${\rm B}_{117}^{\alpha\neq0, \beta},$
${\rm B}_{118}^{\alpha\neq0,1}.$

\end{center}

\subsubsection{Central extensions of ${\mathcal B}_{08}^4$}
Let us use the following notations:
	\begin{longtable}{lllllll} 
$\nabla_1 = [\Delta_{12}],$& $\nabla_2 = [\Delta_{21}],$ & $\nabla_3=[\Delta_{31}],$ &
$\nabla_4=[\Delta_{33}],$ &
$\nabla_5=[\Delta_{23}+\Delta_{41}].$
\end{longtable}	
	
Take $\theta=\sum\limits_{i=1}^{5}\alpha_i\nabla_i\in {\rm H^2}({\mathcal B}_{08}^4).$
The automorphism group of ${\mathcal B}_{08}^4$ consists of invertible matrices of the form
$$\phi=
	\begin{pmatrix}
	x &  0  & 0 & 0\\
	y &  x^2  & 0 & 0\\
	z &  0  & x^2 & 0\\
    u &  x(y+z)  & v & x^3
	\end{pmatrix}.
	$$
	Since
	$$
	\phi^T\begin{pmatrix}
	0 & \alpha_1 & 0  & 0\\
	\alpha_2 & 0  & \alpha_5 & 0\\
	\alpha_3 &	0 & \alpha_4 & 0\\
	\alpha_5 & 0  & 0 & 0
	\end{pmatrix} \phi=\begin{pmatrix}
	\alpha^* & \alpha_1^* & \alpha^{**}  & 0\\
	\alpha_2^*+\alpha^{**} & 0  & \alpha_5^* & 0\\
	\alpha_3^* &	0 & \alpha_4^* & 0\\
	\alpha_5^* & 0  & 0 & 0
	\end{pmatrix},
$$
we have that the action of ${\rm Aut} ({\mathcal B}_{08})$ on the subspace
$\langle \sum\limits_{i=1}^{5}\alpha_i\nabla_i  \rangle$
is given by
$\langle \sum\limits_{i=1}^{5}\alpha_i^{*}\nabla_i\rangle,$
where
\begin{longtable}{lcllcllcl}
$\alpha^*_1$&$=$&$x^3 \alpha _1,$ &
$\alpha_2^*$&$=$&$x^2(x\alpha _2-z \alpha _4+2 z \alpha _5),$&
$\alpha^*_3$&$=$&$x(x^2 \alpha_3+xz\alpha _4+v \alpha _5),$ \\
$\alpha^*_4$&$=$&$x^4 \alpha _4,$ &
$\alpha_5^*$&$=$&$x^4 \alpha_5.$\\
\end{longtable}
We are interested only in the cases with $\alpha_5\neq 0.$ Choosing $v=-\frac{x \left(x \alpha _3+z \alpha _4\right)}{\alpha _5},$  we have $\alpha_3^*=0.$
\begin{enumerate}
    \item $\alpha_4=2\alpha_5,\ \alpha_2=\alpha_1=0,$ then we have the representative $\left\langle 2\nabla_4+ \nabla_5\right\rangle;$
    \item $\alpha_4=2\alpha_5,\ \alpha_2=0,\ \alpha_1\neq0,$ then choosing $x=\frac{\alpha_1}{\alpha_5},$ we have the representative  $\left\langle\nabla_1+2\nabla_4+ \nabla_5\right\rangle;$
    \item $\alpha_4=2\alpha_5,\ \alpha_2\neq0,$ then choosing $x=\frac{\alpha_2}{\alpha_5},$ we have the family of representatives  $\left\langle\alpha\nabla_1+\nabla_2+2\nabla_4+ \nabla_5\right\rangle;$
    \item $\alpha_4\neq2\alpha_5,\ \alpha_1=0,$ then choosing $x=1,\ z=\frac{\alpha _2}{\alpha _4-2 \alpha _5},$ we have the family of representatives  $\left\langle\alpha\nabla_4+ \nabla_5\right\rangle_{\alpha\neq2};$
    \item $\alpha_4\neq2\alpha_5,\ \alpha_1\neq0,$ then choosing $x=\frac{\alpha_1}{\alpha_5},\ z=\frac{\alpha_1\alpha_2}{(\alpha_4-2\alpha_5)\alpha_5},$ we have the family of representatives $\left\langle\nabla_1+\alpha\nabla_4+ \nabla_5\right\rangle_{\alpha\neq2}.$
\end{enumerate}

Summarizing all cases, we have the following distinct orbits
\begin{center}  
$\langle\nabla_1+\alpha\nabla_4+\nabla_5\rangle,$   $\langle\alpha\nabla_4+\nabla_5\rangle,$  $\langle\alpha\nabla_1+\nabla_2+2\nabla_4+\nabla_5\rangle,$ 
\end{center}
which gives the following new algebras (see section \ref{secteoA}):

\begin{center}
${\rm B}_{119}^{\alpha},$
${\rm B}_{120}^{\alpha},$
${\rm B}_{121}^{\alpha}.$

\end{center}

\subsubsection{Central extensions of ${\mathcal B}_{09}^4$}
Let us use the following notations:
	\begin{longtable}{lllllll} 
$\nabla_1 = [\Delta_{13}],$& $\nabla_2 = [\Delta_{21}],$ & $\nabla_3=[\Delta_{31}],$ &
$\nabla_4=[\Delta_{33}],$ &
$\nabla_5=[\Delta_{14}+\Delta_{32}].$
\end{longtable}	
	
Take $\theta=\sum\limits_{i=1}^{5}\alpha_i\nabla_i\in {\rm H^2}({\mathcal N}^4_{09}).$
The automorphism group of ${\mathcal N}^4_{09}$ consists of invertible matrices of the form
$$\phi=
	\begin{pmatrix}
	x &  0  & 0 & 0\\
	y &  x^2  & 0 & 0\\
	z &  0  & x^2 & 0\\
    u &  x(y+z)  & v & x^3
	\end{pmatrix}.
	$$
	Since
	$$
	\phi^T\begin{pmatrix}
	0 & 0 & \alpha_1  & \alpha_5\\
	\alpha_2 & 0  & 0 & 0\\
	\alpha_3 &	\alpha_5 & \alpha_4 & 0\\
	0 & 0  & 0 & 0
	\end{pmatrix} \phi=\begin{pmatrix}
	\alpha^* &  \alpha^{**}  & \alpha_1^* & \alpha_5^*\\
	\alpha_2^* & 0  & 0 & 0\\
	\alpha_3^*+\alpha^{**} & \alpha_5^* & \alpha_4^* & 0\\
	0 & 0 & 0 & 0
	\end{pmatrix},
$$
we have that the action of ${\rm Aut} ({\mathcal N}^4_{09})$ on the subspace
$\langle \sum\limits_{i=1}^{5}\alpha_i\nabla_i  \rangle$
is given by
$\langle \sum\limits_{i=1}^{5}\alpha_i^{*}\nabla_i\rangle,$
where
\begin{longtable}{lcllcllcl}
$\alpha^*_1$&$=$&$x(x^2\alpha_1+xz\alpha_4+v\alpha_5),$ &
$\alpha_2^*$&$=$&$x^3 \alpha _2,$\\
$\alpha^*_3$&$=$&$x^2(x\alpha_3+z(\alpha _4-2 \alpha _5)),$ &
$\alpha^*_4$&$=$&$x^4 \alpha _4,$&
$\alpha_5^*$&$=$&$x^4 \alpha_5.$\\
\end{longtable}
We are interested only in the cases with $\alpha_5\neq 0.$ Choosing $v=-\frac{x \left(x \alpha _1+z \alpha _4\right)}{\alpha _5},$ we have $\alpha_1^*=0.$
\begin{enumerate}
    \item $\alpha_4=2\alpha_5,\ \alpha_3=\alpha_2=0,$ then we have the representative  $\left\langle2\nabla_4+\nabla_5\right\rangle;$
    \item $\alpha_4=2\alpha_5,\ \alpha_3=0,\ \alpha_2\neq0,$ then choosing $x=\frac{\alpha_2}{\alpha_5},$ we have the representative  $\left\langle\nabla_2+2\nabla_4+\nabla_5\right\rangle;$
    \item $\alpha_4=2\alpha_5,\ \alpha_3\neq0,$ then choosing $x=\frac{\alpha_3}{\alpha_5},$ we have the family of representatives $\left\langle\alpha\nabla_2+\nabla_3+ 2\nabla_4+ \nabla_5\right\rangle;$
    \item $\alpha_4\neq2\alpha_5,\ \alpha_2=0,$ then choosing $x=1, \ z=-\frac{\alpha _3}{\alpha _4-2 \alpha _5},$ we have the family of representatives $\left\langle\alpha\nabla_4+ \nabla_5\right\rangle_{\alpha\neq2};$
    \item $\alpha_4\neq2\alpha_5,\ \alpha_2\neq0,$ then choosing $x=\frac{\alpha_2}{\alpha_5}, \ z=-\frac{\alpha_2\alpha _3}{(\alpha_4-2\alpha_5) \alpha_5},$ we have the family of  representatives $\left\langle\nabla_2+\alpha\nabla_4+ \nabla_5\right\rangle_{\alpha\neq2}.$
\end{enumerate}

Summarizing all cases, we have the following distinct orbits
\begin{longtable} {llll}
$\langle\alpha\nabla_2+\nabla_3+2\nabla_4+\nabla_5\rangle,$ & $\langle\alpha\nabla_4+\nabla_5\rangle,$ & $\langle\nabla_2+\alpha\nabla_4+\nabla_5\rangle,$\\
\end{longtable} 
which gives the following new algebras (see section \ref{secteoA}):

\begin{center}
${\rm B}_{122}^{\alpha},$
${\rm B}_{123}^{\alpha},$
${\rm B}_{124}^{\alpha}.$

\end{center}

\subsubsection{Central extensions of ${\mathcal B}_{11}^4$}
Let us use the following notations:
	\begin{longtable}{lllllll} 
$\nabla_1 = [\Delta_{11}],$& $\nabla_2 = [\Delta_{21}],$ & $\nabla_3=[\Delta_{22}],$ &
$\nabla_4=[\Delta_{32}],$ &
$\nabla_5=[\Delta_{14}+\Delta_{33}+\Delta_{42}].$
\end{longtable}	
	
Take $\theta=\sum\limits_{i=1}^{5}\alpha_i\nabla_i\in {\rm H^2}({\mathcal B}_{11}^4).$
The automorphism group of ${\mathcal B}_{11}^4$ consists of invertible matrices of the form
$$\phi=
	\begin{pmatrix}
	x &  0  & 0 & 0\\
	y &  x^2  & 0 & 0\\
	z &  0  & x^2 & 0\\
    u &  x(y+z)  & v & x^3
	\end{pmatrix}.
	$$
	Since
	$$
	\phi^T\begin{pmatrix}
	\alpha_1 & 0 & 0  & \alpha_5\\
	\alpha_2 & \alpha_3  & 0 & 0\\
	0 &	\alpha_4 & \alpha_5 & 0\\
	0 & \alpha_5  & 0 & 0
	\end{pmatrix} \phi=\begin{pmatrix}
	\alpha_1^* & \alpha^* & \alpha^{**}  & \alpha_5^*\\
	\alpha_2^* & \alpha_3^*  & 0 & 0\\
	0 &	\alpha_4^*+\alpha^{**} & \alpha_5^* & 0\\
	0 & \alpha_5^*  & 0 & 0
	\end{pmatrix},
$$
we have that the action of ${\rm Aut} ({\mathcal B}_{11})$ on the subspace
$\langle \sum\limits_{i=1}^{5}\alpha_i\nabla_i  \rangle$
is given by
$\langle \sum\limits_{i=1}^{5}\alpha_i^{*}\nabla_i\rangle,$
where
\begin{longtable}{lcllcllcl}
$\alpha^*_1$&$=$&$x \left(x \alpha _1+y \alpha _5\right),$ &
$\alpha_2^*$&$=$&$x^2 \alpha _2,$\\
$\alpha^*_3$&$=$&$x \left(x \alpha _3+z \alpha _5\right),$ &
$\alpha^*_4$&$=$&$x^3 \alpha _4,$&
$\alpha_5^*$&$=$&$x^4 \alpha_5.$\\
\end{longtable}
We are interested only in the cases with $\alpha_5\neq 0.$ Choosing $y=-\frac{x \alpha _1}{\alpha _5},\ z=-\frac{x \alpha _3}{\alpha _5},$  we have $\alpha_1^*=\alpha_3^*=0.$
\begin{enumerate}
    \item $\alpha_4=\alpha_2=0,$ then we have the representative $\left\langle\nabla_5\right\rangle;$
    \item $\alpha_4=0,\ \alpha_2\neq0,$ then choosing $x=\sqrt{\frac{\alpha_2}{\alpha_5}},$ we have the representative  $\left\langle\nabla_2+\nabla_5\right\rangle;$
    \item $\alpha_4\neq0,$ then choosing $x=\frac{\alpha_4}{\alpha_5},$ we have the family of representatives  $\left\langle\alpha\nabla_2+\nabla_4+\nabla_5\right\rangle.$
\end{enumerate}

Summarizing all cases, we have the following distinct orbits
\begin{center}  
$\langle\nabla_5\rangle,$   $\langle\nabla_2+\nabla_5\rangle,$  $\langle\alpha\nabla_2+\nabla_4+\nabla_5\rangle,$\\
\end{center}
which gives the following new algebras (see section \ref{secteoA}):

\begin{center}

${\rm B}_{125},$
${\rm B}_{126},$
${\rm B}_{127}.$

\end{center}

\subsubsection{Central extensions of ${\mathcal B}_{16}^4$}
Let us use the following notations:
	\begin{longtable}{lllllll} 
$\nabla_1 = [\Delta_{11}],$& $\nabla_2 = [\Delta_{21}],$ & $\nabla_3=[\Delta_{22}],$ &
$\nabla_4=[\Delta_{14}+\Delta_{23}],$ &
$\nabla_5=[\Delta_{32}].$
\end{longtable}	
	
Take $\theta=\sum\limits_{i=1}^{5}\alpha_i\nabla_i\in {\rm H^2}({\mathcal B}_{16}^4).$
The automorphism group of ${\mathcal B}_{16}^4$ consists of invertible matrices of the form
$$\phi=
	\begin{pmatrix}
	x &  0  & 0 & 0\\
	0 &  x^2  & 0 & 0\\
	0 &  y  & x^3 & 0\\
    u &  v  & xy & x^4
	\end{pmatrix}.
	$$
	Since
	$$
	\phi^T\begin{pmatrix}
	\alpha_1 & 0 & 0  & \alpha_4\\
	\alpha_2 & \alpha_3  & \alpha_4 & 0\\
	0 &	\alpha_5 & 0 & 0\\
	0 & 0  & 0 & 0
	\end{pmatrix} \phi=\begin{pmatrix}
    \alpha_1^* & \alpha^* & \alpha^{**}  & \alpha_4^*\\
	\alpha_2^* & \alpha_3^*+\alpha^{**}  & \alpha_4^* & 0\\
	0 &	\alpha_5^* & 0 & 0\\
	0 & 0  & 0 & 0
	\end{pmatrix},
$$
we have that the action of ${\rm Aut} ({\mathcal B}_{16}^4)$ on the subspace
$\langle \sum\limits_{i=1}^{5}\alpha_i\nabla_i  \rangle$
is given by
$\langle \sum\limits_{i=1}^{5}\alpha_i^{*}\nabla_i\rangle,$
where
\begin{longtable}{lcllcllcl}
$\alpha^*_1$&$=$&$x(x\alpha_1+u\alpha_4),$ &
$\alpha_2^*$&$=$&$x^3 \alpha_2,$\\
$\alpha^*_3$&$=$&$x^2(x^2 \alpha _3+  y \alpha _5),$ &
$\alpha^*_4$&$=$&$x^5 \alpha _4,$ &
$\alpha_5^*$&$=$&$x^5 \alpha_5.$
\end{longtable}
We are interested only in the cases with $\alpha_4\neq 0.$ Choosing $u=-\frac{x\alpha_1}{\alpha_4},$ we have $\alpha_1^*=0.$
\begin{enumerate}
    \item $\alpha_5=\alpha_3=\alpha_2=0,$ then we have the representative  $\left\langle\nabla_4\right\rangle;$
    \item $\alpha_5=\alpha_3=0,\ \alpha_2\neq0,$ then choosing $x=\sqrt{\frac{\alpha_2}{\alpha_4}},$ we have the representative  $\left\langle\nabla_2+\nabla_4\right\rangle;$
    \item $\alpha_5=0,\ \alpha_3\neq0,$ then choosing $x=\frac{\alpha_3}{\alpha_4},$ we have the family of representatives  $\left\langle\alpha\nabla_2+\nabla_3+ \nabla_4\right\rangle;$
    \item $\alpha_5\neq0,\ \alpha_2=0,$ then choosing $x=1,\ y=-\frac{\alpha_3}{\alpha_5},$ we have the family of representatives  $\left\langle\nabla_4+\alpha\nabla_5\right\rangle_{\alpha\neq0};$
    \item $\alpha_5\neq0,\ \alpha_2\neq0,$ then choosing $x=\sqrt{\frac{\alpha_2}{\alpha_4}},\  y=-\frac{\alpha_2\alpha_3}{\alpha_4\alpha_5},$ we have the family of representatives  $\left\langle\nabla_2+\nabla_4+ \alpha\nabla_5\right\rangle_{\alpha\neq0}.$
\end{enumerate}

Summarizing all cases, we have the following distinct orbits
\begin{center}  
$\langle\nabla_4+\alpha\nabla_5\rangle,$   $\langle\nabla_2+\nabla_4+\alpha\nabla_5\rangle,$  $\langle\alpha\nabla_2+\nabla_3+\nabla_4\rangle,$ 
\end{center}
which gives the following new algebras (see section \ref{secteoA}):

\begin{center}

${\rm B}_{128}^{\alpha},$
${\rm B}_{129}^{\alpha},$
${\rm B}_{130}^{\alpha}.$

\end{center}

\subsubsection{Central extensions of ${\mathcal B}_{17}^4$}
Let us use the following notations:
	\begin{longtable}{lllllll} 
$\nabla_1 = [\Delta_{11}],$& $\nabla_2 = [\Delta_{14}],$ & $\nabla_3=[\Delta_{21}],$ &
$\nabla_4=[\Delta_{22}],$ &
$\nabla_5=[\Delta_{32}].$
\end{longtable}	
	
Take $\theta=\sum\limits_{i=1}^{5}\alpha_i\nabla_i\in {\rm H^2}({\mathcal B}_{17}^4).$
The automorphism group of ${\mathcal B}_{17}^4$ consists of invertible matrices of the form
$$\phi=
	\begin{pmatrix}
	x &  0  & 0 & 0\\
	0 &  y  & 0 & 0\\
	0 &  z  & xy & 0\\
    u &  v  & xz & x^2y
	\end{pmatrix}.
	$$
	Since
	$$
	\phi^T\begin{pmatrix}
	\alpha_1 & 0 & 0  & \alpha_2\\
	\alpha_3 & \alpha_4  & 0 & 0\\
	0 &	\alpha_5 & 0 & 0\\
	0 & 0  & 0 & 0
	\end{pmatrix} \phi=\begin{pmatrix}
    \alpha_1^* & \alpha^* & \alpha^{**}  & \alpha_2^*\\
	\alpha_3^* & \alpha_4^*  & 0 & 0\\
	0 &	\alpha_5^* & 0 & 0\\
	0 & 0  & 0 & 0
	\end{pmatrix},
$$
we have that the action of ${\rm Aut} ({\mathcal B}_{17}^4)$ on the subspace
$\langle \sum\limits_{i=1}^{5}\alpha_i\nabla_i  \rangle$
is given by
$\langle \sum\limits_{i=1}^{5}\alpha_i^{*}\nabla_i\rangle,$
where
\begin{longtable}{lcllcllcl}
$\alpha^*_1$&$=$&$x \left(x \alpha _1+u \alpha _2\right),$ &
$\alpha_2^*$&$=$&$x^3y \alpha _2,$&
$\alpha^*_3$&$=$&$x y \alpha _3,$ \\
$\alpha^*_4$&$=$&$y \left(y \alpha _4+z \alpha _5\right),$&
$\alpha_5^*$&$=$&$xy^2 \alpha_5.$
\end{longtable}
We are interested only in the cases with $\alpha_2\neq 0.$ Choosing $v=-\frac{x \alpha _1}{\alpha _2},$ we have $\alpha_1^*=0.$
\begin{enumerate}
    \item $\alpha_5=\alpha_4=\alpha_3=0,$ then we have the representative  $\left\langle\nabla_2\right\rangle;$
    \item $\alpha_5=\alpha_4=0,\ \alpha_3\neq0,$ then choosing $x=\sqrt{\frac{\alpha_3}{\alpha_2}}, \ y=1,$ we have the representative  $\left\langle\nabla_2+\nabla_3\right\rangle;$
    \item $\alpha_5=0,\ \alpha_4\neq0,\ \alpha_3=0,$ then choosing $x=1,\ y=\frac{\alpha_2}{\alpha_4},$ we have the representative $\left\langle\nabla_2+\nabla_4\right\rangle;$
    \item $\alpha_5=0,\ \alpha_4\neq0,\ \alpha_3\neq0,$ then choosing $x=\sqrt{\frac{\alpha_3}{\alpha_2}},\ y=\frac{\alpha_3}{\alpha_4}\sqrt{\frac{\alpha_3}{\alpha_2}},$ we have the representative $\left\langle\nabla_2+\nabla_3+ \nabla_4\right\rangle;$
    \item $\alpha_5\neq0,\ \alpha_3=0,$ then choosing $x=1,\ y=\frac{\alpha_2}{\alpha_5},\ z=-\frac{\alpha_2\alpha_4}{\alpha_5^2},$ we have the representative $\left\langle\nabla_2+\nabla_5\right\rangle;$
    \item $\alpha_5\neq0,\ \alpha_3\neq0,$ then choosing $x=\sqrt{\frac{\alpha_3}{\alpha_2}},\ y=\frac{\alpha_3 } {\alpha_5},\ z=-\frac{\alpha_3\alpha_4}{\alpha_5^2},$ we have the representative $\left\langle\nabla_2+\nabla_3+ \nabla_5\right\rangle.$
\end{enumerate}

Summarizing all cases, we have the following distinct orbits
\begin{center}  
$\langle\nabla_2\rangle,$  
$\langle\nabla_2+\nabla_4\rangle,$  
$\langle\nabla_2+\nabla_3\rangle,$   $\langle\nabla_2+\nabla_3+\nabla_4\rangle,$  
$\langle\nabla_2+\nabla_5\rangle,$   $\langle\nabla_2+\nabla_3+\nabla_5\rangle,$  
\end{center}
which gives the following new algebras (see section \ref{secteoA}):

\begin{center}

${\rm B}_{131},$
${\rm B}_{132},$
${\rm B}_{133},$
${\rm B}_{134},$
${\rm B}_{135},$
${\rm B}_{136}.$

\end{center}

\subsubsection{Central extensions of ${\mathcal B}_{18}^4$}
Let us use the following notations:
	\begin{longtable}{lllllll} 
$\nabla_1= [\Delta_{13}],$& $\nabla_2 = [\Delta_{21}],$ & $\nabla_3=[\Delta_{22}],$ &
$\nabla_4=[\Delta_{32}],$ &
$\nabla_5=[\Delta_{31}+\Delta_{42}].$
\end{longtable}	
	
Take $\theta=\sum\limits_{i=1}^{5}\alpha_i\nabla_i\in {\rm H^2}({\mathcal B}_{18}^4).$
The automorphism group of ${\mathcal B}_{18}^4$ consists of invertible matrices of the form
$$\phi=
	\begin{pmatrix}
	x^2 &  0  & 0 & 0\\
	0 &  x  & 0 & 0\\
	y &  0  & x^2 & 0\\
    u &  z  & xy & x^4
	\end{pmatrix}.
	$$
	Since
	$$
	\phi^T\begin{pmatrix}
	0 & 0 & \alpha_1  & 0\\
	\alpha_2 & \alpha_3 & 0 & 0\\
	\alpha_5 &	\alpha_4 & 0 & 0\\
	0 & \alpha_5  & 0 & 0
	\end{pmatrix} \phi=\begin{pmatrix}
    \alpha^{**} & \alpha^* & \alpha_1^*  & 0\\
	\alpha_2^* & \alpha_3^* & 0 & 0\\
	\alpha_5^* & \alpha_4^*+\alpha^{**} & 0 & 0\\
	0 & \alpha_5^*  & 0 & 0
	\end{pmatrix},
$$
we have that the action of ${\rm Aut} ({\mathcal B}_{18}^4)$ on the subspace
$\langle \sum\limits_{i=1}^{5}\alpha_i\nabla_i  \rangle$
is given by
$\langle \sum\limits_{i=1}^{5}\alpha_i^{*}\nabla_i\rangle,$
where
\begin{longtable}{lcllcllcl}
$\alpha^*_1$&$=$&$x^5 \alpha _1,$ &
$\alpha_2^*$&$=$&$x^3 \alpha _2,$ &
$\alpha^*_3$&$=$&$x \left(x \alpha _3+z \alpha _5\right),$ \\
$\alpha^*_4$&$=$&$-x^2 y \alpha _1+x^4 \alpha _4,$ &
$\alpha_5^*$&$=$&$x^5 \alpha _5.$ 
\end{longtable}
We are interested only in the cases with $\alpha_5\neq 0.$ Choosing $z=-\frac{x \alpha _3}{\alpha _5},$ we have $\alpha_3^*=0.$
\begin{enumerate}
    \item $\alpha_1=\alpha_4=\alpha_2=0,$ then we have the representative $\left\langle\nabla_5\right\rangle;$
    \item $\alpha_1=\alpha_4=0,\ \alpha_2\neq0,$ then choosing $x=\sqrt{\frac{\alpha_2}{\alpha_5}},$ we have the representative  $\left\langle\nabla_2+\nabla_5\right\rangle;$
    \item $\alpha_1=0,\ \alpha_4\neq0,$ then choosing $x=\frac{\alpha_4}{\alpha_5},$ we have the family of representatives $\left\langle\alpha\nabla_2+\nabla_4+\nabla_5\right\rangle;$
    \item $\alpha_1\neq0,\ \alpha_2=0,$ then choosing $x=1,\ y=\frac{\alpha_4}{\alpha_1},$ we have the representative $\left\langle\alpha\nabla_1+\nabla_5\right\rangle_{\alpha\neq0};$
    \item $\alpha_1\neq0,\ \alpha_2\neq0,$ then choosing $x=\sqrt{\frac{\alpha_2}{\alpha_5}},\ y=\frac{\alpha_2\alpha_4}{\alpha_1\alpha_5},$ we have the family of representatives $\left\langle\alpha\nabla_1+\nabla_2+ \nabla_5\right\rangle_{\alpha\neq0}.$
\end{enumerate}

Summarizing all cases, we have the following distinct orbits
\begin{center} 
$\langle\alpha\nabla_2+\nabla_4+\nabla_5\rangle,$   $\langle\alpha\nabla_1+\nabla_5\rangle,$  $\langle\alpha\nabla_1+\nabla_2+\nabla_5\rangle,$ 
\end{center}
which gives the following new algebras (see section \ref{secteoA}):
\begin{center}
${\rm B}_{137}^{\alpha},$
${\rm B}_{138}^{\alpha},$
${\rm B}_{139}^{\alpha}.$
\end{center}

\subsubsection{Central extensions of ${\mathcal B}_{19}^4$}
Let us use the following notations:
	\begin{longtable}{lllllll} 
$\nabla_1=[\Delta_{11}],$& $\nabla_2 = [\Delta_{13}],$ & $\nabla_3=[\Delta_{21}],$ &
$\nabla_4=[\Delta_{22}],$ &
$\nabla_5=[\Delta_{42}].$
\end{longtable}	
	
Take $\theta=\sum\limits_{i=1}^{5}\alpha_i\nabla_i\in {\rm H^2}({\mathcal B}_{19}^4).$
The automorphism group of ${\mathcal B}_{19}^4$ consists of invertible matrices of the form
$$\phi=
	\begin{pmatrix}
	x &  0  & 0 & 0\\
	0 &  y  & 0 & 0\\
	z &  0  & xy & 0\\
    u &  v  & yz & xy^2
	\end{pmatrix}.
	$$
	Since
	$$
	\phi^T\begin{pmatrix}
	\alpha_1 & 0 & \alpha_2  & 0\\
	\alpha_3 & \alpha_4 & 0 & 0\\
	0 &	0 & 0 & 0\\
	0 & \alpha_5  & 0 & 0
	\end{pmatrix} \phi=\begin{pmatrix}
    \alpha_1^* & \alpha^* & \alpha_2^*  & 0\\
	\alpha_3^* & \alpha_4^* & 0 & 0\\
	0 &	\alpha^{**} & 0 & 0\\
	0 & \alpha_5^*  & 0 & 0
	\end{pmatrix},
$$
we have that the action of ${\rm Aut} ({\mathcal B}_{19}^4)$ on the subspace
$\langle \sum\limits_{i=1}^{5}\alpha_i\nabla_i  \rangle$
is given by
$\langle \sum\limits_{i=1}^{5}\alpha_i^{*}\nabla_i\rangle,$
where
\begin{longtable}{lcllcllcl}
$\alpha^*_1$&$=$&$x \left(x \alpha _1+z \alpha _2\right),$ &
$\alpha_2^*$&$=$&$x^2y \alpha _2,$ &
$\alpha^*_3$&$=$&$x y \alpha _3,$ \\
$\alpha^*_4$&$=$&$y \left(y \alpha _4+v \alpha _5\right),$&
$\alpha_5^*$&$=$&$x y^3 \alpha _5.$
\end{longtable}
We are interested only in the cases with $\alpha_5\neq 0.$ Choosing $v=-\frac{y \alpha _4}{\alpha _5},$ we have $\alpha_4^*=0.$
\begin{enumerate}
    \item $\alpha_2=\alpha_1=\alpha_3=0,$ then we have the representative $\left\langle\nabla_5\right\rangle;$
    \item $\alpha_2=\alpha_1=0,\ \alpha_3\neq0,$ then choosing $y=\sqrt{\frac{\alpha_3}{\alpha_5}},$ we have the representative  $\left\langle\nabla_3+\nabla_5\right\rangle;$
    \item $\alpha_2=0,\ \alpha_1\neq0,\ \alpha_3=0,$ then choosing $x=\frac{\alpha_5}{\alpha_1},\ y=1,$ we have the representative $\left\langle\nabla_1+\nabla_5\right\rangle;$
    \item $\alpha_2=0,\ \alpha_1\neq0,\ \alpha_3\neq0,$ then choosing $x=\frac{\alpha_3\sqrt{\alpha_3}}{\alpha_1\sqrt{\alpha_5}},\ y=\sqrt{\frac{\alpha_3}{\alpha_5}},$ we have the representative $\left\langle\nabla_1+\nabla_3+\nabla_5\right\rangle;$
    \item $\alpha_2\neq0,\ \alpha_3=0,$ then choosing $x=1,\ y=\sqrt{\frac{\alpha_2}{\alpha_5}},\ z=-\frac{\alpha_1}{\alpha_2},$ we have the representative $\left\langle\nabla_2+\nabla_5\right\rangle;$
    \item $\alpha_2\neq0,\ \alpha_3\neq0,$ then choosing $x=\frac{\alpha_3}{\alpha_2},\  y=\sqrt{\frac{\alpha_3}{\alpha_5}},\ z=-\frac{\alpha_1\alpha_3}{\alpha_2^2},$ we have the representative $\left\langle\nabla_2+\nabla_3+ \nabla_5\right\rangle.$
\end{enumerate}

Summarizing all cases, we have the following distinct orbits
\begin{center}  
$\langle\nabla_5\rangle,$  
$\langle\nabla_3+\nabla_5\rangle,$   $\langle\nabla_1+\nabla_5\rangle,$ $\langle\nabla_1+\nabla_3+\nabla_5\rangle,$ 
$\langle\nabla_2+\nabla_5\rangle,$ $\langle\nabla_2+\nabla_3+\nabla_5\rangle,$ 
\end{center}
which gives the following new algebras (see section \ref{secteoA}):

\begin{center}

${\rm B}_{140},$
${\rm B}_{141},$
${\rm B}_{142},$
${\rm B}_{143},$
${\rm B}_{144},$
${\rm B}_{145}.$

\end{center}

\subsection{$2$-dimensional central extensions of $3$-dimensional  nilpotent bicommutative algebras}

 \subsubsection{The description of second cohomology spaces of  $3$-dimensional nilpotent bicommutative algebras:}

\
In the following table we give the description of the second cohomology space of   two-generated $3$-dimensional nilpotent bicommutative algebras

\begin{longtable}{ll llllll}

%\multicolumn{8}{c}{{\bf The list of nilpotent 3-dimensional bicommutative algebras}}  \\
\hline
 
{${\mathcal B}_{01}^{3*}$} &$:$ &  $e_1e_1 = e_2$ &&&&\\ 
\multicolumn{8}{l}{
${\rm H}_{com}^2({\mathcal B}_{01}^{3*})= 
\Big\langle 
 [\Delta_{ 1 2}+\Delta_{21}], 
[\Delta_{1 3}+\Delta_{ 31}], 
[\Delta_{3 3}]  

\Big\rangle $} \\
\multicolumn{8}{l}{
${\rm H}_{bicom}^2({\mathcal B}_{01}^{3*})=  {\rm H}_{com}^2({\mathcal B}_{01}^*) \oplus \Big\langle 
 [\Delta_{21}],[\Delta_{31}] 
\Big\rangle  $}\\
 
\hline
{${\mathcal B}_{02}^{3*}$} &$:$ & $e_1e_1 = e_3$& $e_2e_2=e_3$  &&&  \\ 

\multicolumn{8}{l}{
${\rm H}_{bicom}^2({\mathcal B}_{02}^{3*})=\Big\langle [\Delta_{12}],[\Delta_{21}],[\Delta_{22}]
 \Big\rangle $}\\
\hline {${\mathcal B}_{03}^{3*}$} &$:$ &  $e_1e_2=e_3$ & $e_2e_1=-e_3$ &&& \\ 
\multicolumn{8}{l}{
${\rm H}^2({\mathcal B}_{03}^{3*})=
\Big\langle [\Delta_{11}],[\Delta_{12}],[\Delta_{22}]
\Big\rangle $}\\
\hline
${\mathcal B}_{04}^{3*}(\alpha\neq  0)$ &$:$ & $e_1e_1=\alpha e_3$ & $e_2e_1=e_3$& $e_2e_2=e_3$  &&\\  
\multicolumn{8}{l}{
${\rm H}^2({\mathcal B}_{04}^{3*}(\alpha\neq  0))=
\Big\langle  [\Delta_{11}],[\Delta_{12}],[\Delta_{21}]  \Big\rangle $ }\\

\hline
${\mathcal B}_{04}^{3*}(0)$ &$:$ & $e_1e_2=e_3$ & &&\\ 

\multicolumn{8}{l}{
${\rm H}^2({\mathcal B}_{04}^*(0))= \Big\langle [\Delta_{11}], [\Delta_{13}], [\Delta_{21}], [\Delta_{22}], [\Delta_{32}] \Big\rangle $}\\

%\hline
%${\mathcal B}_{01}^3$ &$:$ & $e_1e_1 =e_2$& $e_2e_1 =e_3$   &&&\\ 

%\multicolumn{8}{l}{
%${\rm H}^2({\mathcal B}_{01}^3)=\Big\langle [\Delta_{12}], [\Delta_{31}] \Big\rangle $}\\

%\hline
%{${\mathcal B}_{03}^3(\alpha)$} &$:$ & $e_1e_1 =e_2$ & $e_1e_2 = e_3$ &  $e_2e_1 =\alpha e_3$ &&\\ 

 %\multicolumn{8}{l}{
%${\rm H}^2({\mathcal B}_{03}^3(\alpha))=
%\Big\langle  [\Delta_{21}],[\Delta_{13}+\alpha\Delta_{22}+\alpha\Delta_{31}].
%\Big\rangle $}\\

\hline
\end{longtable}

\subsubsection{Central extensions of ${\mathcal B}^{3*}_{01}$}
Let us use the following notations:
\[
\nabla_1=[\Delta_{12}+\Delta_{21}], \quad \nabla_2=[\Delta_{13}+\Delta_{31}], \quad \nabla_3=[\Delta_{21}], \quad \nabla_4=[\Delta_{31}], \quad  \nabla_5=[\Delta_{33}]. \]

The automorphism group of ${\mathcal B}^{3*}_{01}$ consists of invertible matrices of the form

\[\phi=\begin{pmatrix}
x & 0 & 0\\
u & x^2 & w\\
z & 0 & y
\end{pmatrix}. \]

Since

\[ \phi^T\begin{pmatrix}
0 & \alpha_1 & \alpha_2\\
\alpha_1+\alpha_3 & 0 & 0\\
\alpha_2+\alpha_4 & 0 & \alpha_5
\end{pmatrix}\phi =
\begin{pmatrix}
\alpha^* & \alpha^*_1 & \alpha^*_2 \\
\alpha^*_1+ \alpha^*_3& 0 & 0 \\
\alpha^*_2+\alpha^*_4 & 0 & \alpha^*_5
\end{pmatrix},
\]
the action of $\operatorname{Aut} (\mathcal{B}_{01}^{3*})$ on subspace
$\Big\langle \sum\limits_{i=1}^5 \alpha_i\nabla_i \Big\rangle$ is given by
$\Big\langle \sum\limits_{i=1}^5 \alpha_i^*\nabla_i \Big\rangle,$
where
\begin{longtable}{lcllcllcl}
$\alpha^*_1$&$=$&$x^3 \alpha_1,$ &
$\alpha^*_2$&$=$&$w x \alpha _1+x y \alpha _2+y z \alpha _5,$ &
$\alpha^*_3$&$=$&$x^3 \alpha_3,$\\
$\alpha^*_4$&$=$&$ x(w\alpha_3+y \alpha_4),$&
$\alpha^*_5$&$=$&$y^2 \alpha_5.$
\end{longtable}

We are interested  only in $2$-dimensional central extensions and consider the vector space generated by the following two cocycles:
\begin{center}
$\theta_1=\alpha_1\nabla_1+\alpha_2\nabla_2+\alpha_3\nabla_3+\alpha_4\nabla_4+\alpha_5\nabla_5 \ \ \text{and} \ \  \theta_2=\beta_1\nabla_1+\beta_2\nabla_2+\beta_4\nabla_4+\beta_5\nabla_5.$
\end{center}
Our aim is to find only central extensions with $(\alpha_3,\alpha_4, \beta_3,\beta_4)\neq 0.$
Hence, we have the following cases.

\begin{enumerate}
    \item $\alpha_3\neq0,$ then we have
    \begin{longtable}{lcllcl}
    $\alpha^*_1$&$=$&$x^3 \alpha_1,$ &   $\beta^*_1$&$=$&$x^3\beta_1,$ \\
    $\alpha^*_2$&$=$&$w x \alpha _1+x y \alpha _2+y z \alpha _5,$ &   $\beta^*_2$&$=$&$w x \beta_1+xy \beta_2+y z \beta_5,$\\ 
    $\alpha^*_3$&$=$&$x^3 \alpha_3,$   & $\beta^*_3$&$=$&$0,$\\ 
    $\alpha_4^*$&$=$&$x(w\alpha_3+y \alpha_4),$   & $\beta_4^*$&$=$&$x y\beta_4,$ \\
    $\alpha_5^*$&$=$&$y^2\alpha_5,$ &   $\beta_5^*$&$=$&$y^2\beta_5.$ \\
    \end{longtable}
    \begin{enumerate}
        \item $\beta_5\neq0,$ then we can suppose $\alpha_5=0$ and choosing 
        $w=-\frac{y\alpha_4}{\alpha_3},$ 
        $z=-\frac{x(\alpha_4\beta_1-\alpha_3\beta_2)}{\alpha_3\beta_5}$, 
        we have $\alpha_4^*=\beta_2^*=0.$ Thus, we can assume $\alpha_4=\beta_2=0$ and consider following subcases:
        \begin{enumerate}
            \item $\alpha_2=\beta_4=\beta_1=0,$ then we have the family of  representatives $ \left\langle \alpha\nabla_1+ \nabla_3,\nabla_5 \right\rangle;$ 
 
 \item $\alpha_2=\beta_4=0,\ \beta_1\neq0,$ then choosing $x=\sqrt[3]{{\beta_5}{\beta_1^{-1}}},\ y=1,$ we have the family of representatives $ \left\langle  \alpha \nabla_1+\nabla_3,\nabla_1+\nabla_5 \right\rangle
 ;$ 
            \item $\alpha_2=0,\ \beta_4\neq0,\ \beta_1=0$ then choosing $x={\beta_5}{\beta_4^{-1}}, \ y=1,$ we have the family of representatives $ \left\langle  \alpha \nabla_1+\nabla_3,\nabla_4+\nabla_5 \right\rangle;$ 
            
            \item $\alpha_2=0,\ \beta_4\neq0,\ \beta_1\neq0,$ then choosing $x={\beta_4^2}{\beta_1^{-1}\beta_5^{-1}},\ y={\beta_4^3}{\beta_1^{-1}\beta_5^{-2}},$ we have the family of representatives $ \left\langle  \alpha \nabla_1+\nabla_3,\nabla_1+\nabla_4+\nabla_5 \right\rangle;$ 
            
            \item $\alpha_2\neq0,\ \beta_4=\beta_1=0,$ then choosing $x=1,\ y={\alpha_3}{\alpha_2^{-1}},$ we have the family of  representatives $ \left\langle  \alpha \nabla_1+ \nabla_2+\nabla_3,\nabla_5\right\rangle;$ 
            
            \item $\alpha_2\neq0,\ \beta_4=0,\ \beta_1\neq0,$ then choosing 
            $x={\alpha_2^2\beta_1}{\alpha_3^{-2}\beta_5^{-1}},$ 
            $y={\alpha_2^3\beta_1^2}{\alpha_3^{-3}\beta_5^{-2}},$ we have the family of representatives $\left\langle  \alpha \nabla_1+ \nabla_2 +\nabla_3,\nabla_1+\nabla_5\right\rangle;$ 

\item $\alpha_2\neq0,\ \beta_4\neq0,$ then choosing $x={\alpha_2\beta_4}{\alpha_3^{-1}\beta_5^{-1}},$ $y={\alpha_2\beta_4^2}{\alpha_3^{-1}\beta_5^{-2}},$ 
we have the family of  representatives $ \left\langle  \alpha \nabla_1+ \nabla_2+ \nabla_3,\beta \nabla_1+ \nabla_4+ \nabla_5\right\rangle;$ 
           \end{enumerate}

\item $\beta_5=0, \beta_4\neq0.$ 
        \begin{enumerate}
            \item $\alpha_5=\beta_1=0, $ $ \alpha_1\beta_4\neq\alpha_3\beta_2,$ then choosing $y=1,\ w=\frac{\alpha_4\beta_2-\alpha_2\beta_4}{\alpha_1\beta_4-\alpha_3\beta_2},$ we have
            the family of representatives $ \left\langle \alpha\nabla_1+ \nabla_3,\beta\nabla_2+\nabla_4 \right\rangle_{\alpha\neq \beta};$ 

\item $\alpha_5=\beta_1=0,$ $\alpha_1\beta_4=\alpha_3\beta_2,$ $\alpha_2\alpha_3=\alpha_1\alpha_4,$ 
then choosing $y=1,\ w=-{\alpha_4}{\alpha_3^{-1}},$ we have the family of representatives $\left\langle \alpha\nabla_1+ \nabla_3,\alpha\nabla_2+\nabla_4 \right\rangle;$ 

            \item $\alpha_5=\beta_1=0,$ $ \alpha_1\beta_4=\alpha_3\beta_2,$ $ \alpha_2\alpha_3\neq\alpha_1\alpha_4$ then choosing 
            $x=\alpha_4 \beta_2-\alpha_2 \beta_4,$
            $y=-\alpha_3 \beta_4 (\alpha_4 \beta_2-\alpha_2 \beta_4),$ and 
            $w=\alpha_4 \beta_4 (\alpha_4 \beta_2-\alpha_2 \beta_4),$ we have the family of  representatives 
            $\left\langle \alpha\nabla_1+\nabla_2+\nabla_3,\alpha\nabla_2+\nabla_4 \right\rangle;$ 
            
            \item $\alpha_5=0,$ $\beta_1\neq0,$ then choosing 
            \begin{center}
                $x=1, \ y=\frac{\alpha_3}{\beta_4},$ 
            $w=\frac{(\alpha_3\beta_2+ \alpha_4 \beta_1-\alpha_1\beta_4)- \sqrt{(\alpha_3\beta_2+ \alpha_4 \beta_1-\alpha_1\beta_4)^2- 4\alpha_3 \beta_1(\alpha_4\beta_2-\alpha_2 \beta_4)}} {\alpha_3\beta_4},$
            \end{center} we have the family of representatives 
            $ \left\langle \alpha\nabla_1+ \nabla_3,\nabla_1+ \beta \nabla_2 +\nabla_4 \right\rangle;$
            
            \item $\alpha_5\neq0,$ $\beta_1=0,$ then choosing 
      \begin{center}      $x=\alpha_5,$
            $y=-\sqrt{\alpha_3}\alpha_5,$
            $z=0$
            and
            $w=-\frac{\sqrt{\alpha_3}\alpha_5(\alpha_4\beta_2-\alpha_2\beta_4)}{\alpha_3\beta_2-\alpha_1\beta_4},$\end{center}  
            we have the family of representatives 
            $\left\langle \alpha\nabla_1+\nabla_3+\nabla_5,\beta \nabla_2 +\nabla_4 \right\rangle;$ 
            
            \item $\alpha_5\neq0,$ $\beta_1\neq0,$ then choosing 
        \begin{center}
        $x=\frac{\alpha_3\beta_4^2}{\alpha_5\beta_1^2}, $ 
            $y=\frac{\alpha_3^2\beta_4^3}{\alpha_5^2\beta_1^3},$ 
            $w=-\frac{\alpha_3^2 \beta_2 \beta_4^3}{\alpha_5^2 \beta_1^4}$
            and
            $z=\frac{\alpha_3 (\alpha_1 \beta_2-\alpha_2 \beta_1) \beta_4^2}{\alpha_5^2 \beta_1^3},$
            \end{center}we have the family of representatives   $\left\langle \alpha\nabla_1+ \nabla_3+ \nabla_5, \nabla_1+\nabla_4\right \rangle;$ 
        \end{enumerate}
        
        \item $\beta_5=0, \beta_4=0,\ \beta_1\neq0,$ then we can suppose $\alpha_1=0$ and consider following subcases:
        \begin{enumerate}
            \item $\alpha_5=0,$ then choosing $w=-\frac{y\beta_2}{\beta_1}$, we have $\beta_2^*=0.$ 
            \begin{enumerate}
                \item if $\alpha_2=\alpha_4=0,$ then we have a   split algebra;
                \item if $\alpha_2=0,$ $ \alpha_4\neq0,$ then choosing $x=1,$ $ y=\frac{\alpha_3}{\alpha_4},$ we have the representative $\left\langle \nabla_3+\nabla_4, \nabla_1\right \rangle;$

\item if $\alpha_2\neq0,$ then choosing $x=1,$ $ y=\frac{\alpha_3}{\alpha_2},$ we have the family of representatives  $\left\langle \nabla_2+\nabla_3+\alpha \nabla_4, \nabla_1\right \rangle;$
            \end{enumerate}

\item $\alpha_5\neq0, $ $ \beta_2=0,$ then choosing 
$x=\alpha_5 ,$
$y=\sqrt{\alpha_3}  \alpha_5,$
$z=-\alpha_2$ and
$w= 0,$
 we have the representative $\left\langle \nabla_3+ \nabla_5, \nabla_1\right \rangle;$

\item $\alpha_5\neq0, $ $ \beta_2\neq0,$ then choosing 
\begin{center}
    $x=\frac{\alpha_3 \beta_2^2}{\alpha_5 \beta_1^2},$
$y=\frac{\alpha_3^2 \beta_2^3}{\alpha_5^2 \beta_1^3},$
$z=\frac{\alpha_3 \beta_2^2 (\alpha_2 \beta_1-\alpha_1 \beta_2)}{\alpha_5^2 \beta_1^3}$ and
$w=0,$
\end{center}
 we have the representative $\left\langle \nabla_3+\nabla_5, \nabla_1+\nabla_2\right \rangle.$

        \end{enumerate}
        \item $\beta_5=\beta_4=\beta_1=0, \beta_2\neq0,$ then we can suppose $\alpha_2=0$ and choosing $w=-\frac{y\alpha_4}{\alpha_3},$ we have $\alpha_4^*=0.$  Thus, we have following subcases:
        \begin{enumerate}

\item if $\alpha_5=0,$ then we have the family of representatives 
$\left\langle \alpha \nabla_1+ \nabla_3, \nabla_2\right \rangle;$
 
 \item if $\alpha_5\neq0,$ then choosing $x=1,\ y=\sqrt{{\alpha_3}{\alpha_5^{-1}}},$ we have the family of representatives $\left\langle \alpha \nabla_1+ \nabla_3+\nabla_5, \nabla_2\right\rangle.$
        \end{enumerate}
    \end{enumerate}

\item $\alpha_3=0,$ $\alpha_4\neq0,$ then we can suppose $\beta_4=0.$
    \begin{longtable}{lcllcl}
    $\alpha^*_1$&$=$&$x^3 \alpha_1,$   & $\beta^*_1$&$=$&$x^3\beta_1,$ \\
    $\alpha^*_2$&$=$&$w x \alpha _1+x y \alpha _2+y z \alpha _5,$ &   $\beta^*_2$&$=$&$w x \beta_1+xy \beta_2+y z \beta_5,$\\ 
    $\alpha^*_3$&$=$&$0,$ &  $\beta^*_3$&$=$&$0,$\\ 
    $\alpha_4^*$&$=$&$x y \alpha_4,$   & $\beta_4^*$&$=$&$0,$ \\
    $\alpha_5^*$&$=$&$y^2\alpha_5,$ & $\beta_5^*$&$=$&$y^2\beta_5.$ \\
    \end{longtable}
    \begin{enumerate}
        \item $\beta_5\neq0,$ then we can suppose $\alpha_5=0$ and choosing $z=-\frac{x(w\beta_1+y\beta_2)}{y\beta_5},$ we have $\beta_2^*=0.$ Thus, we have following subcases:
        \begin{enumerate}
            \item if $\beta_1=0,$ then $\alpha_1\neq0$ and  choosing $x=1,\ y=\frac{\alpha_1}{\alpha_4},\ w=-\frac{\alpha_2}{\alpha_4}$ we have the representative $\left\langle \nabla_1+ \nabla_4, \nabla_5\right\rangle;$
            
            \item if $\beta_1\neq0,$ $ \alpha_1=0$ then choosing $x=1,\ y=\sqrt{\frac{\beta_1} {\beta_5}},$ we have the family of representatives              $\left\langle \alpha \nabla_2+\nabla_4, \nabla_1+\nabla_5\right\rangle;$

\item if $\beta_1\neq0,$ $ \alpha_1\neq0$ then choosing $x=\frac{\alpha_4^2\beta_1} {\alpha_1^2\beta_5},$ 
$ y=\frac{\alpha_2^3 \beta_1^2}{\alpha_1^3\beta_5^2}, $ 
$ w=-\frac{\alpha_2^4\beta_1^2}{\alpha_1^4\beta_5^2},$ we have the representative  
            $\left\langle \nabla_1+ \nabla_4, \nabla_1+\nabla_5\right\rangle.$
        \end{enumerate}
        \item $\beta_5=0,$ $\beta_1\neq0,$ then we can suppose $\alpha_1=0$ and choosing $w=-\frac{y\beta_2}{\beta_1},$ we have $\beta_2^*=0.$ Thus, we have following subcases:
        \begin{enumerate}
            \item if $\alpha_5=0,$ then we have the family of representatives $\left\langle \alpha \nabla_2+ \nabla_4, \nabla_1\right\rangle;$
            \item if $\alpha_5\neq0,$ then choosing $x=1,\ y=\frac{\alpha_4}{\alpha_5},\ z=-\frac{\alpha_2}{\alpha_5},$ we have the representative 
            $\left\langle \nabla_4+\nabla_5, \nabla_1 \right\rangle.$
        \end{enumerate}
    \item $\beta_5=\beta_1=0,$ $ \beta_2\neq0,$ then we can suppose $\alpha_2=0$. Since in case of $\alpha_1 =0,$ we have a split extension, we can assume $\alpha_1  \neq 0,$  
        Thus, we have following subcases:
        \begin{enumerate}
            \item if $\alpha_5=0,$ then choosing $x=1,\ y=\frac{\alpha_1}{\alpha_4},$ we have the representative 
            $\left\langle \nabla_1+\nabla_4, \nabla_2 \right\rangle;$
            \item if $\alpha_5\neq0,$ then choosing $x=\frac{\alpha_4^2}{\alpha_1\alpha_5},\ y=\frac{\alpha_4^3}{\alpha_1\alpha_5},$ we have the representative 
            $\left\langle \nabla_1+\nabla_4+\nabla_5, \nabla_2 \right\rangle.$
        \end{enumerate}
    \end{enumerate}
\end{enumerate}

Now we have the following distinct orbits:
\begin{center}$ \left\langle \alpha\nabla_1+ \nabla_3,\nabla_5 \right\rangle,$ \ 
$ \left\langle  \alpha \nabla_1+\nabla_3,\nabla_1+\nabla_5 \right\rangle,$  \ 
$ \left\langle  \alpha \nabla_1+\nabla_3,\nabla_4+\nabla_5 \right\rangle,$ \
$ \left\langle  \alpha \nabla_1+\nabla_3,\nabla_1+\nabla_4+\nabla_5 \right\rangle,$ \ 
$ \left\langle  \alpha \nabla_1+ \nabla_2+\nabla_3,\nabla_5\right\rangle,$  \ 
$\left\langle  \alpha \nabla_1+ \nabla_2 +\nabla_3,\nabla_1+\nabla_5\right\rangle,$ \
$ \left\langle  \alpha \nabla_1+ \nabla_2+ \nabla_3,\beta \nabla_1+ \nabla_4+ \nabla_5\right\rangle,$ \ 
$ \left\langle \alpha\nabla_1+ \nabla_3,\beta\nabla_2+\nabla_4 \right\rangle,$ \ 
$\left\langle \alpha\nabla_1+\nabla_2+\nabla_3,\alpha\nabla_2+\nabla_4 \right\rangle,$  \ 
$ \left\langle \alpha\nabla_1+ \nabla_3,\nabla_1+ \beta \nabla_2 +\nabla_4 \right\rangle,$ \ 
$\left\langle \alpha\nabla_1+\nabla_3+\nabla_5,\beta \nabla_2 +\nabla_4 \right\rangle,$  \ 
$\left\langle \alpha\nabla_1+ \nabla_3+ \nabla_5, \nabla_1+\nabla_4\right \rangle,$  \ 
$\left\langle \nabla_1, \nabla_3+\nabla_4 \right \rangle,$  \ 
$\left\langle \nabla_1, \nabla_2+\nabla_3+\alpha \nabla_4\right \rangle,$  \ 
$\left\langle \nabla_1, \nabla_3+\nabla_5, \right \rangle,$ \ 
$\left\langle \nabla_1+\nabla_2, \nabla_3+\nabla_5, \right \rangle,$ \ 
$\left\langle \alpha \nabla_1+ \nabla_3, \nabla_2\right \rangle,$  \ 
$\left\langle \alpha \nabla_1+ \nabla_3+\nabla_5, \nabla_2\right\rangle,$ \ 
$\left\langle \nabla_1+ \nabla_4, \nabla_5\right\rangle,$ \  
$\left\langle \nabla_1+\nabla_5, \alpha \nabla_2+\nabla_4\right\rangle,$ \ 
$\left\langle \nabla_1+ \nabla_4, \nabla_1+\nabla_5\right\rangle,$ \ 
$\left\langle \nabla_1, \alpha \nabla_2+ \nabla_4\right\rangle,$ \ 
$\left\langle \nabla_1, \nabla_4+\nabla_5\right\rangle,$ \ 
$\left\langle \nabla_1+\nabla_4, \nabla_2 \right\rangle,$ \ 
$\left\langle \nabla_1+\nabla_4+\nabla_5, \nabla_2 \right\rangle.$
\end{center}

Hence, we have the following new $5$-dimensional nilpotent bicommutative algebras (see section \ref{secteoA}):

\begin{center}
${\rm B}_{146}^{\alpha},$ 
${\rm B}_{147}^{\alpha},$
${\rm B}_{148}^{\alpha},$
${\rm B}_{149}^{\alpha},$
${\rm B}_{150}^{\alpha},$
${\rm B}_{151}^{\alpha},$
${\rm B}_{152}^{\alpha,\beta},$
${\rm B}_{153}^{\alpha,\beta},$
${\rm B}_{154}^{\alpha},$
${\rm B}_{155}^{\alpha,\beta},$
${\rm B}_{156}^{\alpha,\beta},$
${\rm B}_{157}^{\alpha},$
${\rm B}_{158},$
${\rm B}_{159}^{\alpha},$
${\rm B}_{160},$
${\rm B}_{161},$

${\rm B}_{162}^{\alpha},$
${\rm B}_{163}^{\alpha},$
${\rm B}_{164},$
${\rm B}_{165}^{\alpha},$
${\rm B}_{166},$
${\rm B}_{167}^{\alpha},$
${\rm B}_{168},$
${\rm B}_{169},$
${\rm B}_{170}.$
\end{center}

\subsubsection{Central extensions of ${\mathcal B}^{3*}_{04}(0)$}

Let us use the following notations:
\[\nabla_1=[\Delta_{11}], \quad  \nabla_2=[\Delta_{13}], \quad \nabla_3=[\Delta_{21}], \quad \nabla_4=[\Delta_{22}], \quad
\nabla_5=[\Delta_{32}].\]

The automorphism group of ${\mathcal B}^{3*}_{04}(0)$ consists of invertible matrices of the form

\[\phi=\left(
                         \begin{array}{ccc}
                               x & 0 & 0   \\
                               0 & y & 0  \\
                               z & t & xy                               \end{array}\right)
                               .\]
Since
\[\phi^T\left(\begin{array}{ccc}
\alpha_ 1& 0 & \alpha_2  \\
\alpha_3 & \alpha_4 & 0  \\
0 & \alpha_5 & 0 \\
\end{array}\right)\phi=
\left(\begin{array}{ccc}
\alpha_1^*& \alpha^* & \alpha_2^*  \\
\alpha_3^* & \alpha_4^* & 0  \\
0 & \alpha_5^* & 0 \\
\end{array}\right),\]

the action of $\operatorname{Aut} ({\mathcal B}^{3*}_{04}(0))$ on the subspace
$\langle  \sum\limits_{i=1}^5\alpha_i \nabla_i \rangle$
is given by
$\langle  \sum\limits_{i=1}^5\alpha^*_i \nabla_i \rangle,$ where

\begin{longtable}{lcllcllcl}
$\alpha^*_1$&$=$&$x(x\alpha_1+z\alpha_2),$&
$\alpha^*_2$&$=$&$x^2y\alpha_2,$&
$\alpha^*_3$&$=$&$x y \alpha _3,$\\
$\alpha^*_4$&$=$&$y(y\alpha _4+t \alpha _5),$& 
$\alpha^*_5$&$=$&$xy^2\alpha_5.$
\end{longtable}

We are interested  only in $(\alpha_2,\alpha_5)\neq (0,0)$
  and consider the vector space generated by the following two cocycles:
$$ \theta_1=\alpha_1\nabla_1+\alpha_2\nabla_2+\alpha_3\nabla_3+\alpha_4\nabla_4+\alpha_5\nabla_5 \ \ \text{and} \ \  \theta_2=\beta_1\nabla_1+\beta_3\nabla_3+\beta_4\nabla_4+\beta_5\nabla_5.$$

\begin{enumerate}
\item $\alpha_2\neq0,$ then we have
\begin{longtable}{lcllcl}
$\alpha^*_1$&$=$&$x(x\alpha_1+z\alpha_2),$ & $\beta^*_1$&$=$&$x^2\beta_1,$ \\
$\alpha^*_2$&$=$&$x^2y\alpha_2,$ &  $\beta^*_2$&$=$&$0,$\\ 
$\alpha^*_3$&$=$&$x y \alpha _3,$ & $\beta^*_3$&$=$&$x y \beta_3,$\\ 
$\alpha_4^*$&$=$&$y(y\alpha _4+t \alpha _5),$ &  $\beta_4^*$&$=$&$y(y\beta_4+t\beta_5),$ \\
$\alpha_5^*$&$=$&$xy^2\alpha_5,$  & $\beta_5^*$&$=$&$x y^2\beta_5.$ \\
\end{longtable}
\begin{enumerate}
    \item $\beta_5\neq0,$ then we can suppose $\alpha_5^*=0$ and choosing $z=-\frac{x\alpha_1}{\alpha_2},\ t=-\frac{y\beta_4}{\beta_5}$, we have $\alpha_1^*=\beta_4^*=0.$ Thus, we have following subcases:
    \begin{enumerate}
        \item $\alpha_3=\alpha_4=\beta_3=\beta_1=0,$ then we have the   representative $\left\langle \nabla_2,\nabla_5 \right\rangle;$
        \item $\alpha_3=\alpha_4=\beta_3=0, \ \beta_1\neq0,$ then choosing $x=\frac{\beta_5}{\beta_1}, \ y=1,$ we have the representative $\left\langle \nabla_2,\nabla_1+\nabla_5 \right\rangle;$
        \item $\alpha_3=\alpha_4=0\ \beta_3\neq0,$ $\beta_1=0,$ then choosing $y=\frac{\beta_3}{\beta_5},$ we have the representative $\left\langle \nabla_2,\nabla_3+\nabla_5 \right\rangle;$
        \item $\alpha_3=\alpha_4=0,\ \beta_3\neq0,\ \beta_1\neq0,$ then choosing $x=\frac{\beta_3^2}{\beta_1\beta_5},\ y=\frac{\beta_3}{\beta_5},$ we have the representative $\left\langle \nabla_2,\nabla_1+\nabla_3+\nabla_5 \right\rangle;$
        \item $\alpha_3=0,\ \alpha_4\neq0,\ \beta_3=\beta_1=0,$ then choosing $x=1,\ y=\frac{\alpha_2}{\alpha_4},$ we have the representative $\left\langle \nabla_2+\nabla_4,\nabla_5 \right\rangle;$
        \item $\alpha_3=0,\ \alpha_4\neq0,\ \beta_3=0,\ \beta_1\neq0,$ then choosing $x=\sqrt[3]{\frac{\alpha_4^2\beta_1}{\alpha_2^2\beta_5}},\ y=\sqrt[3]{\frac{\alpha_4\beta_1^2}{\alpha_2\beta_5^2}},$
        we have the representative $\left\langle \nabla_2+\nabla_4, \nabla_1+ \nabla_5\right\rangle;$
        \item $\alpha_3=0,\ \alpha_4\neq0,\ \beta_3\neq0,$ then choosing 
        $x=\sqrt{\frac{\alpha_4\beta_3}{\alpha_2\beta_5}},\ y=\frac{\beta_3}{\beta_5},$
        we have the family of representatives $\left\langle \nabla_2+\nabla_4,\alpha\nabla_1+ \nabla_3+\nabla_5 \right\rangle;$
        \item $\alpha_3\neq0,\ \alpha_4=\beta_1=\beta_3=0,$ then choosing $x=\frac{\alpha_3}{\alpha_2},$ we have the representative $\left\langle \nabla_2+\nabla_3,\nabla_5\right\rangle;$
        \item $\alpha_3\neq0,\ \alpha_4=\beta_1=0,\ \beta_3\neq0,$ then choosing $x=\frac{\alpha_3}{\alpha_2},\ y=\frac{\beta_3}{\beta_5},$ we have the representative $\left\langle \nabla_2+\nabla_3,\nabla_3+ \nabla_5\right\rangle;$
        \item $\alpha_3\neq0,\ \alpha_4=0,\ \beta_1\neq0,$ then choosing $x=\frac{\alpha_3}{\alpha_2},\ y=\sqrt{\frac{\alpha_3\beta_1}{\alpha_2\beta_5}},$ we have the family of representatives $\left\langle\nabla_2+\nabla_3,\nabla_1+\alpha \nabla_3 +\nabla_5\right\rangle^{O(\alpha)\simeq O(-\alpha)};$
        \item $\alpha_3\neq0,\ \alpha_4\neq0,$ then choosing $x=\frac{\alpha_3}{\alpha_2},\ y=\frac{\alpha_3^2}{\alpha_2\alpha_4},$ we have the family of representatives $\left\langle \nabla_2+\nabla_3+\nabla_4,\alpha \nabla_1+\beta\nabla_4+\nabla_5 \right\rangle.$
    \end{enumerate}
    \item $\beta_5=0, \beta_4\neq0,$ then choosing $t=0,\ z=\frac{x(\alpha_4\beta_1-\alpha_1\beta_4)}{\alpha_2\beta_4},$ we can suppose $\alpha_1^*=\alpha_4^*=0$ and  have following subcases: 
    \begin{enumerate}
        \item $\alpha_3=\alpha_5=\beta_3=\beta_1=0,$ then we have the representative $\left\langle \nabla_2,\nabla_4 \right\rangle;$
        \item $\alpha_3=\alpha_5=\beta_3=0, \beta_1\neq0,$ then choosing $x=1, \ y=\sqrt{\frac{\beta_1}{\beta_4}},$ we have the representative $\left\langle \nabla_2,\nabla_1+\nabla_4 \right\rangle;$
        \item $\alpha_3=\alpha_5=0,\ \beta_3\neq0,$ then choosing $x=1, \ y=\frac{\beta_3}{\beta_4},$ we have the family of representatives $\left\langle \nabla_2,\alpha \nabla_1+ \nabla_3+\nabla_4 \right\rangle;$
        \item $\alpha_3=0,\ \alpha_5\neq0,$ then choosing $x=\frac{\alpha_5}{\alpha_2},\ y=1,$ we have the family of representatives 
        $\left\langle \nabla_2+\nabla_5,\alpha \nabla_1+\beta \nabla_3+\nabla_4 \right\rangle;$
        \item $\alpha_3\neq0,\ \alpha_5=\beta_3=\beta_1=0,$ then choosing $x=\frac{\alpha_3}{\alpha_2},$ we have the representative $\left\langle \nabla_2+\nabla_3,\nabla_4 \right\rangle;$
        \item $\alpha_3\neq0,\ \alpha_5=\beta_3=0,\ \beta_1\neq0,$ then choosing $x=\frac{\alpha_3}{\alpha_2}, \ y=\frac{\alpha_3\sqrt{\beta_1}}{\alpha_2\sqrt{\beta_4}},$ we have the representative $\left\langle \nabla_2+\nabla_3,\nabla_1+\nabla_4 \right\rangle;$
        \item $\alpha_3\neq0,\ \alpha_5=0,\ \beta_3\neq0,$ then choosing $x=\frac{\alpha_3}{\alpha_2}, \ y=\frac{\alpha_3\beta_3}{\alpha_2\beta_4},$
         we have the family of representatives $\left\langle \nabla_2+\nabla_3,\alpha \nabla_1+ \nabla_3+\nabla_4 \right\rangle;$
        \item $\alpha_3\neq0,\ \alpha_5\neq0,$ then choosing $x=\frac{\alpha_3}{\alpha_2}, \ y=\frac{\alpha_3}{\alpha_5},$ we have the family of representatives
        $\left\langle \nabla_2+\nabla_3+\nabla_5,\alpha \nabla_1+\beta \nabla_3+\nabla_4\right\rangle.$
    \end{enumerate}
    \item $\beta_5=\beta_4=0,\ \beta_3\neq0,$ then choosing $z=\frac{x(\alpha_3\beta_1-\alpha_1\beta_3)}{\alpha_2\beta_3},$ we can suppose $\alpha_1^*=\alpha_3^*=0$ and  have following subcases: 
        \begin{enumerate}
        \item $\alpha_5=\beta_1=\alpha_4=0,$ then we have the representative $\left\langle \nabla_2,\nabla_3 \right\rangle;$
        \item $\alpha_5=\beta_1=0,\ \alpha_4\neq0,$ then choosing $x=1,\ y=\frac{\alpha_2}{\alpha_4},$ we have the representative $\left\langle \nabla_2+\nabla_4, \nabla_3 \right\rangle;$
        \item $\alpha_5=0,\ \beta_1\neq0,\ \alpha_4=0$ then choosing $x=1,\ y=\frac{\beta_1}{\beta_3},$ we have the representative $\left\langle \nabla_2,\nabla_1+\nabla_3 \right\rangle;$
        \item $\alpha_5=0,\ \beta_1\neq0,\ \alpha_4\neq0,$ then choosing $x=\frac{\alpha_4\beta_1}{\alpha_2\beta_3},\ y=\frac{\alpha_4\beta_1^2}{\alpha_2\beta_3^2},$ we have the representative $ \left\langle \nabla_2+\nabla_4, \nabla_1+\nabla_3 \right\rangle;$
        \item $\alpha_5\neq0,$ then choosing $y=1,\ x=\frac{\alpha_5}{\alpha_2},\  t=-\frac{\alpha_4}{\alpha_5},$ we have the family of representatives $\left\langle \nabla_2+\nabla_5,\alpha\nabla_1+\nabla_3\right\rangle.$
    \end{enumerate}

 \item $\beta_5=\beta_4=\beta_3=0,\ \beta_1\neq0,$ then we can suppose $\alpha_1^*=0$ and consider following subcases: 
        \begin{enumerate}
        \item $\alpha_5=\alpha_4=\alpha_3=0,$ then we have the representative $\left\langle \nabla_2, \nabla_1 \right\rangle;$  \item $\alpha_5=\alpha_4=0, \ \alpha_3\neq0,$ then choosing $x=\frac{\alpha_3}{\alpha_2},$ we have the representative $\left\langle \nabla_2+\nabla_3,\nabla_1 \right\rangle;$
        \item $\alpha_5=0,\ \alpha_4\neq0, \ \alpha_3=0,$ then choosing $x=1,\ y=\frac{\alpha_2}{\alpha_4},$ we have the representative $\left\langle \nabla_2+\nabla_4,\nabla_1\right\rangle;$
        \item $\alpha_5=0,\ \alpha_4\neq0, \ \alpha_3\neq0,$ then choosing $x=\frac{\alpha_3}{\alpha_2},\ y=\frac{\alpha_3^2}{\alpha_2\alpha_4},$ we have the representative $ \left\langle \nabla_2+\nabla_3+\nabla_4,\nabla_1 \right\rangle;$
        \item $\alpha_5\neq0,\ \alpha_3=0,$ then choosing $x=\frac{\alpha_5}{\alpha_2},\ y=1,\ t=-\frac{\alpha_4}{\alpha_5},$ we have the representative $ \left\langle \nabla_2+\nabla_5,\nabla_1 \right\rangle;$
        \item $\alpha_5\neq0,\ \alpha_3\neq0,$ then choosing $x=\frac{\alpha_3}{\alpha_2},\ y=\frac{\alpha_3}{\alpha_5},\ t=-\frac{\alpha_3\alpha_4}{\alpha_5^2},$ we have the representative $ \left\langle \nabla_2+\nabla_3+\nabla_5,\nabla_1 \right\rangle.$
    \end{enumerate}

\end{enumerate}
\item $\alpha_2=0,$ then $\alpha_5\neq0$ and we have
\begin{longtable}{lcllcl}
$\alpha^*_1$&$=$&$x^2\alpha_1,$ &  $\beta^*_1$&$=$&$x^2\beta_1,$ \\
$\alpha^*_2$&$=$&$0,$ &  $\beta^*_2$&$=$&$0,$\\ 
$\alpha^*_3$&$=$&$x y\alpha_3,$ &   $\beta^*_3$&$=$&$x y \beta_3,$\\ 
$\alpha_4^*$&$=$&$y(y\alpha _4+t \alpha _5),$ &   $\beta_4^*$&$=$&$y^2\beta_4,$ \\
$\alpha_5^*$&$=$&$xy^2\alpha_5,$ &   $\beta_5^*$&$=$&$0.$ \\
\end{longtable}

\begin{enumerate}
    \item $\beta_1\neq0,$ then choosing $t=\frac{y(\alpha_1\beta_4-\alpha_4\beta_1)}{\alpha_5\beta_1},$ we can suppose $\alpha_1^*=\alpha_4^*=0$ and  have following subcases:
            \begin{enumerate}
            \item $\alpha_3=\beta_4= \beta_3=0,$ then we have the representative $ \left\langle \nabla_5,\nabla_1 \right\rangle;$
            \item $\alpha_3=\beta_4=0, \ \beta_3\neq0,$ then choosing $x=1, y = \frac{\beta_1}{\beta_3},$
            we have the representative $ \left\langle \nabla_5,\nabla_1 +\nabla_3 \right\rangle;$
            \item $\alpha_3=0, \ \beta_4\neq 0,$ then choosing $x=1, y = \sqrt{\frac{\beta_1}{\beta_4}},$
            we have the family of representatives $ \left\langle \nabla_5,\nabla_1 +\alpha \nabla_3 + \nabla_4\right\rangle^{O(\alpha)\simeq O(-\alpha)};$
        \item $\alpha_3\neq0, \ \beta_4= \beta_3=0,$ then  choosing $x=1, y = \frac{\alpha_3}{\alpha_5},$ we have the representative $ \left\langle \nabla_3+\nabla_5, \nabla_1 \right\rangle;$
         \item $\alpha_3\neq0, \ \beta_4=0, \ \beta_3\neq 0,$ then  choosing $x=\frac{\alpha_3\beta_3}{\alpha_5\beta_1}, y = \frac{\alpha_3}{\alpha_5},$ we have the representative $ \left\langle \nabla_3+\nabla_5, \nabla_1 + \nabla_3 \right\rangle;$
            \item $\alpha_3\neq0, \ \beta_4\neq 0,$ then  choosing $x=\frac{\alpha_3}{\alpha_5}\sqrt{\frac{\beta_3}{\beta_5}}, y = \frac{\alpha_3}{\alpha_5},$ we have the family of representatives $ \left\langle \nabla_3+\nabla_5, \nabla_1 + \alpha \nabla_3 + \nabla_4\right\rangle^{O(\alpha)\simeq O(-\alpha)}.$
           \end{enumerate}    
\item $\beta_1= 0, \ \beta_3\neq 0,$ then choosing $t=\frac{y(\alpha_3\beta_4-\alpha_4\beta_3)}{\alpha_5\beta_3},$ we can suppose $\alpha_3^*=\alpha_4^*=0$ and  have following subcases:
          \begin{enumerate}
            \item $\alpha_1=\beta_4= 0,$ then we have the representative $ \left\langle \nabla_5,\nabla_3 \right\rangle;$
            \item $\alpha_1=0, \ \beta_4\neq 0,$ then choosing $x=1, y =\frac{\beta_3}{\beta_4},$ we have the representative $ \left\langle \nabla_5,\nabla_3 + \nabla_4\right\rangle;$
            \item $\alpha_1\neq 0, \ \beta_4= 0,$ then choosing $x=\frac{\alpha_5}{\alpha_1}, y = 1,$ we have the representative $ \left\langle \nabla_1 +\nabla_5,\nabla_3\right\rangle;$
            \item $\alpha_1\neq 0, \ \beta_4\neq 0,$ then choosing $x=\frac{\alpha_1\beta_4^2}{\alpha_5\beta_3^2}, y = \frac{\alpha_1\beta_4}{\alpha_5\beta_3},$ we have the representative $ \left\langle \nabla_1 +\nabla_5,\nabla_3+\nabla_4\right\rangle.$
        \end{enumerate}
\item $\beta_1= \beta_3= 0,$ then $\beta_4\neq 0,$ and we can suppose $\alpha_4^*=0.$ Consider following subcases:
          \begin{enumerate}
            \item $\alpha_1=\alpha_3= 0,$ then we have the representative $ \left\langle \nabla_5,\nabla_4 \right\rangle;$
            \item $\alpha_1=0, \ \alpha_3\neq  0,$ then choosing $x=1, y =\frac{\alpha_3}{\alpha_5},$ we have the representative $ \left\langle \nabla_3 + \nabla_5,\nabla_4 \right\rangle;$
            \item $\alpha_1\neq 0, \ \alpha_3= 0,$ then choosing $x=\frac{\alpha_5}{\alpha_1}, y =1,$ we have the representative $ \left\langle \nabla_1 + \nabla_5,\nabla_4 \right\rangle;$
            \item $\alpha_1\neq 0, \ \alpha_3 \neq 0,$ then choosing $x=\frac{\alpha_3^2}{\alpha_1\alpha_5}, y =\frac{\alpha_3}{\alpha_5},$ we have the representative $ \left\langle \nabla_1 + \nabla_3 + \nabla_5,\nabla_4 \right\rangle.$
            
        \end{enumerate}
            
        \end{enumerate}
       
\end{enumerate}

Now we have the following distinct orbits:
\begin{center} $\left\langle \nabla_2,\nabla_5 \right\rangle,$ 
$\left\langle \nabla_2,\nabla_1+\nabla_5 \right\rangle,$ 
$\left\langle \nabla_2,\nabla_3+\nabla_5 \right\rangle,$ 
$\left\langle \nabla_2,\nabla_1+\nabla_3+\nabla_5 \right\rangle,$
$\left\langle \nabla_2+\nabla_4,\nabla_5 \right\rangle,$ 
$\left\langle \nabla_2+\nabla_4, \nabla_1+ \nabla_5\right\rangle,$ 
$\left\langle \nabla_2+\nabla_4,\alpha\nabla_1+ \nabla_3+\nabla_5 \right\rangle,$ 
$\left\langle \nabla_2+\nabla_3,\nabla_5\right\rangle,$ 
$\left\langle \nabla_2+\nabla_3,\nabla_3+ \nabla_5\right\rangle,$
$\left\langle\nabla_2+\nabla_3,\nabla_1+\alpha \nabla_3 +\nabla_5\right\rangle^{O(\alpha)\simeq O(-\alpha)},$  
$\left\langle \nabla_2+\nabla_3+\nabla_4,\alpha \nabla_1+\beta\nabla_4+\nabla_5 \right\rangle,$ 
$\left\langle \nabla_2,\nabla_4 \right\rangle,$  
$\left\langle \nabla_2,\nabla_1+\nabla_4 \right\rangle,$ 
$\left\langle \nabla_2,\alpha \nabla_1+ \nabla_3+\nabla_4 \right\rangle,$ 
$\left\langle \nabla_2+\nabla_5,\alpha \nabla_1+\beta \nabla_3+\nabla_4 \right\rangle,$ 
$\left\langle \nabla_2+\nabla_3,\nabla_4 \right\rangle,$ 
$\left\langle \nabla_2+\nabla_3,\nabla_1+\nabla_4 \right\rangle,$
$\left\langle \nabla_2+\nabla_3,\alpha \nabla_1+ \nabla_3+\nabla_4 \right\rangle,$ 
$\left\langle \nabla_2+\nabla_3+\nabla_5,\alpha \nabla_1+\beta \nabla_3+\nabla_4\right\rangle,$ 
$\left\langle \nabla_2,\nabla_3 \right\rangle,$ 
$\left\langle \nabla_2+\nabla_4, \nabla_3 \right\rangle,$
$\left\langle \nabla_2,\nabla_1+\nabla_3 \right\rangle,$
$ \left\langle \nabla_2+\nabla_4, \nabla_1+\nabla_3 \right\rangle,$ 
$\left\langle \nabla_2+\nabla_5,\alpha\nabla_1+\nabla_3\right\rangle,$ 
$\left\langle \nabla_2, \nabla_1 \right\rangle,$ 
$\left\langle \nabla_2+\nabla_3,\nabla_1 \right\rangle,$ 
$\left\langle \nabla_2+\nabla_4,\nabla_1\right\rangle,$ 
$ \left\langle \nabla_2+\nabla_3+\nabla_4,\nabla_1 \right\rangle,$ 
$ \left\langle \nabla_2+\nabla_5,\nabla_1 \right\rangle,$ 
$ \left\langle \nabla_2+\nabla_3+\nabla_5,\nabla_1 \right\rangle,$ 
$ \left\langle \nabla_5,\nabla_1 \right\rangle,$  
$ \left\langle \nabla_5,\nabla_1 +\nabla_3 \right\rangle,$ 
$ \left\langle \nabla_5,\nabla_1 +\alpha \nabla_3 + \nabla_4\right\rangle^{O(\alpha)\simeq O(-\alpha)},$ 
$ \left\langle \nabla_3+\nabla_5, \nabla_1 \right\rangle,$ 
$ \left\langle \nabla_3+\nabla_5, \nabla_1 + \nabla_3 \right\rangle,$ 
$ \left\langle \nabla_3+\nabla_5, \nabla_1 + \alpha \nabla_3 + \nabla_4\right\rangle^{O(\alpha)\simeq O(-\alpha)},$ 
$ \left\langle \nabla_5,\nabla_3 \right\rangle,$ 
$ \left\langle \nabla_5,\nabla_3 + \nabla_4\right\rangle,$ 
$ \left\langle \nabla_1 +\nabla_5,\nabla_3\right\rangle,$ 
$ \left\langle \nabla_1 +\nabla_5,\nabla_3+\nabla_4\right\rangle,$ 
$ \left\langle \nabla_5,\nabla_4 \right\rangle,$ 
$ \left\langle \nabla_3 + \nabla_5,\nabla_4 \right\rangle,$ 
$ \left\langle \nabla_1 + \nabla_5,\nabla_4 \right\rangle,$ 
$ \left\langle \nabla_1 + \nabla_3 + \nabla_5,\nabla_4 \right\rangle.$

\end{center}

Hence, we have the following new $5$-dimensional nilpotent bicommutative algebras (see section \ref{secteoA}):

\begin{center}
${\rm B}_{171},$
${\rm B}_{172},$
${\rm B}_{173},$
${\rm B}_{174},$
${\rm B}_{175},$
${\rm B}_{176},$
${\rm B}_{177}^{\alpha},$
${\rm B}_{178},$
${\rm B}_{179},$
${\rm B}_{180}^{\alpha},$
${\rm B}_{181}^{\alpha,\beta},$
${\rm B}_{182},$

${\rm B}_{183},$
${\rm B}_{184}^{\alpha},$
${\rm B}_{185}^{\alpha,\beta},$
${\rm B}_{186},$
${\rm B}_{187},$
${\rm B}_{188}^{\alpha},$
${\rm B}_{189}^{\alpha,\beta},$
${\rm B}_{78}^{0},$
${\rm B}_{79}^{0},$
${\rm B}_{190},$
${\rm B}_{191},$
${\rm B}_{192}^{\alpha},$
${\rm B}_{193},$
${\rm B}_{194},$
${\rm B}_{195},$
${\rm B}_{196},$
${\rm B}_{197},$
${\rm B}_{198},$
${\rm B}_{199},$
${\rm B}_{200},$
${\rm B}_{201}^{\alpha},$
${\rm B}_{202},$
${\rm B}_{203},$
${\rm B}_{204}^{\alpha},$
${\rm B}_{205},$
${\rm B}_{57}^{0,0,0},$
${\rm B}_{206},$
${\rm B}_{56}^{0},$
${\rm B}_{101}^{0},$
${\rm B}_{207},$
${\rm B}_{104}^{0},$
${\rm B}_{105}^{0}.$
\end{center}

\section{Classification theorem for $5$-dimensional bicommutative algebras}\label{secteoA}
The algebraic classification of complex $5$-dimensional bicommutative   algebras consists of two parts:
\begin{enumerate}
    \item $5$-dimensional algebras with identity $xyz=0$ (also known as $2$-step nilpotent algebras) are the intersection of all varieties of algebras defined by a family of polynomial identities of degree three or more; for example, it is in the intersection of associative, Zinbiel, Leibniz,
    Novikov, bicommutative, etc, algebras. All these algebras can be obtained as central extensions of zero-product algebras. The geometric classification of $2$-step nilpotent algebras is given in \cite{ikp20}. It is the reason why we are not interested in it.
    
     \item $5$-dimensional nilpotent (non-$2$-step nilpotent) bicommutative algebras, which are central extensions of  nilpotent bicommutative algebras with nonzero products of a smaller dimension. These algebras are classified by several steps:

    \begin{enumerate}
        \item complex split   $5$-dimensional bicommutative    algebras are classified in  \cite{kpv20};

        \item complex non-split   $5$-dimensional nilpotent commutative associative algebras are listed in \cite{krs20};
                
        \item  complex  one-generated $5$-dimensional nilpotent bicommutative algebras  are classified in \cite{kpv21};
               
        \item complex non-split non-one-generated $5$-dimensional nilpotent non-commutative bicommutative algebras are classified in Theorem A (see below).
            \end{enumerate}
  \end{enumerate}

\begin{theoremA}%\label{teorA}
Let ${\mathbb B}$ be a complex non-split  non-one-generated $5$-dimensional  nilpotent (non-$2$-step nilpotent) non-commutative bicommutative algebra.
Then ${\mathbb B}$ is isomorphic to one algebra from the following list:
{\tiny
\begin{longtable}{llllllll}
${\rm B}_{01}$ & $: $ & $e_1e_1=e_2$  & $e_1e_4=e_5$ & $e_2e_1=e_5$\\ & & $e_3e_3=e_5$  & $e_4e_1=e_5$ \\

${\rm B}_{02}$ & $: $ & $e_1e_1=e_2$ & $e_2e_1=e_5$ & $e_3e_4=e_5$ & $e_4e_3=-e_5$ \\

${\rm B}_{03}$ & $: $ & $e_1e_1=e_2$ & $e_1e_3=e_5$ & $e_2e_1=e_5$ & $e_3e_1=e_5$ \\ && $e_4e_3=e_5$ & 
$e_4e_4=e_5$ \\

${\rm B}_{04}^{\alpha}$ & $: $ & $e_1e_1=e_2$ & $e_1e_2=e_5$ & $e_2e_1=\alpha e_5$ \\ && $e_3e_3=e_5$   & $e_4e_1=e_5$ \\

${\rm B}_{05}$ & $: $ & $e_1e_1=e_2$ & $e_1e_2=e_5$ & $e_2e_1= e_5$\\ & & $e_3e_1=e_5$  & $e_3e_3=e_5$ & $e_4e_4=e_5$ \\

${\rm B}_{06}$ & $: $ & $e_1e_1=e_2$ & $e_1e_2=e_5$ & $e_2e_1= e_5$\\ & & $e_3e_1=e_5$   & $e_3e_4=e_5$ & $e_4e_3=e_5$\\

${\rm B}_{07}^{\alpha}$ & $: $ & $e_1e_1=  e_2$ & $e_1e_2= \alpha e_5$ & $e_2e_1=(\alpha+1)e_5$ \\ && $e_3e_3= e_5$   & $e_4e_4=e_5$\\

${\rm B}_{08}^{\alpha}$ & $: $ & $e_1e_1=  e_2$ & $e_1e_2=  e_5$ & $e_2e_1=\alpha e_5$ \\ && $e_3e_4= e_5$  & $e_4e_3=-e_5$\\ 

${\rm B}_{09}$ & $: $ & $e_1e_1=  e_2$ & $e_1e_2=  e_5$ & $e_2e_1= - e_5$ \\ && $e_3e_1= e_5$   & $e_3e_4= e_5$ & $e_4e_3=-e_5$\\

${\rm B}_{10}^{\alpha, \lambda}$ & $: $ & $e_1e_1=  e_2$ & $e_1e_2=  \alpha e_5$ & $e_2e_1= (\alpha+1) e_5$ \\ && $e_3e_3= \lambda e_5$   & $e_4e_3= e_5$ & $e_4e_4=e_5$\\

${\rm B}_{11}^{\lambda}$ & $: $ & $e_1e_1=  e_2$ & $e_1e_2=  \frac{-1 + \sqrt{1-4\lambda}}{2} e_5$ & $e_3e_3= \lambda e_5$  & $e_4e_3= e_5$ \\
& & $e_4e_4=e_5$ & $e_2e_1= \frac{1 + \sqrt{1-4\lambda}}{2} e_5$ & $e_3e_1=e_5$ &\\

${\rm B}_{12}^{\lambda\neq \frac{1}{4}}$ & $: $ & $e_1e_1=  e_2$ & $e_1e_2=  \frac{-1 -\sqrt{1-4\lambda}}{2} e_5$ &  $e_3e_3= \lambda e_5$  & $e_4e_3= e_5$ \\
& & $e_4e_4=e_5$ & $e_2e_1= \frac{1 - \sqrt{1-4\lambda}}{2} e_5$ & $e_3e_1=e_5$ &\\

${\rm B}_{13}$ & $: $ & $e_1e_1=e_3$  & $e_1e_3=e_5$ &  $e_2e_1=e_5$ &$e_2e_2=e_4$  \\ 
& & $e_2e_4=e_5$& $e_2e_4=e_5$ & $e_3e_1=e_5$ \\

${\rm B}_{14}^{\alpha}$ & $: $ & $e_1e_1=e_3$  & $e_1e_2=e_5$ &  $e_1e_3=\alpha e_5$ &$e_2e_1=e_5$  \\ & & $e_2e_2=e_4$  & $e_2e_4=-(1+\alpha)e_5$ & $e_3e_1=(1+\alpha)e_5$ & $e_4e_2=-\alpha e_5$\\

${\rm B}_{15}^{\alpha}$ & $: $ & $e_1e_1=e_3$  & $e_1e_3=\alpha e_5$ & $e_2e_2=e_4$  \\ & & $e_2e_4=e_5$   & $e_3e_1=(1+\alpha)e_5$ & $e_4e_2=e_5$\\

${\rm B}_{16}^{\alpha, \beta}$ & $:$ & $e_1e_1=e_3$  & $e_1e_3=\alpha e_5$ & $e_2e_2=e_4$ \\ & & $e_2e_4=\beta e_5$    & $e_3e_1=(1+\alpha)e_5$ & $e_4e_2=(1+\beta)e_5$\\

${\rm B}_{17}$ & $: $ & $e_1e_1=e_3$  & $e_1e_2=e_3$ &  $e_1e_3=e_5$ &$e_4e_4=e_5$\\

${\rm B}_{18}$ & $: $ & $e_1e_1=e_3$  & $e_1e_2=e_3$ &  $e_1e_3=e_5$ \\ && $e_2e_4=e_5$   & $e_4e_4=e_5$\\

${\rm B}_{19}$ & $: $ & $e_1e_1=e_3$  & $e_1e_2=e_3$ &  $e_1e_3=e_5$\\ & &$e_2e_1=e_5$   & $e_4e_4=e_5$\\

${\rm B}_{20}$ & $: $ & $e_1e_1=e_3$  & $e_1e_2=e_3$ &  $e_1e_3=e_5$ \\ &&$e_2e_1=e_5$  & $e_2e_4=e_5$ &$e_4e_4=e_5$\\

${\rm B}_{21}^{\alpha}$ & $: $ & $e_1e_1=e_3$  & $e_1e_2=e_3$ &  $e_1e_3=e_5$ &$e_2e_1=e_5$ \\ &&$e_2e_2=e_5$  &$e_2e_4=\alpha e_5$ &$e_4e_4=e_5$\\

${\rm B}_{22}^{\alpha}$ & $: $ & $e_1e_1=e_3$  & $e_1e_2=e_3$ &  $e_1e_3=e_5$ \\ && $e_2e_2=e_5$   &$e_2e_4=\alpha e_5$ &$e_4e_4=e_5$\\

${\rm B}_{23}$ & $: $ & $e_1e_1=e_3$  & $e_1e_2=e_3$ &  $e_1e_3=e_5$ \\ && $e_2e_2=e_5$  & $e_4e_1=e_5$ &$e_4e_2=e_5$\\

${\rm B}_{24}$ & $: $ & $e_1e_1=e_3$  & $e_1e_2=e_3$ &  $e_1e_3=e_5$ & $e_2e_2=e_5$ \\ && $e_2e_4=-e_5$ & $e_4e_1=e_5$ &$e_4e_2=e_5$\\

${\rm B}_{25}^{\alpha}$ & $: $ & $e_1e_1=e_3$  & $e_1e_2=e_3$ &  $e_1e_3=e_5$ \\ & & $e_2e_4=\alpha e_5$ & $e_4e_1=e_5$ &$e_4e_2=e_5$\\

${\rm B}_{26}$ & $: $ & $e_1e_1=e_3$  & $e_1e_2=e_3$ &  $e_1e_3=e_5$ \\ && $e_2e_1=e_5$ & $e_4e_2=e_5$\\

${\rm B}_{27}$ & $: $ & $e_1e_1=e_3$  & $e_1e_2=e_3$ &  $e_1e_3=e_5$ \\ && $e_2e_2=e_5$  & $e_2e_4=-e_5$ & $e_4e_2=e_5$\\

${\rm B}_{28}^{\alpha}$ & $: $ & $e_1e_1=e_3$  & $e_1e_2=e_3$ &  $e_1e_3=e_5$ \\ && $e_2e_4=\alpha e_5$  & $e_4e_2=e_5$\\

${\rm B}_{29}$ & $: $ & $e_1e_1=e_3$  & $e_1e_2=e_3$ &  $e_1e_3=e_5$ &  $e_4e_1=e_5$\\

${\rm B}_{30}$ & $: $ & $e_1e_1=e_3$  & $e_1e_2=e_3$ &  $e_1e_3=e_5$\\ & & $e_2e_2=e_5$  & $e_4e_1=e_5$\\

${\rm B}_{31}$ & $: $ & $e_1e_1=e_3$  & $e_1e_2=e_3$ &  $e_1e_3=e_5$ \\ & & $e_2e_4=e_5$  & $e_4e_1=e_5$\\

${\rm B}_{32}$ & $: $ & $e_1e_1=e_3$  & $e_1e_2=e_3$ &  $e_1e_3=e_5$ & $e_2e_4=e_5$ \\

${\rm B}_{33}$ & $: $ & $e_1e_1=e_3$  & $e_1e_2=e_3$ &  $e_3e_1=e_5$ \\ & & $e_3e_2=e_5$  & $e_4e_4=e_5$ \\

${\rm B}_{34}$ & $: $ & $e_1e_1=e_3$  & $e_1e_2=e_3$ &  $e_1e_4=e_5$ \\ && $e_3e_1=e_5$  & $e_3e_2=e_5$ & $e_4e_4=e_5$ \\

${\rm B}_{35}^{\alpha}$ & $: $ & $e_1e_1=e_3$  & $e_1e_2=e_3+e_5$ & $e_1e_4=\alpha e_5$\\ & & $e_3e_1=e_5$  & $e_3e_2=e_5$ & $e_4e_4=e_5$ \\

${\rm B}_{36}$ & $: $ & $e_1e_1=e_3$  & $e_1e_2=e_3+e_5$ & $e_2e_2=e_5$  \\ && $e_3e_1=e_5$ & $e_3e_2=e_5$ & $e_4e_4=e_5$ \\

${\rm B}_{37}$ & $: $ & $e_1e_1=e_3$  & $e_1e_2=e_3+e_5$ & $e_1e_4=e_5$ & $e_2e_2=e_5$  \\ & & $e_3e_1=e_5$ & $e_3e_2=e_5$ & $e_4e_4=e_5$ \\

${\rm B}_{38}^{\alpha}$ & $: $ & $e_1e_1=e_3$  & $e_1e_2=e_3$ & $e_1e_4=\alpha e_5$ & $e_2e_2=e_5$  \\ && $e_3e_1=e_5$ & $e_3e_2=e_5$ & $e_4e_4=e_5$ \\

${\rm B}_{39}$ & $: $ & $e_1e_1=e_3$  & $e_1e_2=e_3$ & $e_1e_3=e_5$ \\ &&  $e_3e_1=e_5$  & $e_3e_2=e_5$ & $e_4e_4=e_5$ \\

${\rm B}_{40}$ & $: $ & $e_1e_1=e_3$  & $e_1e_2=e_3$ & $e_1e_3=e_5$ & $e_1e_4=e_5$ \\ &&  $e_3e_1=e_5$ & $e_3e_2=e_5$ & $e_4e_4=e_5$ \\

${\rm B}_{41}^{\alpha}$ & $: $ & $e_1e_1=e_3$  & $e_1e_2=e_3$ & $e_1e_3=e_5$ & $e_1e_4=\alpha e_5$ \\ && $e_2e_2=e_5$ & $e_3e_1=e_5$  & $e_3e_2=e_5$ & $e_4e_4=e_5$ \\

${\rm B}_{42}$ & $: $ & $e_1e_1=e_3$  & $e_1e_2=e_3$ & $e_1e_4=e_5$ \\ && $e_2e_4=e_5$  & $e_3e_1=e_5$ & $e_3e_2=e_5$ \\

${\rm B}_{43}$ & $: $ & $e_1e_1=e_3$  & $e_1e_2=e_3+e_5$ & $e_1e_4=e_5$ \\ && $e_2e_4=e_5$    & $e_3e_1=e_5$ & $e_3e_2=e_5$ \\

${\rm B}_{44}$ & $: $ & $e_1e_1=e_3$  & $e_1e_2=e_3$ & $e_1e_4=e_5$ & $e_2e_4=e_5$ \\ && $e_3e_1=e_5$ & $e_3e_2=e_5$ & $e_4e_1=e_5$ \\

${\rm B}_{45}$ & $: $ & $e_1e_1=e_3$  & $e_1e_2=e_3+e_5$ & $e_2e_4=e_5$ \\ && $e_3e_1=e_5$    & $e_3e_2=e_5$ \\

${\rm B}_{46}$ & $: $ & $e_1e_1=e_3$  & $e_1e_2=e_3+e_5$ & $e_2e_4=e_5$ \\ && $e_3e_1=e_5$  & $e_3e_2=e_5$ & $e_4e_1=e_5$  \\

${\rm B}_{47}^{\alpha}$ & $: $ & $e_1e_1=e_3$  & $e_1e_2=e_3$ & $e_1e_3=e_5$ & $e_1e_4=e_5$ \\ && $e_2e_4=e_5$ & $e_3e_1=e_5$ & $e_3e_2=e_5$ & $e_4e_1=\alpha e_5$  \\

${\rm B}_{48}^{\alpha, \beta}$ & $: $ & $e_1e_1=e_3$  & $e_1e_2=e_3$ & $e_1e_3=\beta e_5$ & $e_2e_4=e_5$ \\ && $e_3e_1=e_5$ & $e_3e_2=e_5$ & $e_4e_1=\alpha e_5$  \\

${\rm B}_{49}$ & $: $ & $e_1e_1=e_3$  & $e_1e_2=e_3+e_5$ &  $e_1e_4=-e_5$ \\ &&  $e_3e_1=e_5$   &   $e_3e_2=e_5$ & $e_4e_1=e_5$ \\

${\rm B}_{50}^{\alpha}$ & $: $ & $e_1e_1=e_3$  & $e_1e_2=e_3$ &  $e_1e_4=\alpha e_5$ \\ &&  $e_3e_1=e_5$  & $e_3e_2=e_5$ & $e_4e_1=e_5$ \\

${\rm B}_{51}^{\alpha}$ & $: $ & $e_1e_1=e_3$  & $e_1e_2=e_3$ & $e_1e_3=e_5$ & $e_1e_4=\alpha e_5$ \\ && $e_3e_1=e_5$ & $e_3e_2=e_5$ & $e_4e_1=e_5$ \\

${\rm B}_{52}$ & $: $ & $e_1e_1=e_3$  & $e_1e_2=e_3$ & $e_1e_4=e_5$ \\ &&  $e_3e_1=e_5$  & $e_3e_2=e_5$  \\

${\rm B}_{53}$ & $: $ & $e_1e_1=e_3$  & $e_1e_2=e_3$ & $e_1e_4=e_5$ \\ && $e_2e_2=e_5$  &  $e_3e_1=e_5$ & $e_3e_2=e_5$  \\

${\rm B}_{54}$ & $: $ & $e_1e_1=e_3$  & $e_1e_2=e_3$ & $e_1e_3=e_5$ \\ & & $e_1e_4=e_5$  &  $e_3e_1=e_5$ & $e_3e_2=e_5$  \\

${\rm B}_{55}$ & $: $ & $e_1e_1=e_3$  & $e_1e_2=e_3$ & $e_1e_3=e_5$  & $e_1e_4=e_5$ \\ &&  $e_2e_2=e_5$ &  $e_3e_1=e_5$ & $e_3e_2=e_5$  \\

${\rm B}_{56}^{\alpha}$ & $: $ &$e_1e_1 = e_5$ & $e_1e_2 = e_3$ & $e_2e_1 = e_4$ \\ &&  $e_2e_2 = -e_3$  &  $e_2e_4 = \alpha e_5$ &  $e_3e_2 = e_5$\\

${\rm B}_{57}^{\alpha, \beta, \gamma}$ & $: $ & $e_1e_2 = e_3$ & $e_1e_3 = \gamma e_5$ & $e_2e_1 = e_4$ &  $e_2e_2 = -e_3$ \\ && $e_2e_3 = - \gamma e_5$  &  $e_2e_4 = \alpha e_5$ &  $e_3e_2 = e_5$ & $e_4e_1 = \beta e_5$\\

${\rm B}_{58}^{(\alpha, \beta)\neq (0,0)}$ & $: $ & $e_1e_2 = e_3$ & $e_1e_3 = e_5$ & $e_2e_1 = e_4$ &  $e_2e_2 = -e_3$ \\ && $e_2e_3 = - e_5$  &  $e_2e_4 = \alpha e_5$ &  $e_4e_1 = \beta e_5$\\

${\rm B}_{59}^{\alpha \neq 0}$ & $: $ & $e_1e_2 = e_3$ & $e_1e_3 = e_5$ & $e_2e_1 = e_4$ \\ &&  $e_2e_2 = -e_3 + e_5$   & $e_2e_3 = - e_5$ & $e_4e_1 = \alpha e_5$\\

${\rm B}_{60}^{\alpha \neq 1, \beta\neq0}$ & $: $ &$e_1e_1 =e_3$ & $e_1e_2 =e_4$ & $e_1e_3 = e_5$ & $e_1e_4=\beta e_5$ \\ & & $e_2e_1 =-\alpha e_3$ &  $e_2e_2 = -e_4$  &   $e_2e_3 =-\alpha e_5$ & $e_2e_4 =-\beta e_5$ &\\

${\rm B}_{61}^{\alpha \neq 1, \beta\neq0}$ & $: $ &$e_1e_1 =e_3$ & $e_1e_2 =e_4+e_5$  & $e_1e_4=\beta e_5$ & $e_2e_1 =-\alpha e_3$ \\ &&  $e_2e_2 = -e_4$   &    $e_2e_4 =-\beta e_5$  & $e_3e_1 = e_5$ & \\

${\rm B}_{62}^{\alpha \neq 1, (\beta,\gamma)\neq(0,0),\delta}$ & $: $ &$e_1e_1 =e_3$ & $e_1e_2=e_4$  & $e_1e_3 =\delta e_5$ &  $e_1e_4=\gamma e_5$ \\ & & $e_2e_1 =-\alpha e_3$   &  $e_2e_2 = -e_4$  &  $e_2e_3 =-\alpha \delta e_5$  & $e_2e_4 =-\gamma e_5$ \\ && $e_3e_1 = e_5$ & $e_4e_2 =\beta e_5$ \\

${\rm B}_{63}^{\alpha\neq0}$ & $: $ &$e_1e_1 =e_3$ & $e_1e_2=e_4$  & $e_1e_3 =\alpha e_5$ \\ & & $e_2e_1 =e_5$  &  $e_2e_2 = -e_4$   &  $e_4e_2 = e_5$ \\

${\rm B}_{64}^{\alpha \neq 0,\beta}$ & $: $ &$e_1e_1 =e_3$ & $e_1e_2=e_4$  & $e_1e_3 = \alpha e_5$ &  $e_1e_4=\beta e_5$ \\ &&   $e_2e_2 = -e_4$    &  $e_2e_4 =-\beta e_5$ & $e_4e_2 = e_5$ \\

${\rm B}_{65}$ & $: $ &$e_1e_1 =e_3$ & $e_1e_2=e_4$  & $e_1e_3 =e_5$ &  $e_1e_4=-e_5$ \\ && $e_2e_1 =-e_3$   &  $e_2e_2 = -e_4$  &  $e_2e_3 =-e_5$  & $e_2e_4 =e_5$ \\

${\rm B}_{66}$ & $: $ &$e_1e_1 =e_3$ & $e_1e_2=e_4+e_5$  & $e_1e_3 =e_5$ &  $e_1e_4=-e_5$ \\ && $e_2e_1 =-e_3$   &  $e_2e_2 = -e_4$  &  $e_2e_3 =-e_5$  & $e_2e_4 =e_5$ \\

${\rm B}_{67}$ & $: $ &$e_1e_1 =e_3$ & $e_1e_2=e_4+e_5$  & $e_1e_3 =e_5$ &  $e_1e_4=-e_5$ \\ && $e_2e_1 =-e_3-e_5$   &  $e_2e_2 = -e_4$  &  $e_2e_3 =-e_5$  & $e_2e_4 =e_5$ \\

${\rm B}_{68}^{\alpha}$ & $: $ &$e_1e_1 =e_3$ & $e_1e_2=e_4$  & $e_1e_3 =\alpha e_5$ &  $e_1e_4=-\alpha e_5$ \\ && $e_2e_1 =-e_3$   &  $e_2e_2 = -e_4$  &  $e_2e_3 =-\alpha e_5$  & $e_2e_4 =\alpha e_5$ \\ && $e_3e_1 =e_5$ & $e_3e_2 =-e_5$   & $e_4e_1 =-e_5$ & $e_4e_2 =e_5$\\

${\rm B}_{69}^{\alpha}$ & $: $ &$e_1e_1 =e_3$ & $e_1e_2=e_4$  & $e_1e_3 =\alpha e_5$ &  $e_1e_4=-\alpha e_5$ \\ & & $e_2e_1 =-e_3+e_5$   &  $e_2e_2 = -e_4$  &  $e_2e_3 =-\alpha e_5$  & $e_2e_4 =\alpha e_5$ \\ & & $e_3e_1 =e_5$ & $e_3e_2 = -e_5$   & $e_4e_1 =-e_5$ & $e_4e_2 =e_5$\\

${\rm B}_{70}$ & $: $ &$e_1e_1 =e_3$ & $e_1e_2=e_4$  & $e_1e_3 =e_5$ &  $e_2e_1 =-e_3$ \\ &&  $e_2e_2 = -e_4$   &  $e_2e_3 =- e_5$  & $e_3e_1 =e_5$ & $e_3e_2 = -e_5$ \\ && $e_4e_1 =-e_5$ & $e_4e_2 =e_5$\\

${\rm B}_{71}$ & $: $ &$e_1e_1 =e_3$ & $e_1e_2=e_4$  & $e_1e_3 =e_5$ &  $e_2e_1 =-e_3+e_5$  \\ &&  $e_2e_2 = -e_4$    &  $e_2e_3 =- e_5$  & $e_3e_1 =e_5$ & $e_3e_2 = -e_5$ \\ && $e_4e_1 =-e_5$ & $e_4e_2 =e_5$\\

${\rm B}_{72}^{\alpha}$ & $: $ &$e_1e_1 =e_3$ & $e_1e_2=e_4$  & $e_1e_3 =\alpha e_5$ \\ &&  $e_1e_4=e_5$   & $e_2e_1 =-e_3$   &  $e_2e_2 = -e_4$  \\ &&  $e_2e_3 =-\alpha e_5$  & $e_2e_4 =-e_5$ & $e_3e_1 =e_5$ \\

${\rm B}_{73}^{\alpha}$ & $: $ &$e_1e_1 =e_3$ & $e_1e_2=e_4+e_5$  & $e_1e_3 =\alpha e_5$ \\ & &  $e_1e_4=e_5$   & $e_2e_1 =-e_3$  &  $e_2e_2 = -e_4$ \\ & &  $e_2e_3 =-\alpha e_5$  & $e_2e_4 =-e_5$ & $e_3e_1 =e_5$ \\

${\rm B}_{74}^{\alpha,\beta}$ & $: $ &$e_1e_1 =e_3$ & $e_1e_2=e_4$  & $e_1e_3 =\alpha e_5$ &  $e_1e_4=\beta e_5$ \\ && $e_2e_1 =-e_3$  &  $e_2e_2 = -e_4$  &  $e_2e_3 =-\alpha e_5$  & $e_2e_4 =-\beta e_5$ \\ && $e_3e_1 =e_5$ & $e_3e_2 =e_5$  & $e_4e_1 =e_5$\\

${\rm B}_{75}^{\alpha}$ & $: $ &
$e_1e_1 =e_3$ & $e_1e_2=e_4$  & $e_1e_3 =\alpha e_5$ \\
&&  $e_2e_1 =-e_3$ &  $e_2e_2 = -e_4$  &  $e_2e_3 =-\beta e_5$  \\& & $e_3e_1 =-2e_5$ & $e_3e_2 =e_5$ & $e_4e_1 =e_5$\\ 

${\rm B}_{76}^{\alpha}$ & $: $ &$e_1e_1 =e_3$ & $e_1e_2=e_4$  & $e_1e_3 =\alpha e_5$ &  $e_1e_4= e_5$ \\ && $e_2e_1 =-e_3$   &  $e_2e_2 = -e_4$  &  $e_2e_3 =-\alpha e_5$  & $e_2e_4 =-e_5$ \\ && $e_3e_1 =-2e_5$ & $e_3e_2 =e_5$ & $e_4e_1 =e_5$\\

${\rm B}_{77}$ & $: $ &$e_1e_2 =e_3$ & $e_2e_1=e_4$  & $e_3e_2 =e_5$ &  $e_4e_1= e_5$ \\

${\rm B}_{78}^{\alpha}$ & $: $ & $e_1e_2 =e_3$ & $e_2e_1=e_4$  & $e_2e_4 =e_5$ &  $e_3e_2=\alpha e_5$ \\

${\rm B}_{79}^{\alpha}$ & $: $ & $e_1e_1=e_5$ & $e_1e_2 =e_3$ & $e_2e_1=e_4$  \\ && $e_2e_4 =e_5$   &  $e_3e_2=\alpha e_5$ \\

${\rm B}_{80}$ & $: $ & $e_1e_2 =e_3$ & $e_1e_3 =e_5$ & $e_2e_1=e_4$  \\ && $e_2e_4 =e_5$  &  $e_3e_2=\alpha e_5$ \\

${\rm B}_{81}^{(\alpha,\beta)\neq (0,0)}$ & $: $ & $e_1e_2 =e_3$ & $e_1e_3 =\beta e_5$ & $e_2e_1=e_4$ \\ & & $e_2e_4 =e_5$  &  $e_3e_2=\alpha e_5$  &  $e_4e_1=e_5$ \\

${\rm B}_{82}$ & $: $ & $e_1e_1 =e_4$ & $e_1e_2 = e_3$ & $e_1e_3=e_5$  & $e_1e_4 =-2e_5$ \\ & &  $e_2e_1=-e_3$   &  $e_2e_2=2e_3+e_4$ & $e_2e_4 =-e_5$ \\

${\rm B}_{83}^{\alpha}$ & $: $ & $e_1e_1 =e_4$ & $e_1e_2 = e_3$ & $e_1e_3=\alpha e_5$  & $e_1e_4 =(1-2\alpha)e_5$ \\ &&  $e_2e_1=-e_3$  &  $e_2e_2=2e_3+e_4$ & $e_2e_3 = e_5$ & $e_2e_4 =-\alpha e_5$ \\

${\rm B}_{84}$ & $: $ & $e_1e_1 =e_4$ & $e_1e_2 = e_3$ & $e_1e_3= e_5$  & $e_1e_4 =-e_5$ \\ & &  $e_2e_1=-e_3+e_5$   &  $e_2e_2=2e_3+e_4$ & $e_2e_3 = e_5$ & $e_2e_4 =-e_5$ \\

${\rm B}_{85}$ & $: $ & $e_1e_1 =e_4$ & $e_1e_2 = e_3$ &  $e_1e_4 =-\frac{1}{2}e_5$ &  $e_2e_1=-e_3$  \\ & &  $e_2e_2=2e_3+e_4+e_5$  & $e_2e_3 =-\frac12 e_5$ & $e_3e_1 =e_5$ & $e_3e_2 =(i-\frac32)e_5$ \\ && $e_4e_1 =(-\frac12 - i)e_5$ & $e_4e_2 =e_5$ \\

${\rm B}_{86}$ & $: $ & $e_1e_1 =e_4$ & $e_1e_2 = e_3$ &  $e_1e_4 =-\frac{1}{2}e_5$ &  $e_2e_1=-e_3$  \\ & &  $e_2e_2=2e_3+e_4+e_5$  & $e_2e_3 =-\frac12 e_5$ & $e_3e_1 =e_5$ & $e_3e_2 =(-\frac32 - i)e_5$ \\ & & $e_4e_1 =(i-\frac12)e_5$ & $e_4e_2 =e_5$ \\

${\rm B}_{87}^{\alpha,\beta,\gamma}$ & $: $ & $e_1e_1 =e_4$ & $e_1e_2 = e_3$ &  $e_1e_3 =\alpha e_5$ & $e_1e_4=(\gamma-2\alpha)e_5$ \\ &&  $e_2e_1=-e_3$    &  $e_2e_2=2e_3+e_4$& $e_2e_3 =\gamma e_5$ & $e_2e_4 =-\alpha e_5$ \\ && $e_3e_1=e_5$ & $e_3e_2 =(\beta-2)e_5$   & $e_4e_1 =-\beta e_5$ & $e_4e_2 = e_5$ \\

${\rm B}_{88}^{\alpha}$ & $: $ &$e_1e_2 =e_4$ & $e_2e_1 =-\frac12 e_4$ & $e_2e_2 = e_3$ & $e_2e_3=\alpha e_5$ \\ & & $e_3e_2 = e_5$  &  $e_3e_1 =-\frac12 e_5$  &   $e_4e_2 =e_5$ \\

${\rm B}_{89}$ & $: $ &$e_1e_2 =e_4$ & $e_1e_3 =e_5$ & $e_2e_1 =e_4+e_5$ & $e_2e_2 = e_3$  \\ & & $e_2e_4=e_5$  &  $e_3e_1 =e_5$  &   $e_4e_2 =e_5$ \\

${\rm B}_{90}$ & $: $ &$e_1e_2 =e_4$ & $e_1e_3 =e_5$ & $e_2e_1 =e_4+e_5$ & $e_2e_2 = e_3$ \\ & & $e_2e_4=e_5$ &  $e_3e_1 =e_5$  & $e_3e_2 =e_5$  &   $e_4e_2 =e_5$ \\

${\rm B}_{91}^{\alpha, \beta}$ & $: $ &$e_1e_2 =e_4$ & $e_1e_3 =\alpha e_5$ & $e_1e_4 =e_5$ & $e_2e_2 = e_3$ \\ && $e_2e_3=e_5$  &  $e_2e_4 =\alpha e_5$  & $e_4e_2 =\beta e_5$\\

${\rm B}_{92}^{\alpha}$ & $: $ &$e_1e_2 =e_4$ &  $e_1e_4 =e_5$ & $e_2e_1 = e_5$ \\ && $e_2e_2 = e_3$ & $e_3e_2=e_5$  &  $e_4e_2 =\alpha e_5$ \\

${\rm B}_{93}^{\alpha}$ & $: $ &$e_1e_2 =e_4$ &  $e_1e_4 =e_5$ &  $e_2e_2 = e_3$ \\ && $e_3e_2=e_5$   &  $e_4e_2 =\alpha e_5$ \\

${\rm B}_{94}$ & $: $ &$e_1e_1=e_5$ & $e_1e_2 =e_4$ &  $e_2e_2 =e_3$ \\ &&  $e_2e_3 = e_5$   &  $e_4e_2 = e_5$ \\

${\rm B}_{95}$ & $: $ &$e_1e_1=e_5$ & $e_1e_2 =e_4$ &  $e_1e_3 =2e_5$ \\ &&  
$e_2e_1 = -2e_4$   &  $e_2e_2 = e_3$   &  $e_2e_3 = e_5$ \\ &&
$e_2e_4 = 2e_5$ &  $e_3e_1 =-2e_5$  &  $e_4e_2 = e_5$ \\

${\rm B}_{96}^{\alpha}$ & $: $ & $e_1e_2=e_4$ & $e_1e_3 =e_5$ &  $e_2e_1 =-2e_4$ &  $e_2e_2 = e_3$  \\ &&  $e_2e_3 = e_5$   &  $e_2e_4 = e_5$ &  $e_3e_1 = -2\alpha e_5$ &  $e_4e_2 = \alpha e_5$ \\

${\rm B}_{97}^{\alpha}$ & $: $ & $e_1e_2=e_4$ & $e_1e_3 =e_5$ &  $e_2e_1 =-2e_4$ &  $e_2e_2 = e_3$ \\ &&  $e_2e_3 = \alpha e_5$    &  $e_2e_4 = e_5$ &  $e_3e_2=e_5$ \\

${\rm B}_{98}^{\alpha}$ & $: $ & $e_1e_2=e_4$ & $e_2e_1 =\alpha e_4$ &  $e_2e_2 =e_3$ \\ & &  $e_3e_1 = \alpha e_5$    &  $e_4e_2 =  e_5$ \\

${\rm B}_{99}^{\alpha,\beta}$ & $: $ & $e_1e_2=e_4$ & $e_1e_3=e_5$ & $e_2e_1 =\alpha e_4$ &  $e_2e_2 =e_3$ \\ & &  $e_2e_4 = e_5$  &  $e_3e_1 = \alpha\beta e_5$  &  $e_4e_2 = \beta e_5$ \\

${\rm B}_{100}^{\alpha\neq-2, \beta}$ & $: $ & $e_1e_2=e_4$ & $e_1e_3=e_5$ & $e_2e_1 =\alpha e_4$ &  $e_2e_2 =e_3$ \\ &&  $e_2e_4 = e_5$  &  $e_3e_1 = \alpha\beta e_5$  & $e_3e_2 = e_5$  &  $e_4e_2 = \beta e_5$ \\

${\rm B}_{101}^{\alpha}$ & $: $ & $e_1e_2=e_4$ &  $e_2e_1 =\alpha e_4$ &  $e_2e_2 =e_3$\\ & &   $e_3e_1 = \alpha e_5$    &  $e_4e_2 = e_5$ \\

${\rm B}_{102}^{\alpha\neq0}$ & $: $ & $e_1e_2=e_4$ & $e_1e_3=e_5$ & $e_2e_1 =\alpha e_4$ &  $e_2e_2 =e_3$ \\ & &  $e_2e_4 = e_5$  &  $e_3e_1 =-e_5$  &  $e_4e_2 =-\frac{1}{\alpha} e_5$ \\

${\rm B}_{103}^{\alpha\neq0}$ & $: $ & $e_1e_2=e_4$ & $e_1e_3=e_5$ & $e_2e_1 =\alpha e_4$ &  $e_2e_2 =e_3$ \\ &&  $e_2e_4 = e_5$  &  $e_3e_1 =-e_5$  & $e_3e_2 =e_5$  &  $e_4e_2 =-\frac{1}{\alpha} e_5$ \\

${\rm B}_{104}^{\alpha}$ & $: $ & $e_1e_1=e_5$ & $e_1e_2=e_4$ & $e_1e_3=\alpha e_5$ & $e_2e_1 =\alpha e_4$ \\ &&  $e_2e_2 =e_3$   &  $e_2e_4 =\alpha e_5$&  $e_3e_1 =-\alpha e_5$ &  $e_4e_2 =- e_5$ \\

${\rm B}_{105}^{\alpha\neq-2}$ & $: $ & $e_1e_1=e_5$ & $e_1e_2=e_4$ & $e_1e_3=\alpha e_5$\\ & & $e_2e_1 =\alpha e_4$   &  $e_2e_2 =e_3$  &  $e_2e_4 =\alpha e_5$ \\ &&  $e_3e_1 =-\alpha e_5$  & $e_3e_2 =e_5$ &  $e_4e_2 =- e_5$\\

${\rm B}_{106}$ & $: $ & $e_1e_1=e_2$ & $e_2e_1=e_3$ & $e_3e_1=e_5$ & $e_4e_4 = e_5$ \\

${\rm B}_{107}$ & $: $ & $e_1e_1=e_2$ &  $e_1e_2=e_5$ & $e_2e_1=e_3$ \\ & & $e_3e_1=e_5$  & $e_4e_4 = e_5$ \\

${\rm B}_{108}$ & $: $ & $e_1e_1=e_2$ &  $e_1e_4=e_5$ & $e_2e_1=e_3$ & $e_3e_1=e_5$ \\

${\rm B}_{109}$ & $: $ & $e_1e_1=e_2$ &  $e_1e_2=e_5$ &  $e_1e_4=e_5$ \\ & & $e_2e_1=e_3$ & $e_3e_1=e_5$ \\

${\rm B}_{110}$ & $: $ & $e_1e_1=e_2$ &  $e_1e_2=e_3$ &  $e_1e_3=e_5$ \\ && $e_2e_1=e_5$   & $e_4e_1=e_5$ \\

${\rm B}_{111}$ & $: $ & $e_1e_1=e_2$ &  $e_1e_2=e_3$ &  $e_1e_3=e_5$ \\ && $e_2e_1=e_5$  & $e_4e_4=e_5$ \\

${\rm B}_{112}^{\alpha}$ & $: $ & $e_1e_1=e_2$ &  $e_1e_2=e_3$ &  $e_1e_3=e_5$ & $e_2e_1=e_3+\alpha e_5$ \\ && $e_2e_2=e_5$   & $e_3e_1=e_5$ & $e_4e_1=e_5$ & $e_4e_4=e_5$ \\

${\rm B}_{113}$ & $: $ & $e_1e_1=e_2$ &  $e_1e_2=e_3$ &  $e_1e_3=e_5$ & $e_2e_1=e_3+e_5$ \\ && $e_2e_2=e_5$  & $e_3e_1=e_5$ & $e_4e_1=e_5$  \\

${\rm B}_{114}$ & $: $ & $e_1e_1=e_2$ &  $e_1e_2=e_3$ &  $e_1e_3=e_5$ & $e_2e_1=e_3+e_5$ \\ && $e_2e_2=e_5$   & $e_3e_1=e_5$ & $e_4e_4=e_5$ \\

${\rm B}_{115}^{\alpha\neq 1}$ & $: $ & $e_1e_1=e_2$ &  $e_1e_2=e_3$ &  $e_1e_3=e_5$ & $e_2e_1=\alpha e_3$ \\ && $e_2e_2=\alpha e_5$  & $e_3e_1=\alpha e_5$ & $e_4e_1=e_5$ \\

${\rm B}_{116}^{\alpha\neq 1}$ & $: $ & $e_1e_1=e_2$ &  $e_1e_2=e_3$ &  $e_1e_3=e_5$ & $e_2e_1=\alpha e_3$ \\ && $e_2e_2=\alpha e_5$  & $e_3e_1=\alpha e_5$ & $e_4e_4=e_5$ \\

${\rm B}_{117}^{\alpha\neq0,   \beta}$ & $: $ & $e_1e_1=e_2$ &  $e_1e_2=e_4$ &  $e_1e_3=e_4$ \\ && $e_1e_4=e_5$  & $e_2e_1=\alpha e_4$  & $e_2e_2=\alpha e_5$ \\ & & $e_2e_3=\alpha e_5$  & $e_3e_3=\beta e_5$  & $e_4e_1=\alpha e_5$ \\ 

${\rm B}_{118}^{\alpha\neq0,1}$ & $: $ & $e_1e_1=e_2$ &  $e_1e_2=e_4$ &  $e_1e_3=e_4$ & $e_1e_4=e_5$ \\ && $e_2e_1=\alpha e_4$   & $e_2e_2=\alpha e_5$ & $e_2e_3=\alpha e_5$ & $e_3e_1=e_5$ \\ && $e_3e_3=\frac{2\alpha^2}{(\alpha-1)^2} e_5$ & $e_4e_1=\alpha e_5$ \\ 

${\rm B}_{119}^{\alpha}$ & $: $ & $e_1e_1=e_2$ &  $e_1e_2=e_5$ &  $e_1e_3=e_4$ & $e_2e_1=e_4$ \\ && $e_2e_3=e_5$  & $e_3e_3=\alpha e_5$ & $e_4e_1=e_5$ \\ 

${\rm B}_{120}^{\alpha}$ & $: $ & $e_1e_1=e_2$ &  $e_1e_3=e_4$ & $e_2e_1=e_4$ \\ &&
$e_2e_3=e_5$   & $e_3e_3=\alpha e_5$   & $e_4e_1=e_5$ \\ 

${\rm B}_{121}^{\alpha}$ & $: $ & $e_1e_1=e_2$ &  $e_1e_2=\alpha e_5$ &  $e_1e_3=e_4$ & $e_2e_1=e_4+e_5$ \\ && $e_2e_3=e_5$   & $e_3e_3=2e_5$ & $e_4e_1=e_5$ \\

${\rm B}_{122}^{\alpha}$ & $: $ & $e_1e_1=e_2$  & $e_1e_2=e_4$ & $e_1e_4=e_5$ & $e_2e_1=\alpha e_5$ \\ && $e_3e_1=e_4+e_5$  & $e_3e_2=e_5$ & $e_3e_3=2e_5$\\

${\rm B}_{123}^{\alpha}$ & $: $ & $e_1e_1=e_2$  & $e_1e_2=e_4$ & $e_1e_4=e_5$ \\ && $e_3e_1=e_4$  & $e_3e_2=e_5$  & $e_3e_3=\alpha e_5$\\

${\rm B}_{124}^{\alpha}$ & $: $ & $e_1e_1=e_2$  & $e_1e_2=e_4$ & $e_1e_4=e_5$ & $e_2e_1=e_5$ \\ && $e_3e_1=e_4$  & $e_3e_2=e_5$ & $e_3e_3=\alpha e_5$ \\

${\rm B}_{125}$ & $: $ & $e_1e_2=e_3$  & $e_1e_3=e_4$ & $e_1e_4=e_5$ \\ && 
$e_3e_2=e_4$  & $e_3e_3=e_5$ & $e_4e_2=e_5$ \\

${\rm B}_{126}$ & $: $ & $e_1e_2=e_3$  & $e_1e_3=e_4$ & $e_1e_4=e_5$ & $e_2e_1=e_5$ \\ & & $e_3e_2=e_4$ & $e_3e_3=e_5$ & $e_4e_2=e_5$ \\

${\rm B}_{127}^{\alpha}$ & $: $ & $e_1e_2=e_3$  & $e_1e_3=e_4$ & $e_1e_4=e_5$ & $e_2e_1=\alpha e_5$ \\ && $e_3e_2=e_4+e_5$ & $e_3e_3=e_5$ & $e_4e_2=e_5$ \\

${\rm B}_{128}^{\alpha}$ & $: $ & $e_1e_2=e_3$  & $e_1e_3=e_4$ & $e_1e_4=e_5$ \\ && $e_2e_2=e_4$  & $e_2e_3=e_5$ & $e_3e_2=\alpha e_5$ \\

${\rm B}_{129}^{\alpha}$ & $: $ & $e_1e_2=e_3$  & $e_1e_3=e_4$ & $e_1e_4=e_5$ & $e_2e_1=e_5$ \\ && $e_2e_2=e_4$ & $e_2e_3=e_5$ & $e_3e_2=\alpha e_5$ \\

${\rm B}_{130}^{\alpha}$ & $: $ & $e_1e_2=e_3$  & $e_1e_3=e_4$ & $e_1e_4=e_5$ & $e_2e_1=\alpha e_5$ \\ && $e_2e_2=e_4+e_5$ & $e_2e_3=e_5$  \\

${\rm B}_{131}$ & $: $ & $e_1e_2=e_3$  & $e_1e_3=e_4$ & $e_1e_4=e_5$ \\

${\rm B}_{132}$ & $: $ & $e_1e_2=e_3$  & $e_1e_3=e_4$ & $e_1e_4=e_5$ & $e_2e_2=e_5$ \\

${\rm B}_{133}$ & $: $ & $e_1e_2=e_3$  & $e_1e_3=e_4$ & $e_1e_4=e_5$ & $e_2e_1=e_5$\\

${\rm B}_{134}$ & $: $ & $e_1e_2=e_3$  & $e_1e_3=e_4$ & $e_1e_4=e_5$ \\ && $e_2e_1=e_5$  & $e_2e_2=e_5$\\

${\rm B}_{135}$ & $: $ & $e_1e_2=e_3$  & $e_1e_3=e_4$ & $e_1e_4=e_5$ & $e_3e_2=e_5$ \\

${\rm B}_{136}$ & $: $ & $e_1e_2=e_3$  & $e_1e_3=e_4$ & $e_1e_4=e_5$ \\ && $e_2e_1=e_5$  & $e_3e_2=e_5$ \\

${\rm B}_{137}^{\alpha}$ & $: $ & $e_1e_1=e_4$  &  $e_1e_2=e_3$  & $e_2e_1=\alpha e_5$ \\ && $e_3e_1=e_5$  & $e_3e_2=e_4+e_5$ & $e_4e_2=e_5$ \\

${\rm B}_{138}^{\alpha}$ & $: $ & $e_1e_1=e_4$  &  $e_1e_2=e_3$  & $e_1e_3=\alpha e_5$ \\ && $e_3e_1=e_5$  & $e_3e_2=e_4$ & $e_4e_2=e_5$ \\

${\rm B}_{139}^{\alpha}$ & $: $ & $e_1e_1=e_4$  &  $e_1e_2=e_3$  & $e_1e_3=\alpha e_5$ & $e_2e_1=e_5$ \\ && $e_3e_1=e_5$ & $e_3e_2=e_4$ & $e_4e_2=e_5$ \\

${\rm B}_{140}$ & $: $ & $e_1e_2=e_3$  &  $e_3e_2=e_4$  & $e_4e_2= e_5$ \\

${\rm B}_{141}$ & $: $ & $e_1e_2=e_3$  & $e_2e_1= e_5$ &  $e_3e_2=e_4$  & $e_4e_2= e_5$ \\

${\rm B}_{142}$ & $: $ & $e_1e_1=e_5$  & $e_1e_2=e_3$  &  $e_3e_2=e_4$  & $e_4e_2= e_5$ \\

${\rm B}_{143}$ & $: $ & $e_1e_1=e_5$  & $e_1e_2=e_3$ & $e_2e_1= e_5$ \\ &&  $e_3e_2=e_4$  & $e_4e_2= e_5$ \\

${\rm B}_{144}$ & $: $ & $e_1e_2=e_3$  &  $e_1e_3=e_5$ & $e_3e_2=e_4$  & $e_4e_2= e_5$ \\

${\rm B}_{145}$ & $: $ & $e_1e_2=e_3$  &  $e_1e_3=e_5$ & $e_2e_1= e_5$ \\ && $e_3e_2=e_4$  & $e_4e_2= e_5$ \\

${\rm B}_{146}^{\alpha}$ & $:$& $e_1 e_1 = e_2$  & $e_1e_2=\alpha e_4$  & $e_2 e_1= (\alpha+1)e_4$  & $e_3e_3=e_5$  \\

${\rm B}_{147}^{\alpha}$ & $:$& $e_1 e_1 = e_2$  & $e_1e_2=\alpha e_4+e_5$  & $e_2 e_1= (\alpha+1)e_4+e_5$  & $e_3e_3=e_5$  \\

${\rm B}_{148}^{\alpha}$ & $:$& $e_1 e_1 = e_2$  & $e_1e_2=\alpha e_4$  & $e_2 e_1= (\alpha+1)e_4$  \\ && $e_3e_1=e_5$  & $e_3e_3=e_5$  \\

${\rm B}_{149}^{\alpha}$ & $:$& $e_1 e_1 = e_2$  & $e_1e_2=\alpha e_4+e_5$  & $e_2 e_1= (\alpha+1)e_4+e_5$  \\ & & $e_3e_1=e_5$  & $e_3e_3=e_5$  \\

${\rm B}_{150}^{\alpha}$ & $:$& $e_1 e_1 = e_2$  & $e_1e_2=\alpha e_4$  & $e_1e_3=e_4$  \\ & & $e_2 e_1= (\alpha+1)e_4$ & $e_3e_1=e_4$ & $e_3e_3=e_5$  \\

${\rm B}_{151}^{\alpha}$ & $:$& $e_1 e_1 = e_2$  & $e_1e_2=\alpha e_4 + e_5$  & $e_1e_3=e_4$   \\ && $e_2 e_1= (\alpha+1)e_4+e_5$  & $e_3e_1=e_4$ & $e_3e_3=e_5$  \\

${\rm B}_{152}^{\alpha, \beta}$ & $:$& $e_1 e_1 = e_2$  & $e_1e_2=\alpha e_4 + \beta e_5$  & $e_1e_3=e_4$  \\ & & $e_2 e_1= (\alpha+1)e_4+\beta e_5$ & $e_3e_1=e_4+e_5$ & $e_3e_3=e_5$  \\

${\rm B}_{153}^{\alpha, \beta}$ & $:$& $e_1 e_1 = e_2$  & $e_1e_2=\alpha e_4 $  & $e_1e_3=\beta e_5$   \\ && $e_2 e_1= (\alpha+1)e_4$  & $e_3e_1=\beta e_5$ & $e_3e_3=e_5$  \\

${\rm B}_{154}^{\alpha}$ & $:$& $e_1 e_1 = e_2$  & $e_1e_2=\alpha e_4 $  & $e_1e_3=e_4+\alpha e_5$  \\ & & $e_2 e_1= (\alpha+1)e_4$ & $e_3e_1=e_4+(\alpha+1) e_5$ \\

${\rm B}_{155}^{\alpha, \beta}$ & $:$& $e_1 e_1 = e_2$  & $e_1e_2=\alpha e_4 + e_5 $  & $e_1e_3=\beta e_5$  \\ && $e_2 e_1= (\alpha+1)e_4+e_5$ & $e_3e_1=(\beta+1) e_5$ \\

${\rm B}_{156}^{\alpha, \beta}$ & $:$& $e_1 e_1 = e_2$  & $e_1e_2=\alpha e_4 $  & $e_1e_3=\beta e_5$  \\ && $e_2 e_1= (\alpha+1)e_4$  & $e_3e_1=(\beta+1) e_5$ & $e_3e_3= e_4$\\

${\rm B}_{157}^{\alpha}$ & $:$& $e_1 e_1 = e_2$  & $e_1e_2=\alpha e_4+e_5 $  & $e_2 e_1= (\alpha+1)e_4+e_5$  \\ & & $e_3e_1=e_5$   & $e_3e_3= e_4 $\\

${\rm B}_{158}$ & $:$& $e_1 e_1 = e_2$  & $e_1e_2=e_4 $  & $e_2 e_1= e_4+e_5$ & $e_3e_1=e_5$ \\

${\rm B}_{159}^{\alpha}$ & $:$& $e_1 e_1 = e_2$  & $e_1e_2=e_4 $ & $e_1e_3=e_5$\\ && $e_2 e_1= e_4+e_5$  & $e_3e_1=\alpha e_5$ \\

${\rm B}_{160}$ & $:$& $e_1 e_1 = e_2$  & $e_1e_2=e_4 $ & $e_2 e_1= e_4+e_5$ & $e_3e_3= e_5$ \\

${\rm B}_{161}$ & $:$& $e_1 e_1 = e_2$  & $e_1e_2=e_4 $ & $e_1e_3=e_4 $  \\ & & $e_2 e_1= e_4+e_5$ & $e_3 e_1= e_4$ & $e_3e_3= e_5$ \\

${\rm B}_{162}^{\alpha}$ & $:$& $e_1 e_1 = e_2$  & $e_1e_2=\alpha e_4 $ & $e_1e_3=e_5 $  \\ && $e_2 e_1= (\alpha+1)e_4$ & $e_3 e_1= e_5$ \\

${\rm B}_{163}^{\alpha}$ & $:$& $e_1 e_1 = e_2$  & $e_1e_2=\alpha e_4 $ & $e_1e_3=e_5 $ \\ & & $e_2 e_1= (\alpha+1)e_4$  & $e_3 e_1= e_5$ & $e_3 e_3= e_4$\\

${\rm B}_{164}$ & $:$& $e_1 e_1 = e_2$  & $e_1e_2= e_4 $ & $e_2 e_1= e_4$  \\ & & $e_3 e_1= e_4$   & $e_3 e_3= e_5$\\

${\rm B}_{165}^{\alpha}$ & $:$& $e_1 e_1 = e_2$  & $e_1e_2= e_4 $ & $e_1 e_3= \alpha e_5$ \\ & & $e_2 e_1= e_4$    & $e_3 e_1= (\alpha+1)e_5$ & $e_3 e_3= e_4$\\

${\rm B}_{166}$ & $:$& $e_1 e_1 = e_2$  & $e_1e_2= e_4 +e_5 $ & $e_2 e_1= e_4+e_5$ \\ & & $e_3 e_1= e_4$ & $e_3 e_3= e_5$\\

${\rm B}_{167}^{\alpha}$ & $:$& $e_1 e_1 = e_2$  & $e_1e_2= e_4$ & $e_1e_3= \alpha e_5$ \\ & & $e_2 e_1= e_4$   & $e_3 e_1= (\alpha+1)e_5$\\

${\rm B}_{168}$ & $:$& $e_1 e_1 = e_2$  & $e_1e_2= e_4$ & $e_2 e_1= e_4$ \\ & & $e_3 e_1= e_5$  & $e_3 e_3= e_5$\\

${\rm B}_{169}$ & $:$& $e_1 e_1 = e_2$  & $e_1e_2= e_4$ & $e_1 e_3= e_5$ \\ & & $e_2 e_1= e_4$  & $e_3 e_1= e_4+ e_5$ \\

${\rm B}_{170}$ & $:$& $e_1 e_1 = e_2$  & $e_1e_2= e_4$ & $e_1 e_3= e_5$  \\ && $e_2 e_1= e_4$   & $e_3 e_1= e_4+ e_5$ & $e_3 e_3= e_4$\\

${\rm B}_{171}$ & $:$& $e_1 e_2 = e_3$  & $e_1e_3=e_4$  & $e_3 e_2= e_5$  \\

${\rm B}_{172}$ & $:$& $e_1 e_1 = e_5$ & $e_1 e_2 = e_3$  & $e_1e_3=e_4$  & $e_3 e_2= e_5$  \\

${\rm B}_{173}$ & $:$&  $e_1 e_2 = e_3$  & $e_1e_3=e_4$  & $e_2 e_1 = e_5$ & $e_3 e_2= e_5$  \\

${\rm B}_{174}$ & $:$& $e_1 e_1 = e_5$ & $e_1 e_2 = e_3$  & $e_1e_3=e_4$  \\ && $e_2 e_1 = e_5$  & $e_3 e_2= e_5$  \\

${\rm B}_{175}$ & $:$& $e_1 e_2 = e_3$  & $e_1e_3=e_4$ & $e_2 e_2 = e_4$  & $e_3 e_2= e_5$  \\

${\rm B}_{176}$ & $:$& $e_1 e_1 = e_5$  & $e_1 e_2 = e_3$  & $e_1e_3=e_4$ \\ && $e_2 e_2 = e_4$   & $e_3 e_2= e_5$  \\

${\rm B}_{177}^{\alpha}$ & $:$& $e_1 e_1 = \alpha e_5$  & $e_1 e_2 = e_3$  & $e_1e_3=e_4$\\ & &  $e_2 e_1= e_5$  & $e_2 e_2 = e_4$  & $e_3 e_2= e_5$  \\

${\rm B}_{178}$ & $:$& $e_1 e_2 = e_3$  & $e_1e_3=e_4$  & $e_2e_1=e_4$ & $e_3 e_2= e_5$  \\

${\rm B}_{179}$ & $:$& $e_1 e_2 = e_3$  & $e_1e_3=e_4$  & $e_2e_1=e_4+e_5$ & $e_3 e_2= e_5$  \\

${\rm B}_{180}^{\alpha}$ & $:$& $e_1 e_1 = e_5$  & $e_1 e_2 = e_3$  & $e_1e_3=e_4$ \\ & & $e_2e_1=e_4+\alpha e_5$   & $e_3 e_2= e_5$  \\

${\rm B}_{181}^{\alpha, \beta}$ & $:$& $e_1 e_1 = \alpha e_5$  & $e_1 e_2 = e_3$  & $e_1e_3=e_4$  & $e_2e_1=e_4$ \\ & & $e_2e_2=e_4 + \beta e_5$ & $e_3 e_2= e_5$  \\

${\rm B}_{182}$ & $:$& $e_1 e_2 = e_3$  & $e_1e_3=e_4$  & $e_2 e_2= e_5$  \\

${\rm B}_{183}$ & $:$& $e_1 e_1 = e_5$  & $e_1 e_2 = e_3$  & $e_1e_3=e_4$  & $e_2 e_2= e_5$  \\

${\rm B}_{184}^{\alpha}$ & $:$& $e_1 e_1 = \alpha e_5$  & $e_1 e_2 = e_3$  & $e_1e_3=e_4$\\ & & $e_2 e_1= e_5$  & $e_2 e_2= e_5$  \\

${\rm B}_{185}^{\alpha, \beta}$ & $:$& $e_1 e_1 = \alpha e_5$  & $e_1 e_2 = e_3$  & $e_1e_3=e_4$ \\ && $e_2 e_1= \beta e_5$ & $e_2 e_2= e_5$ & $e_3 e_2= e_4$  \\

${\rm B}_{186}$ & $:$& $e_1 e_2 = e_3$  & $e_1e_3=e_4$  & $e_2e_1=e_4$ & $e_2 e_2= e_5$  \\

${\rm B}_{187}$ & $:$&  $e_1 e_1= e_5$ & $e_1 e_2 = e_3$  & $e_1e_3=e_4$ \\ & & $e_2e_1=e_4$   & $e_2 e_2= e_5$  \\

${\rm B}_{188}^{\alpha}$ & $:$&  $e_1 e_1= \alpha e_5$ & $e_1 e_2 = e_3$  & $e_1e_3=e_4$ \\ & & $e_2e_1=e_4+e_5$ & $e_2 e_2= e_5$  \\

${\rm B}_{189}^{\alpha, \beta}$ & $:$&  $e_1 e_1= \alpha e_5$ & $e_1 e_2 = e_3$  & $e_1e_3=e_4$  \\ && $e_2e_1=e_4+\beta e_5$  & $e_2 e_2= e_5$ & $e_3 e_2= e_4$ \\

%${\rm B}_{190}$ & $:$& $e_1 e_2 = e_3$  & $e_1e_3=e_4$  & $e_2 e_1= e_5$  \\

%${\rm B}_{191}$ & $:$& $e_1 e_2 = e_3$  & $e_1e_3=e_4$  & $e_2 e_1= e_5$ & $e_2 e_2= e_4$ \\

${\rm B}_{190}$ & $:$ & $e_1 e_1 = e_5$ & $e_1 e_2 = e_3$  & $e_1e_3=e_4$  & $e_2 e_1= e_5$  \\

${\rm B}_{191}$ & $:$ & $e_1 e_1 = e_5$ & $e_1 e_2 = e_3$  & $e_1e_3=e_4$  \\ && $e_2 e_1= e_5$  & $e_2 e_2= e_4$\\

${\rm B}_{192}^{\alpha}$ & $:$ & $e_1 e_1 = \alpha e_5$ & $e_1 e_2 = e_3$  & $e_1e_3=e_4$  \\ && $e_2 e_1= e_5$ & $e_3 e_2= e_4$\\

${\rm B}_{193}$ & $:$ & $e_1 e_1 = e_5$ & $e_1 e_2 = e_3$  & $e_1e_3=e_4$ \\

${\rm B}_{194}$ & $:$ & $e_1 e_1 = e_5$ & $e_1 e_2 = e_3$  & $e_1e_3=e_4$ & $e_2 e_1= e_4$ \\

${\rm B}_{195}$ & $:$ & $e_1 e_1 = e_5$ & $e_1 e_2 = e_3$  & $e_1e_3=e_4$ & $e_2 e_2= e_4$ \\

${\rm B}_{196}$ & $:$ & $e_1 e_1 = e_5$ & $e_1 e_2 = e_3$  & $e_1e_3=e_4$ \\ && $e_2 e_1= e_4$  & $e_2 e_2= e_4$ \\

${\rm B}_{197}$ & $:$ & $e_1 e_1 = e_5$ & $e_1 e_2 = e_3$  & $e_1e_3=e_4$ & $e_3 e_2= e_4$ \\

${\rm B}_{198}$ & $:$ & $e_1 e_1 = e_5$ & $e_1 e_2 = e_3$  & $e_1e_3=e_4$ \\ && $e_2 e_1= e_4$  &$e_3 e_2= e_4$ \\

${\rm B}_{199}$ & $:$ & $e_1 e_1 = e_5$ & $e_1 e_2 = e_3$  &$e_3 e_2= e_4$ \\

${\rm B}_{200}$ & $:$ & $e_1 e_1 = e_5$ & $e_1 e_2 = e_3$  & $e_2 e_1= e_5$ &$e_3 e_2= e_4$  \\

${\rm B}_{201}^{\alpha}$ & $:$ & $e_1 e_1 = e_5$ & $e_1 e_2 = e_3$ & $e_2 e_1= \alpha e_5$\\ & & $e_2 e_2= e_5$   &$e_3 e_2= e_4$ \\

${\rm B}_{202}$ & $:$ & $e_1 e_1 = e_5$ & $e_1 e_2 = e_3$  & $e_2 e_1= e_4$ &$e_3 e_2= e_4$  \\

${\rm B}_{203}$ & $:$ & $e_1 e_1 = e_5$ & $e_1 e_2 = e_3$  & $e_2 e_1= e_4+e_5$ &$e_3 e_2= e_4$  \\

${\rm B}_{204}^{\alpha}$ & $:$ & $e_1 e_1 = e_5$ & $e_1 e_2 = e_3$  & $e_2 e_1= e_4+\alpha e_5$ \\ &&$e_2 e_2= e_5$   &$e_3 e_2= e_4$  \\

${\rm B}_{205}$ & $:$  & $e_1 e_2 = e_3$ & $e_2 e_1 = e_5$ &$e_3 e_2= e_4$ \\

%${\rm B}_{208}$ & $:$  & $e_1 e_2 = e_3$ & $e_2 e_1 = e_5$  & $e_2 e_2 = e_5$ &$e_3 e_2= e_4$ \\

${\rm B}_{206}$ & $:$ & $e_1 e_1 = e_4$ & $e_1 e_2 = e_3$ & $e_2 e_1 = e_5$  &$e_3 e_2= e_4$ \\

%${\rm B}_{210}$ & $:$ & $e_1 e_1 = e_4$ & $e_1 e_2 = e_3$ & $e_2 e_1 = e_5$  & $e_2 e_2 = e_5$ \\ &&$e_3 e_2= e_4$ \\

%${\rm B}_{211}$ & $:$  & $e_1 e_2 = e_3$ & $e_2 e_2 = e_5$ &$e_3 e_2= e_4$ \\

${\rm B}_{207}$ & $:$  & $e_1 e_2 = e_3$ & $e_2 e_1 = e_4$ & $e_2 e_2 = e_5$ &$e_3 e_2= e_4$ \\

%${\rm B}_{213}$ & $:$ & $e_1 e_1 = e_4$  & $e_1 e_2 = e_3$ &  $e_2 e_2 = e_5$ &$e_3 e_2= e_4$ \\

%${\rm B}_{214}$ & $:$ & $e_1 e_1 = e_4$  & $e_1 e_2 = e_3$ &  $e_2 e_1 = e_4$ & $e_2 e_2 = e_5$ \\ & &$e_3 e_2= e_4$ \\

\end{longtable}
}

Note that 
${\rm B}_{116}^{1}$ is a commutative algebra and

\begin{center}
${\rm B}_{12}^{\frac 1 4} \simeq {\rm B}_{11}^{\frac 1 4},$
${\rm B}_{16}^{\alpha, \beta} \simeq {\rm B}_{16}^{\beta, \alpha},$
${\rm B}_{35}^{\alpha} \simeq {\rm B}_{35}^{-\alpha},$ 
${\rm B}_{38}^{\alpha} \simeq {\rm B}_{38}^{-\alpha}$
${\rm B}_{41}^{\alpha} \simeq {\rm B}_{41}^{-\alpha},$

${\rm B}_{57}^{0,0,\gamma\neq0} \simeq {\rm B}_{185}^{0, \frac{1}{\gamma}},$
${\rm B}_{58}^{(0,0)} \simeq {\rm B}_{184}^0,$
${\rm B}_{59}^{0} \simeq {\rm B}_{188}^0,$
${\rm B}_{60}^{0,0} \simeq {\rm B}_{156}^{-1,-1},$
${\rm B}_{60}^{1,0} \simeq {\rm B}_{93},$
${\rm B}_{60}^{\alpha\neq 0,1;0} \simeq {\rm B}_{184}^{\frac{1}{1-\alpha}},$
${\rm B}_{61}^{0,0} \simeq {\rm B}_{152}^{0,0},$
${\rm B}_{61}^{\alpha\neq0,0} \simeq {\rm B}_{204}^{1-\alpha},$
${\rm B}_{63}^{0} \simeq {\rm B}_{201}^{1},$
${\rm B}_{64}^{0,\gamma\neq0} \simeq {\rm B}_{185}^{(\alpha^2,\alpha), \alpha\neq0},$
${\rm B}_{64}^{0,0} \simeq {\rm B}_{201}^{1},$

${\rm B}_{81}^{(0,0)} \simeq {\rm B}_{192}^{0},$
${\rm B}_{81}^{\alpha,\beta} \simeq {\rm B}_{81}^{\frac{1}{\beta},\frac{1}{\alpha}},$
${\rm B}_{83}^{\alpha} \simeq {\rm B}_{83}^{\frac{1}{\alpha}},$ ${\rm B}_{87}^{\alpha,\beta,\gamma} \simeq {\rm B}_{87}^{\frac{\gamma}{\beta-2}, \frac{1-2\beta}{2-\beta},\frac{\alpha}{\beta-2}},$

${\rm B}_{91}^{\alpha,\beta} \simeq {\rm B}_{91}^{\alpha,-\beta}$
${\rm B}_{92}^{\alpha} \simeq {\rm B}_{92}^{-\alpha},$ 
${\rm B}_{93}^{\alpha} \simeq {\rm B}_{93}^{-\alpha},$ 
${\rm B}_{100}^{-2,0} \simeq {\rm B}_{97}^{0},$
${\rm B}_{100}^{-2,\beta\neq0} \simeq {\rm B}_{99}^{-2,\beta\neq0},$

${\rm B}_{105}^{-2} \simeq {\rm B}_{104}^{-2}.$
${\rm B}_{112}^{\beta}\simeq {\rm B}_{112}^{-\beta},$
${\rm B}_{180}^{\alpha}\simeq  {\rm B}_{180}^{-\alpha},$ 
${\rm B}_{201}^{\alpha}\simeq  {\rm B}_{201}^{-\alpha},$ 
${\rm B}_{204}^{\alpha}\simeq  {\rm B}_{204}^{-\alpha}.$

\end{center}

\end{theoremA}

 %\newpage


\begin{thebibliography}{99}


 
 
  %\bibitem{akk23}
%  Abdurasulov K.,  Kaygorodov I.,     Khudoyberdiyev A.,  
% The algebraic  classification of nilpotent Novikov   algebras, pre



 \bibitem{fkk223}
Abdelwahab H.,  Barreiro E.,  Calderón A.,  Fernández Ouaridi A.,  
The classification of nilpotent Lie-Yamaguti algebras,
Linear Algebra and Its Applications, 654 (2022), 339--378

 \bibitem{fkk222}
Abdelwahab H.,  Barreiro E.,  Calderón A.,  Fernández Ouaridi A.,  
The algebraic classification and degenerations of nilpotent Poisson algebras, Journal of Algebra, 615 (2023), 243--277.

 \bibitem{fkk22}
Abdelwahab H.,  Barreiro E.,  Calderón A.,  Fernández Ouaridi A.,  
The algebraic and geometric classification of nilpotent Lie triple systems up to dimension four,  Revista de la Real Academia de Ciencias Exactas, Físicas y Naturales. Serie A. Matemáticas, 117 (2023),   1,   11.

 
 \bibitem{bcz21}
Bai Yu.,  Chen Yu., Zhang Z.,
    Gelfand-Kirillov dimension of bicommutative algebras,
    Linear and Multilinear Algebra, 70 (2022), 22,    7623--7649.

\bibitem{akk23}   Beites P.,    Fern\'andez Ouaridi A.,   Kaygorodov I., 
The algebraic and geometric classification of transposed Poisson algebras, Revista de la Real Academia de Ciencias Exactas, Físicas y Naturales. Serie A. Matemáticas,  117 (2023), 2,  Paper 55.

  \bibitem{kb}
Belov A.,   Local finite basis property and local representability of varieties of associative rings,  Izvestiya: Mathematics, 74 (2010),  1, 1--126 


 \bibitem{bdd} 
Burde D., Dekimpe K., Deschamps S., 
    Affine actions on nilpotent Lie groups, 
    Forum Mathematicum, 21 (2009),  5, 921--934.

\bibitem{DKS09} Burde D., Dekimpe K., Deschamps S.,
LR-algebras,  
New Developments in Lie Theory and Geometry, Amer. Math. Soc., Providence, RI, 
Contemporary Mathematics, 491 (2009), 125--140.


\bibitem{DKV10} Burde D., Dekimpe K., Vercammen K.,
Complete {\rm LR}-structures on solvable Lie algebras, 
Journal of Group Theory, 13 (2010),  5, 703--719.
 
\bibitem{cfk18}
Calder\'on A.,  Fern\'andez Ouaridi A., Kaygorodov I.,  
On the classification of bilinear maps with   radical of a fixed codimension,  
   Linear and Multilinear Algebra,   70  (2022), 18, 3553--3576.

 
\bibitem{cayley} Cayley A.,
 On the theory of analytical forms called trees,
Philosophical Magazine, 13 (1857), 19--30;
 Collected Math. Papers, University Press, Cambridge, 3 (1890), 242--246.

 

  
 
 

 


\bibitem{drensky2}
Drensky V.,
 Varieties of bicommutative algebras,
 Serdica Mathematical Journal, 45 (2019), 2, 167--188.


\bibitem{drensky}
Drensky V.,
Invariant theory of free bicommutative algebras,
arXiv:2210.08317  
 

\bibitem{drensky1}
Drensky V., Zhakhayev B.,
 Noetherianity and Specht problem for varieties of bicommutative algebras,
 Journal of Algebra, 499 (2018), 1, 570--582.

\bibitem{dt03}
Dzhumadildaev A., Tulenbaev K.,
 Bicommutative algebras,
 Russian Mathematical Surveys, 58 (2003), 6, 1196--1197.

\bibitem{dit11}
Dzhumadildaev A., Ismailov N., Tulenbaev K.,
 Free bicommutative algebras,
 Serdica Mathematical Journal, 37 (2011),  1, 25--44.

\bibitem{di}
Dzhumadildaev A., Ismailov N.,
 Polynomial identities of bicommutative algebras, Lie and Jordan elements,
 Communications in Algebra, 46 (2018), 12,  5241--5251.




   
\bibitem{fkkv22}
Fern\'andez Ouaridi A.,  Kaygorodov I.,  Khrypchenko M., Volkov Yu., 
    Degenerations of nilpotent algebras,
     Journal of Pure and Applied Algebra,   226 (2022),  3, 106850.

 

    
  
 


 

\bibitem{hac16}
Hegazi A., Abdelwahab H., Calderón Martín A.,
    The classification of $n$-dimensional non-Lie Malcev algebras with $(n-4)$-dimensional annihilator, 
    Linear Algebra and its Applications, 505 (2016), 32--56.
 \bibitem{ikp20}
 Ignatyev M.,  Kaygorodov I., Popov Yu., 
  The geometric classification of $2$-step nilpotent algebras   and applications,     Revista Matemática Complutense, 35 (2022), 3, 907–922. 

 
 
 
  


\bibitem{klp20} Kaygorodov I., Lopes S., P\'{a}ez-Guill\'{a}n P.,   
Non-associative central extensions of null-filiform associative algebras, 
Journal of Algebra,   560  (2020),   1190--1210.


 
 
\bibitem{kpv20}
Kaygorodov I.,  Páez-Guillán P.,  Voronin V., 
    The algebraic and geometric classification of nilpotent bicommutative algebras, 
    Algebras and Represention Theory, 23 (2020), 6, 2331--2347. 

\bibitem{kpv21}
Kaygorodov I.,  Páez-Guillán P.,  Voronin V., 
    One-generated  nilpotent bicommutative algebras, 
     Algebra Colloquium,   29   (2022), 3, 453--474.
    
\bibitem{corr}
 Kaygorodov I.,  Popov Yu.,  Pozhidaev A.,  Volkov Yu., 
 Corrigendum to `Degenerations of Zinbiel and nilpotent Leibniz algebras', Linear and Multilinear Algebra, 70 (2022), 5, 993--995.
 
 
\bibitem{krs20} Kaygorodov I.,  Rakhimov I.,  Said Husain Sh. K., 
The algebraic  classification of nilpotent  associative commutative algebras,
Journal of Algebra and its Applications, 19 (2020), 11,    2050220.
 
 
 

\bibitem{shestak}
Shestakov I., Zhang Z., 
    Automorphisms of finitely generated relatively free bicommutative algebras, 
    Journal of Pure and  Applied Algebra, 225 (2021), 8, 106636.


\bibitem{ss78}
Skjelbred T., Sund T.,
    Sur la classification des algebres de Lie nilpotentes,
    C. R. Acad. Sci. Paris Ser. A-B, 286 (1978), 5,  A241--A242.

  
 

\end{thebibliography}
\end{document}